\documentclass[12pt]{article}

\usepackage{amsmath,amssymb,amsthm,amscd}

\begin{document}


\newcommand{\non}{\nonumber}
\newcommand{\scl}{\scriptstyle}
\newcommand{\wt}{\widetilde}
\newcommand{\wh}{\widehat}
\newcommand{\ot}{\otimes}
\newcommand{\fand}{\quad\text{and}\quad}
\newcommand{\Fand}{\qquad\text{and}\qquad}
\newcommand{\ts}{\,}
\newcommand{\tss}{\hspace{1pt}}
\newcommand{\tpr}{t^{\tss\prime}}
\newcommand{\spr}{s^{\tss\prime}}
\newcommand{\degpr}{{\rm deg}^{\tss\prime}\tss}
\newcommand{\grpr}{{\rm gr}^{\tss\prime}\ts}


\newcommand{\U}{ {\rm U}}
\newcommand{\Z}{ {\rm Z}}
\newcommand{\ZY}{ {\rm ZY}}
\newcommand{\Sr}{ {\rm S}}
\newcommand{\Ar}{ {\rm A}}
\newcommand{\Br}{ {\rm B}}
\newcommand{\Cr}{ {\rm C}}
\newcommand{\Ir}{ {\rm I}}
\newcommand{\Jr}{ {\rm J}}
\newcommand{\Prm}{ {\rm P}}
\newcommand{\X}{ {\rm X}}
\newcommand{\Y}{ {\rm Y}}
\newcommand{\ZX}{ {\rm ZX}}
\newcommand{\SY}{ {\rm SY}}
\newcommand{\Pf}{ {\rm Pf}\ts}
\newcommand{\Hf}{ {\rm Hf}\ts}
\newcommand{\tr}{ {\rm tr}\ts}
\newcommand{\End}{{\rm{End}\ts}}
\newcommand{\Hom}{{\rm{Hom}}}
\newcommand{\sgn}{ {\rm sgn}\ts}
\newcommand{\sign}{ {\rm sign}\ts}
\newcommand{\qdet}{ {\rm qdet}\ts}
\newcommand{\sdet}{ {\rm sdet}\ts}
\newcommand{\Ber}{ {\rm Ber}\ts}
\newcommand{\inv}{ {\rm inv}\ts}
\newcommand{\inva}{ {\rm inv}}
\newcommand{\middle}{ {\rm mid} }
\newcommand{\ev}{ {\rm ev}}
\newcommand{\gr}{ {\rm gr}\ts}
\newcommand{\pa}{\partial}


\newcommand{\CC}{\mathbb{C}}
\newcommand{\ZZ}{\mathbb{Z}}


\newcommand{\Ac}{{\mathcal A}}
\newcommand{\Bc}{{\mathcal B}}
\newcommand{\Oc}{{\mathcal O}}
\newcommand{\Pc}{{\mathcal P}}
\newcommand{\Qc}{{\mathcal Q}}
\newcommand{\Jc}{{\mathcal J}}


\newcommand{\Sym}{\mathfrak S}
\newcommand{\h}{\mathfrak h}
\newcommand{\gl}{\mathfrak{gl}}
\newcommand{\oa}{\mathfrak{o}}
\newcommand{\spa}{\mathfrak{sp}}
\newcommand{\osp}{\mathfrak{osp}}
\newcommand{\g}{\mathfrak{g}^{}}
\newcommand{\ggot}{\mathfrak{g}}
\newcommand{\agot}{\mathfrak{a}}
\newcommand{\sll}{\mathfrak{sl}}


\newcommand{\al}{\alpha}
\newcommand{\be}{\beta}
\newcommand{\ga}{\gamma}
\newcommand{\de}{\delta^{}}
\newcommand{\ve}{\varepsilon}
\newcommand{\vp}{\varphi}
\newcommand{\la}{\lambda}
\newcommand{\La}{\Lambda}
\newcommand{\Ga}{\Gamma}
\newcommand{\si}{\sigma}
\newcommand{\ze}{\zeta}
\newcommand{\om}{\omega^{}}


\newcommand{\rb}{\mathbf r}
\newcommand{\sbf}{\mathbf s}

\renewcommand{\theequation}{\arabic{section}.\arabic{equation}}

\newtheorem{thm}{Theorem}[section]
\newtheorem{lem}[thm]{Lemma}
\newtheorem{prop}[thm]{Proposition}
\newtheorem{cor}[thm]{Corollary}
\newtheorem{conj}[thm]{Conjecture}

\theoremstyle{definition}
\newtheorem{defin}[thm]{Definition}
\newtheorem{example}[thm]{Example}

\theoremstyle{remark}
\newtheorem{remark}[thm]{Remark}

\newcommand{\bth}{\begin{thm}}
\renewcommand{\eth}{\end{thm}}
\newcommand{\bpr}{\begin{prop}}
\newcommand{\epr}{\end{prop}}
\newcommand{\ble}{\begin{lem}}
\newcommand{\ele}{\end{lem}}
\newcommand{\bco}{\begin{cor}}
\newcommand{\eco}{\end{cor}}
\newcommand{\bcon}{\begin{conj}}
\newcommand{\econ}{\end{conj}}
\newcommand{\bde}{\begin{defin}}
\newcommand{\ede}{\end{defin}}
\newcommand{\bex}{\begin{example}}
\newcommand{\eex}{\end{example}}
\newcommand{\bre}{\begin{remark}}
\newcommand{\ere}{\end{remark}}

\newcommand{\bal}{\begin{aligned}}
\newcommand{\eal}{\end{aligned}}
\newcommand{\beq}{\begin{equation}}
\newcommand{\eeq}{\end{equation}}
\newcommand{\ben}{\begin{equation*}}
\newcommand{\een}{\end{equation*}}

\newcommand{\bpf}{\begin{proof}}
\newcommand{\epf}{\end{proof}}

\def\beql#1{\begin{equation}\label{#1}}

\title{\Large\bf  On the $R$-matrix realization of Yangians
and their representations}

\author{{D. Arnaudon, \quad A. Molev\quad and\quad E. Ragoucy}}

\date{} 

\maketitle

\vspace{17 mm}

\begin{abstract}
We study the Yangians $\Y(\agot)$ associated with the simple
Lie algebras $\agot$ of type $B$, $C$ or $D$.
The algebra $\Y(\agot)$ can be regarded as a quotient of
the extended Yangian $\X(\agot)$ whose defining relations are
written in an $R$-matrix form. In this paper we are concerned with
the algebraic structure and representations of the algebra $\X(\agot)$.
We prove an analog of the Poincar\'e--Birkhoff--Witt
theorem for $\X(\agot)$ and show that the Yangian $\Y(\agot)$
can be realized as a subalgebra of $\X(\agot)$. Furthermore, we
give an independent proof of the classification theorem for
the finite-dimensional irreducible representations of $\X(\agot)$
which implies the corresponding theorem of Drinfeld for the Yangians
$\Y(\agot)$. We also give explicit constructions for all
fundamental representation of the Yangians.
\end{abstract}



\vspace{15 mm}

\noindent
A.M.\newline
{\it School of Mathematics and Statistics\newline
University of Sydney,
NSW 2006, Australia\newline}
{\tt alexm@maths.usyd.edu.au}

\vspace{7 mm}

\noindent
E.R.\newline
{\it LAPTH, Chemin de Bellevue, BP 110\newline
F-74941 Annecy-le-Vieux cedex, France\newline}
{\tt ragoucy@lapp.in2p3.fr}

\newpage

\section{Introduction}\label{sec:int}
\setcounter{equation}{0}

For any simple Lie algebra $\agot$ over $\CC$ the corresponding
Yangian $\Y(\agot)$ is a canonical deformation of the universal enveloping
algebra $\U\big(\agot[x]\big)$, $\agot[x]=\agot\ot\CC[x]$
in the class of Hopf algebras;
see Drinfeld~\cite{d:ha, d:qg, d:nr}. In accordance to Drinfeld, each
Yangian $\Y(\agot)$ has at least three different presentations; see also
Chari and Pressley~\cite[Chapter~12]{cp:gq}.
In this paper we are concerned with the one commonly
known as the $RTT$-presentation and which preceded the other two
historically. It goes back to the work of the St.-Petersburg school on the
inverse scattering method; see e.g. Takhtajan and Faddeev~\cite{tf:qi},
Kulish and Sklyanin~\cite{ks:qs}, Tarasov~\cite{t:sq, t:im},
Reshetikhin, Takhtajan and Faddeev~\cite{rtf:ql}. In the case of
$A$ type, i.e., $\agot=\sll_N$, the $RTT$-presentation of the corresponding
Yangian turns out to be particularly useful in the applications
of the $R$-matrix techniques to the classical Lie algebras;
see e.g. the review paper~\cite{m:ya}
and references therein.
Moreover, this presentation is most convenient for the study
of various subalgebras of the $A$ type Yangian which play an important role
in the applications to the quantum spin chain models; see
e.g. Arnaudon {\it et al.\/}~\cite{aacdfr:gb, acdfr:ab, acdfr:abo},
Molev and Ragoucy~\cite{mr:rr}.

In a recent paper by Arnaudon {\it et al.\/}~\cite{aacfr:rp},
the $RTT$-presentation of the Yangian associated with the
$B$, $C$ or $D$ type Lie algebra $\agot$ was studied.
The Yangian $\Y(\agot)$ was presented as a quotient of
a quadratic algebra whose defining relations are written
in the form of an $RTT$-relation. Below we denote this algebra
by $\X(\agot)$ and call it the {\it extended Yangian\/}.
The paper \cite{aacfr:rp} contains an explicit construction
of a formal series $z(u)$ whose coefficients belong
to the center of $\X(\agot)$. As shown in \cite{aacfr:rp},
the quotient of $\X(\agot)$ by the relations $z(u)=1$ is isomorphic
to $\Y(\agot)$. In the orthogonal case $\agot=\oa_N$ ($B$ and $D$ types)
this reproduces an earlier result of Drinfeld~\cite{d:nr}.

Our aim in this paper is to describe the algebraic structure
of the extended Yangian $\X(\agot)$
for each orthogonal and symplectic Lie algebra
$\agot=\oa_N$ and $\agot=\spa_{2n}$
and classify its finite-dimensional
irreducible representations. First, we prove an analog of
the Poincar\'e--Birkhoff--Witt theorem for the algebra $\X(\agot)$.
Then, following the approach of Molev, Nazarov and Olshanski~\cite{mno:yc},
we define the Yangian $\Y(\agot)$ as a subalgebra of $\X(\agot)$.
In \cite{mno:yc}, the $A$ type Yangian $\Y(\sll_N)$ is defined as a subalgebra
of the Yangian $\Y(\gl_N)$ for the general lineal Lie algebra $\gl_N$
so that the algebra $\X(\agot)$ can be regarded as an analog of $\Y(\gl_N)$
for the $B$, $C$ and $D$ types. Furthermore, we show that
the coefficients of the series $z(u)$ are algebraically independent
and generate the center of $\X(\agot)$. This implies that
the finite-dimensional irreducible representations
of the algebras $\X(\agot)$ and $\Y(\agot)$ are essentially the same.
These representations of the Yangian $\Y(\agot)$ were classified
by Drinfeld~\cite{d:nr}; see also Chari and
Pressley~\cite[Chapter~12]{cp:gq}. However, this classification
is given in terms of a different presentation ({\it new realization\/})
of $\Y(\agot)$. At present, no explicit isomorphism
between the new realization of the orthogonal or symplectic
Yangian $\Y(\agot)$ and its $RTT$-presentation is known.
(A detailed construction of such an isomorphism
in the case of $\Y(\sll_N)$ is recently given by
Brundan and Kleshchev~\cite{bk:pp}.)
Therefore, the classification results of \cite{d:nr} do not imply
an immediate description of the finite-dimensional irreducible
representations of the extended Yangian $\X(\agot)$.

We develop an independent approach to the representation theory
for the algebras $\X(\agot)$.
We define Verma modules $M(\la(u))$ over $\X(\agot)$ in a standard way,
where $\la(u)$ is a tuple of formal series which we call the highest weight.
We show that every
finite-dimensional irreducible representation of $\X(\agot)$
is isomorphic to the unique
irreducible quotient $L(\la(u))$ of $M(\la(u))$.
We classify the finite-dimensional
irreducible representations of $\X(\agot)$ by producing
necessary and sufficient conditions on the highest weight $\la(u)$
for the module $L(\la(u))$
to be finite-dimensional; see Theorem~\ref{thm:fdim}.
Reformulating these conditions for representations of the subalgebra
$\Y(\agot)$ of $\X(\agot)$ we thus obtain another proof of Drinfeld's
theorem~\cite{d:nr} for the case of the classical Lie algebras $\agot=\oa_N$
and $\spa_{2n}$.

As a first step, we consider
the low-rank cases and construct explicit isomorphisms
$\Y(\spa_2)\simeq\Y(\sll_2)$, $\Y(\oa_3)\simeq\Y(\sll_2)$ and
$\Y(\oa_4)\simeq\Y(\sll_2)\ot \Y(\sll_2)$. The former is quite
immediate while the remaining two require appropriate versions
of the fusion procedure for $R$-matrices. The representations are then
described by using the known results for the Yangian $\Y(\sll_2)$
which are due to Tarasov~\cite{t:sq, t:im}.
For the sake of completeness, we reproduce a proof
of those results which is a simpler version of the one
contained in~\cite{m:fd}.
Using the above isomorphisms, we also give explicit formulas for the
evaluation homomorphisms from $\X(\agot)$ to the universal enveloping
algebra $\U(\agot)$ for each $\agot=\spa_2$, $\oa_3$ and $\oa_4$.

In order to establish the necessary conditions for $L(\la(u))$
to be finite-dimensional,
we use an induction argument which allows us to get the conditions
for the rank $n$ Lie algebra $\agot$ from those of rank $n-1$.
The sufficient conditions on $\la(u)$ are established by producing
finite-dimensional modules having $\la(u)$ as a highest weight.
We do this first for the so-called fundamental modules and then
employ the Hopf algebra structure on $\X(\agot)$.
In particular, this proves that every finite-dimensional
irreducible representations of $\X(\agot)$ is isomorphic
to a subquotient of a tensor product of the
corresponding fundamental modules. We also give an explicit construction
of all fundamental modules of $\X(\agot)$ basically following the approach
of Chari and Pressley~\cite{cp:fr} but avoiding the use of their
results on the singularities of the $R$-matrices.
For the applications of the
fundamental Yangian modules to the affine Toda field theories
see Chari and Pressley~\cite{cp:yi}.
\medskip

The financial support of the Australian
Research Council and the
{\it Laboratoire d'Annecy-le-Vieux de Physique Th\'eorique\/}
is acknowledged.

\section{Definitions and preliminaries}\label{sec:def}
\setcounter{equation}{0}

We let $\agot$ denote the simple complex Lie algebra of type
$B_n$, $C_n,$ or $D_n$. That is,
\beql{oos}
\agot=\oa_{2n+1},\quad \spa_{2n},\quad\text{or}\quad\oa_{2n},
\end{equation}
respectively.
Whenever possible, we consider the three cases \eqref{oos} simultaneously,
unless otherwise stated. The Lie algebra $\agot$ can be regarded as
a subalgebra of the general linear Lie algebra $\gl_N$, where
$N=2n+1$ or $N=2n$, respectively. It will be convenient
to enumerate the rows and columns of $N\times N$ matrices
by the indices $-n,\dots,-1,1,\dots,n$, if $N=2n$, and
by the indices $-n,\dots,-1,0,1,\dots,n$, if $N=2n+1$.
For $-n\leqslant i,j\leqslant n$ set
\beql{fij}
F_{ij}=E_{ij}-\theta_{ij}\ts E_{-j,-i}
\end{equation}
where the $E_{ij}$ are the elements of the standard basis of $\gl_N$
and
\beql{thetaij}
\theta_{ij}=\begin{cases}
1\qquad&\text{in the orthogonal case},\\
\sgn i\cdot\sgn j\qquad&\text{in the symplectic case}.
\end{cases}
\eeq
The elements $F_{ij}$ span the
Lie algebra $\agot$ and satisfy the relations
\beql{Fsym}
F_{ij}+\theta_{ij}\ts F_{-j,-i}=0
\eeq
for any $-n\leqslant i,j\leqslant n$,
and
\beql{comrel}
[F_{ij},F_{kl}]=
\de_{kj}\ts F_{il}-\de_{il}\ts F_{kj}
-\de_{k,-i}\ts\theta_{ij}\ts F_{-j,l}
+\de_{l,-j}\ts\theta_{ij}\ts F_{k,-i}.
\eeq

For any $n$-tuple of complex
numbers $\mu=(\mu_1,\dots,\mu_n)$
we shall denote by
$V(\mu)$ the irreducible representation of
the Lie algebra $\agot$ with the highest weight $\mu$.
That is, $V(\mu)$ is generated by a nonzero vector $\xi$
such that
\begin{alignat}{2}
F_{ij}\ts\xi&=0 \qquad &&\text{for} \quad
-n\leqslant i<j\leqslant n, \qquad \text{and}
\non\\
F_{ii}\ts\xi&=\mu_i\ts\xi \qquad &&\text{for} \quad 1\leqslant i\leqslant n.
\non
\end{alignat}
The representation $V(\mu)$ is finite-dimensional if and only if
\ben
\mu_{i}-\mu_{i+1}\in\ZZ_+\qquad\text{for}\quad i=1,\dots,n-1
\een
and
\begin{align}
-\mu_1-\mu_2\in\ZZ_+ \qquad&\text{if}\quad \agot=\oa_{2n},
\non\\
-\mu_1\in\ZZ_+ \qquad&\text{if}\quad \agot=\spa_{2n},
\non\\
-2\ts\mu_1\in\ZZ_+ \qquad&\text{if}\quad \agot=\oa_{2n+1}.
\non
\end{align}

Consider the endomorphism algebra $\End\CC^N$ and let $e_{ij}\in\End\CC^N$
be the standard matrix units (we use lower case letters to distinguish
the elements of $\End\CC^N$ from the basis elements of $\gl_N$; the latter
will also be regarded as generators of the universal enveloping
algebra $\U(\gl_N)$). We denote by $F$ the $N\times N$ matrix
whose $ij$-th entry is $F_{ij}$. We shall also regard $F$ as the
element
\beql{F}
F=\sum_{i,j=-n}^n e_{ij}\ot F_{ij}\in \End\CC^N\ot \U(\agot).
\eeq

We shall use the transposition
$t:\End\CC^N\to\End\CC^N$ which is a linear map defined
on the basis elements by the rule
\beql{tr}
(e_{ij})^t=\theta_{ij}\ts e_{-j,-i},
\eeq
and the standard transposition defined by
\beql{sttr}
(e_{ij})'=e_{ji}.
\eeq
The permutation operator $P$ is an element of $\End\CC^N\ot\End\CC^N$
given by
\beql{P}
P=\sum_{i,j=-n}^n e_{ij}\ot e_{ji}.
\eeq
We let $Q$ denote the transposed operator $Q=P^{\ts t_1}=P^{\ts t_2}$
with respect to the first or second copy of $\End\CC^N$,
\beql{Q}
Q=\sum_{i,j=-n}^n \theta_{ij}\ts e_{ij}\ot e_{-i,-j}.
\eeq
Whenever the double sign $\pm$ or $\mp$ occurs, the upper sign
corresponds to the orthogonal case while the lower sign corresponds to
the symplectic case. Note that the operators $P$ and $Q$ satisfy
the relations
\beql{pq}
P^2=1,\qquad PQ=QP=\pm Q,\qquad Q^2=N\ts Q.
\eeq
Set
\beql{kappa}
\kappa=N/2\mp 1.
\eeq
The $R$-matrix $R(u)$ is a rational function in a complex parameter $u$
with values in $\End\CC^N\ot\End\CC^N$ defined by
\beql{Ru}
R(u)=1-\frac{P}{u}+\frac{Q}{u-\kappa}.
\eeq
It is well known that $R(u)$ satisfies the Yang--Baxter equation
\beql{ybe}
R_{12}(u)\ts R_{13}(u+v)\ts R_{23}(v)=R_{23}(v)\ts R_{13}(u+v)\ts R_{12}(u),
\eeq
see \cite{ks:sy}, \cite{zz:rf}. Here both sides take values in
$\End\CC^N\ot\End\CC^N\ot\End\CC^N$ and the subscripts indicate
the copies of $\End\CC^N$ so that $R_{12}(u)=R(u)\ot 1$ etc.

Following the general approach of \cite{d:qg} and \cite{rtf:ql},
we define the {\it extended Yangian\/}
$\X(\agot)$ as an associative algebra with generators
$t_{ij}^{(r)}$, where $-n\leqslant i,j\leqslant n$ and $r=1,2,\dots$
(the zero value of $i$ and $j$ is skipped if $N=2n$),
satisfying certain quadratic relations. In order to write them down,
introduce the formal series
\beql{tiju}
t_{ij}(u)=\sum_{r=0}^{\infty}t_{ij}^{(r)}\ts u^{-r}\in\X(\agot)[[u^{-1}]],\qquad
t_{ij}^{(0)}=\de_{ij},
\eeq
and set
\beql{Tu}
T(u)=\sum_{i,j=-n}^n e_{ij}\ot t_{ij}(u)\in \End\CC^N\ot \X(\agot)[[u^{-1}]].
\eeq
Consider the algebra
$\End\CC^N\ot\End\CC^N\ot \X(\agot)[[u^{-1}]]$
and introduce its elements $T_1(u)$ and $T_2(u)$ by
\beql{T1T2}
T_1(u)=\sum_{i,j=-n}^n e_{ij}\ot 1\ot t_{ij}(u),\qquad
T_2(u)=\sum_{i,j=-n}^n 1\ot e_{ij}\ot t_{ij}(u).
\eeq
The defining relations for the algebra $\X(\agot)$ have the form
of an $RTT$-{\it relation\/}:
\beql{RTT}
R(u-v)\ts T_1(u)\ts T_2(v)=T_2(v)\ts T_1(u)\ts R(u-v).
\eeq
Equivalently, in terms of the series \eqref{tiju} they
can be written as
\beql{defrel}
\bal[]
[\tss t_{ij}(u),t_{kl}(v)]&=\frac{1}{u-v}
\Big(t_{kj}(u)\ts t_{il}(v)-t_{kj}(v)\ts t_{il}(u)\Big)\\
{}&-\frac{1}{u-v-\kappa}
\Big(\de_{k,-i}\sum_{p=-n}^n\theta_{ip}\ts t_{pj}(u)\ts t_{-p,l}(v)-
\de_{l,-j}\sum_{p=-n}^n\theta_{jp}\ts t_{k,-p}(v)\ts t_{ip}(u)\Big).
\eal
\eeq

\bre\label{rem:oonetwo}
The above definition of $\X(\agot)$ can be extended to the
cases $\agot=\oa_1$ and $\oa_2$. However, both
algebras $\X(\oa_1)$ and $\X(\oa_2)$ are commutative. In addition,
in $\X(\oa_2)$ we have $t_{-1,1}(u)=t_{1,-1}(u)=0$.
In what follows, we only deal with the orthogonal
Lie algebras $\oa_N$ for $N\geqslant 3$.
\qed
\ere

Consider an arbitrary formal series $f(u)$ of the form
\beql{fu}
f(u)=1+f_1u^{-1}+f_2 u^{-2}+\cdots\in \CC[[u^{-1}]].
\eeq
Also, let $a\in\CC$ be a constant and let $B$ be a matrix with
entries in $\CC$ such that $BB^t=1$.
It is easily derived from the defining relations
for the algebra $\X(\agot)$ that
each of the mappings
\begin{align}
\label{mult}
\mu_f:T(u)&\mapsto f(u)\ts T(u),\\
\label{shift}
\tau_a:T(u)&\mapsto T(u-a),\\
T(u)&\mapsto B\ts T(u)\ts B^{t}
\non
\end{align}
defines an automorphism of $\X(\agot)$. Furthermore,
each of the mappings
\begin{align}
T(u)&\mapsto T(-u),
\non\\
T(u)&\mapsto T^{\tss t}(u),
\non\\
T(u)&\mapsto T^{-1}(u),
\non
\end{align}
defines an anti-automorphism of $\X(\agot)$; cf. \cite[Section~1]{mno:yc}.
This is easily
verified with the use of the following property of the $R$-matrix
implied by \eqref{pq}:
\beql{rinv}
R(u)\ts R(-u)=1-\frac{1}{u^2},
\eeq
and the fact that $R(u)$ is stable under the composition of the
transpositions in the first and the second copies of $\End\CC^N$.

The extended Yangian $\X(\agot)$ is a Hopf algebra with the coproduct
\beql{Delta}
\Delta: t_{ij}(u)\mapsto \sum_{a=-n}^n t_{ia}(u)\ot t_{aj}(u),
\eeq
the antipode
\ben
\Sr: T(u)\mapsto T^{-1}(u),
\een
and the counit
\ben
\epsilon: T(u)\mapsto 1,
\een
cf. \cite{rtf:ql}, \cite[Section~1]{mno:yc}.

Multiplying both sides of \eqref{RTT} by $u-v-\kappa$,
taking $u=v+\kappa$ and replacing $v$ by $u$ we get
\beql{QTT}
Q\ts T_1(u+\kappa)\ts T_2(u)=T_2(u)\ts T_1(u+\kappa)\ts Q.
\eeq
Since $Q/N$ is a projection operator in $\CC^N\ot \CC^N$ with
a one-dimensional image, the expression on each side of \eqref{QTT}
must be equal to $Q$ times a series $z(u)$ with coefficients
in $\X(\agot)$. Since $Q\ts T_1(u)=Q\ts T^t_2(u)$
and $T_1(u)\ts Q=T^t_2(u)\ts Q$, we have
\beql{cu}
T^t(u+\kappa)\ts T(u)=T(u)\ts T^t(u+\kappa)=z(u)\ts 1,
\eeq
where
\beql{zucoeff}
z(u)=1+z_1u^{-1}+z_2u^{-2}+\dots,\qquad z_i\in \X(\agot).
\eeq
Taking the $kl$-th entries in \eqref{cu} we get the formulas
\beql{zt}
\sum_{i=-n}^n \theta_{ki}\ts t_{-i,-k}(u+\kappa)\ts t_{il}(u)=
\sum_{i=-n}^n \theta_{il}\ts t_{ki}(u)\ts t_{-l,-i}(u+\kappa)=
\de_{kl}\ts z(u).
\eeq
It was shown in \cite{aacfr:rp} that all the coefficients
$z_i$ are central in $\X(\agot)$, and
$z(u)$ has the property\footnote{Note that the $R$-matrix considered
in \cite{aacfr:rp} coincides with our $R(-u)$.}
\beql{Decu}
\Delta:z(u)\mapsto z(u)\ot z(u).
\eeq
By the Hopf algebra axioms, this implies that
the image of $z(u)$ under the antipode $\Sr$ is found by
\beql{scu}
\Sr:z(u)\mapsto z(u)^{-1}.
\eeq
By \eqref{cu}, we have
\ben
\Sr:T(u)\mapsto z(u)^{-1}\ts T^t(u+\kappa).
\een
Hence, since the transposition is involutive, we conclude
that the square of the antipode is the automorphism of $\X(\agot)$
given by
\beql{ssquare}
\Sr^2:T(u)\mapsto \frac{z(u)}{z(u+\kappa)}\ts T(u+2\kappa);
\eeq
cf. \cite[Section~1]{mno:yc}.

We define the {\it Yangian\/} $\Y(\agot)$ {\it associated with
the Lie algebra\/} $\agot$ as the subalgebra of $\X(\agot)$ which
consists of the elements stable under all the automorphisms
of the form \eqref{mult}.
It will follow from \cite{aacfr:rp}
and the results below that this definition is consistent with
the one given by Drinfeld~\cite{d:ha}; cf. \cite[Section~1]{mno:yc}.

\section{Poincar\'e--Birkhoff--Witt theorem and the center
         of the extended Yangian}
         \label{sec:pbw}
\setcounter{equation}{0}

Let us denote by $\ZX(\agot)$ the subalgebra of $\X(\agot)$
generated by all the coefficients $z_i$ of the series $z(u)$; see
\eqref{zucoeff}.

\bth\label{thm:tenpr}
We have the tensor product decomposition
\beql{tenpr}
\X(\agot)=\ZX(\agot)\ot \Y(\agot).
\eeq
\eth

\bpf
We follow the argument of \cite[Section~2.16]{mno:yc}. There exists
a unique series $y(u)$ of the form
\ben
y(u)=1+y_1u^{-1}+y_2u^{-2}+\cdots,\qquad y_i\in\ZX(\agot)
\een
such that
$y(u)y(u+\kappa)=z(u)$. In order to see this,
it suffices to write this relation in terms of the
coefficients,
\beql{ykzk}
z_k=2\tss y_k+A_k(y_1,\dots,y_{k-1}),\qquad k\geqslant 1,
\eeq
where $A_k$ is a quadratic polynomial in $k-1$ variables.
By \eqref{cu}, the image of the series $z(u)$
under the automorphism \eqref{mult} is $f(u)\ts f(u+\kappa)\ts z(u)$.
Hence, the automorphism \eqref{mult} takes $y(u)$ to $f(u)\ts y(u)$.
This implies that the series $\tau_{ij}(u)$ defined by
\beql{tauu}
\tau_{ij}(u)=y(u)^{-1}\ts t_{ij}(u),\qquad i,j=-n,\dots,n,
\eeq
are stable under all automorphisms \eqref{mult}.
Write
\ben
\tau_{ij}(u)=\de_{ij}+\tau_{ij}^{(1)}u^{-1}+\tau_{ij}^{(2)}u^{-2}+\cdots.
\een
So, the coefficients $\tau_{ij}^{(r)}$
of $\tau_{ij}(u)$ belong to the subalgebra $\Y(\agot)$. Now the
decomposition $\X(\agot)=\ZX(\agot)\cdot \Y(\agot)$ follows from
the relation $t_{ij}(u)=y(u)\ts\tau_{ij}(u)$.

It remains to demonstrate that the elements
$z_i$ are algebraically independent over $\Y(\agot)$.
Due to \eqref{ykzk}, it suffices to do this for the elements $y_i$.
Suppose on the contrary, that for some positive integer $n$
there exists a nonzero polynomial
$B$ in $n$ variables with the coefficients in $\Y(\agot)$
such that
\beql{thisequ}
B(y_1,\dots,y_n)=0.
\eeq
Take the minimal $n$ with this property.
The coefficients of $B$ are stable under
any automorphism \eqref{mult}.
Hence, applying
the automorphism \eqref{mult} with
$f(u)=1+a\tss u^{-n}$ and $a\in\CC$
to the equality \eqref{thisequ}
we get
\ben
B(y_1,\dots,y_n+a)=0
\een
for any $a\in\CC$.
This means that the polynomial $B$ does not depend on
its $n$-th variable, which contradicts the choice of $n$.
\epf

\bco\label{cor:quoti}
The Yangian $\Y(\agot)$ is isomorphic to the quotient of $\X(\agot)$
by the ideal generated by the elements $z_1,z_2,\dots,$ i.e.,
\ben
\Y(\agot)\cong \X(\agot)/(z(u)=1).
\een
Equivalently, $\Y(\agot)$ is generated by the elements
$\tau_{ij}^{(r)}$, where $-n\leqslant i,j\leqslant n$ and $r=1,2,\dots$
subject only to the relations
\beql{defreltau}
\bal[]
[\tss \tau_{ij}(u),\tau_{kl}(v)]&=\frac{1}{u-v}
\Big(\tau_{kj}(u)\ts \tau_{il}(v)-\tau_{kj}(v)\ts \tau_{il}(u)\Big)\\
{}&-\frac{1}{u-v-\kappa}
\Big(\de_{k,-i}\sum_{p=-n}^n\theta_{ip}\ts \tau_{pj}(u)\ts \tau_{-p,l}(v)-
\de_{l,-j}\sum_{p=-n}^n\theta_{jp}\ts \tau_{k,-p}(v)\ts \tau_{ip}(u)\Big)
\eal
\eeq
and
\beql{ztau}
\sum_{i=-n}^n\theta_{ki}\ts \tau_{-i,-k}(u+\kappa)\ts \tau_{il}(u)=\de_{kl}.
\eeq
\eco

\bpf
Let $\Ir$ be the ideal of $\X(\agot)$ introduced
in the statement of the corollary. Then Theorem~\ref{thm:tenpr}
implies that $\X(\agot)=\Ir\oplus \Y(\agot)$ proving the first
statement.

Now,
the coefficients $\tau_{ij}^{(r)}$
of the series $\tau_{ij}(u)$
with $i,j=-n,\dots,n$ generate the subalgebra $\Y(\agot)$.
Indeed,
it follows from the proof of Theorem~\ref{thm:tenpr} that any element
$x\in\X(\agot)$ can be uniquely written as
a polynomial $B$ in $y_1,y_2,\dots$ such that the coefficients
of $B$ are elements of the subalgebra of $\X(\agot)$ generated by
the elements $\tau_{ij}^{(r)}$.
On the other hand, if $x$
belongs to the subalgebra $\Y(\agot)$ then $B$
cannot depend on the elements $y_i$ because $x$ is stable
under all automorphisms \eqref{mult}.
Hence, $x$ belongs to the subalgebra of $\X(\agot)$ generated by
the $\tau_{ij}^{(r)}$.

Finally, recall that
the coefficients $y_i$ of the series $y(u)$ are central in
$\X(\agot)$. Hence, we derive from \eqref{tauu} that the relation
\eqref{defrel} will hold if the series $t_{ij}(u)$ are respectively
replaced by $\tau_{ij}(u)$ which gives \eqref{defreltau}.
Furthermore, \eqref{ztau} follows from \eqref{zt}. Conversely,
\eqref{defreltau} and \eqref{ztau} are defining relations for
$\Y(\agot)$ because they are respectively equivalent to
\eqref{defrel} and the relation $z(u)=1$.
\epf

\bpr\label{prop:yhopf}
The subalgebra $\Y(\agot)$ of\/ $\X(\agot)$ is a Hopf algebra
whose coproduct,
antipode and counit are obtained by restricting those from $\X(\agot)$.
\epr

\bpf
The relation \eqref{Decu} implies that
\beql{Deltay}
\Delta: y(u)\mapsto y(u)\ot y(u).
\eeq
Therefore the image of $\Y(\agot)$
under the coproduct on $\X(\agot)$ is
contained in $\Y(\agot)\ot\Y(\agot)$.
By \eqref{scu},
the image of $y(u)$ under the antipode $\Sr$ is $y(u)^{-1}$.
Hence,
\ben
\Sr:y(u)^{-1}\ts T(u)\mapsto y(u)\ts T^{-1}(u).
\een
Any automorphism \eqref{mult}
leaves the product $y(u)\ts T^{-1}(u)$ invariant and so
the subalgebra
$\Y(\agot)$ of $\X(\agot)$ is stable under $\Sr$.
\epf

Introduce an ascending filtration on the extended Yangian $\X(\agot)$
by setting
\beql{degt}
\deg t_{kl}^{(r)}=r-1
\eeq
for any $k,l\in\{-n,\dots,n\}$.
Denote by $\bar{t}_{kl}^{\ts(r)}$ and $\bar z_r$ the images of
the elements $t_{kl}^{(r)}$ and $z_r$, respectively,
in the $(r-1)$-th component of the associated graded algebra
$\gr\X(\agot)$. Then \eqref{zt} gives the relations
\beql{degtrel}
\bar{t}_{kl}^{\ts(r)}+\theta_{kl}\ts
\bar{t}_{-l,-k}^{\ts(r)}=\de_{kl}\ts\bar z_r.
\eeq
Furthermore, \eqref{tauu} implies that the degree of each element
$\tau_{kl}^{(r)}$ does not exceed $r-1$ and its image
$\bar\tau_{kl}^{\ts(r)}$ in the $(r-1)$-th component of $\gr\X(\agot)$
is given by
\beql{imtaur}
\bar\tau_{kl}^{\ts(r)}=\frac12\big(\bar{t}_{kl}^{\ts(r)}-\theta_{kl}\ts
\bar{t}_{-l,-k}^{\ts(r)}\big).
\eeq
The ascending filtration on the Yangian $\Y(\agot)$
is induced by the one on $\X(\agot)$. We denote by
$\gr\Y(\agot)$ the associated graded algebra.

\bpr\label{prop:hom}
The mapping
\beql{hom}
F_{ij}\ts x^{r-1}\mapsto \bar\tau_{ij}^{\ts(r)}
\eeq
defines an algebra homomorphism $\psi:\U\big(\agot[x]\big)\to\gr\Y(\agot)$.
\epr

\bpf
By \eqref{imtaur},
\ben
\bar\tau_{kl}^{(r)}+\theta_{kl}\ts \bar\tau_{-l,-k}^{(r)}=0
\een
for any $-n\leqslant k,l\leqslant n$ and $r\geqslant 1$.
Furthermore, using the expansion
\ben
\frac{1}{u-v}=u^{-1}+u^{-2}\tss v+\cdots,
\een
take the coefficients at $u^{-r}v^{-s}$ on both sides
of the relation \eqref{defreltau}.
Keeping the highest degree terms, we come to
\ben
[\bar\tau_{ij}^{(r)},\bar\tau_{kl}^{(s)}]=
\de_{kj}\ts\bar\tau_{il}^{(r+s-1)}-\de_{il}\ts\bar\tau_{kj}^{(r+s-1)}
-\de_{k,-i}\ts\theta_{ij}\ts\bar\tau_{-j,l}^{(r+s-1)}
+\de_{l,-j}\ts\theta_{ij}\ts\bar\tau_{k,-i}^{(r+s-1)}.
\een
It remains to compare these relations with
\eqref{Fsym} and \eqref{comrel}.
\epf

Since the graded algebra $\gr\Y(\agot)$ is generated
by the elements $\bar\tau_{ij}^{\ts(r)}$, the homomorphism $\psi$
defined in Proposition~\ref{prop:hom} is obviously surjective.
Our aim now is to show that $\psi$ is an algebra isomorphism
(see Theorem~\ref{thm:isom} below).
We shall follow the approach of Nazarov's paper~\cite[Section~2]{n:yq},
where a similar result was established for the Yangian of the queer Lie
superalgebra.

Let $\rho$ be the vector representation of the Lie algebra $\agot$
on the vector space $\CC^N$. So,
\ben
\rho:F_{ij}\mapsto e_{ij}-\theta_{ij}\ts e_{-j,-i}.
\een
For any $c\in\CC$
consider the corresponding
evaluation representation $\rho_c$ of the polynomial
current Lie algebra $\agot\tss[x]$ given by
\ben
\rho_c:F_{ij}\tss x^s\mapsto c^{\tss s}\rho\tss(F_{ij}),
\qquad
s\geqslant0.
\een

For any $c_1,\dots,c_l\in\CC$ consider the tensor product
of the evaluation representations of $\agot\tss[x]$,
\ben
\rho_{c_1,\dots,\ts c_l}=\rho_{c_1}\ot\cdots\ot\rho_{c_l}.
\een

\ble\label{lem:kernels}
Let the parameters $c_1,\dots,c_l$ and integer $l\geqslant0$
vary. Then the intersection in $\U(\agot\tss[x])$ of the
kernels of all representations $\rho_{c_1,\dots,\ts c_l}$ is trivial.
\ele

\bpf
Choose a basis $Y_1,\dots,Y_{M}$ of $\agot$,
where $M=\dim\agot$, and set $y_i=\rho(Y_i)$.
Let $A$ be a nonzero element of $\U(\agot[x])$.
Choose a total ordering on the set of basis elements $Y_ix^s$ of
$\agot[x]$ and
write $A$ as a linear combination of ordered monomials
in the basis elements.
Let $m$ be the maximal length of monomials which occur in $A$.
For each monomial
\beql{monom}
(Y_{a_1}x^{s_1})\cdots(Y_{a_m}x^{s_m})\in\U(\agot[x])
\eeq
occurring in $A$ consider
the corresponding symmetrized elements
\beql{monomsym}
\sum_{q\in\Sym_m}\,
(\tss Y_{a_{q(1)}}x^{s_{q(1)}})
\ot\dots\ot
(\tss Y_{a_{q(m)}}x^{s_{q(m)}})
\in(\tss\agot[x]\tss)^{\ot m}.
\eeq
Regarding $\U(\agot[x])$ as the quotient of the tensor algebra
of $\agot[x]$ we derive that the elements \eqref{monomsym}
are linearly independent.
Identifying the vector spaces
\ben
(\tss\agot[x]\tss)^{\ot m}=\agot^{\tss\ot m}[x_1,\dots,x_m],
\een
we can regard the sum \eqref{monomsym}
as a polynomial function in $m$ independent variables
$x_1,\dots,x_m$ with values in the
vector space $\agot^{\tss\ot m}$,
\beql{monomsympol}
\sum_{q\in\Sym_m} x_{1}^{s_{q(1)}}\cdots x_{m}^{s_{q(m)}}
\ts Y_{\ts a_{q(1)}}\ot\cdots\ot Y_{\ts a_{q(m)}}.
\eeq

Note that
\ben
\rho_{c_1,\dots, c_l}: Y_a x^s\mapsto \sum_{k=1}^l c_k^s\ts y_a^{[k]},
\qquad
y_a^{[k]}=1^{\ot\ts (k-1)}\ot y_a\ot 1^{\ot\ts (l-k)}.
\een
Hence, the image of the monomial \eqref{monom}
under the representation $\rho_{c_1,\dots,\ts c_m}$
is given by
\beql{imagerho}
\sum_{k_1,\dots,k_m=1}^m c_{k_1}^{s_1}\cdots c_{k_m}^{s_m}
\ts y_{a_1}^{[k_1]}\cdots y_{a_m}^{[k_m]}\in \End(\CC^N)^{\ot m}.
\eeq
Let us complete the set of matrices $y_1,\dots,y_M$ to a basis
$y_1,\dots,y_{N^2}$ of $\End\CC^N$ in such a way that
the identity matrix $1\in\End\CC^N$ occurs as a basis vector $y_i$
for some $i\in\{M+1,\dots,N^2\}$. Denote by
$V_m$ the subspace in $(\End\CC^N)^{\ot m}$
spanned by the basis elements
$y_{\tss i_1}\ot\cdots\ot y_{\tss i_m}$ where
at least one of the tensor factors
is $1$. Observe that the image under the representation
$\rho_{c_1,\dots,\ts c_m}$
of any monomial of length $<m$
occurring in $A$ is contained in $V_m$.
Furthermore, modulo elements belonging to $V_m$, the sum \eqref{imagerho}
can be written as
\beql{evalc}
\sum_{q\in\Sym_m} c_{1}^{s_{q(1)}}\cdots c_{m}^{s_{q(m)}}
\ts y_{\ts a_{q(1)}}\ot\cdots\ot y_{\ts a_{q(m)}}.
\eeq
This sum is the value of \eqref{monomsympol}
under the specialization $x_i=c_i$ and
replacement of $Y_i$ with $y_i=\rho(Y_i)$ for all $i=1,\dots,m$.
However, since $\rho$ is faithful and
the elements \eqref{monomsym}
are linearly independent, there exist values of the parameters
$c_1,\dots,c_m$ such that the corresponding sums \eqref{evalc}
are linearly independent modulo the subspace $V_m$
which completes the proof.
\epf

We are now in a position to prove the following.

\bth\label{thm:isom}
The mapping $\psi:\U\big(\agot[x]\big)\to\gr\Y(\agot)$
defined in \eqref{hom} is an algebra isomorphism.
\eth

\bpf
Due to Proposition~\ref{prop:hom}, we only need to show that
the kernel of $\psi$ is trivial.
Let $C$ be a nonzero element of $\U\big(\agot[x]\big)$.
We shall show that $\psi(C)\ne 0$. The universal enveloping algebra
$\U\big(\agot[x]\big)$ has a grading defined on the generators
by declaring the degree of $F_{ij}x^s$ to be equal to $s$.
Then $\psi$ is obviously a homomorphism
of graded algebras. Hence, we may assume that $C$ is homogeneous
of degree, say, $d$. Write
\beql{C}
C=\sum {C\ts}_{i_1j_1,\dots,\ts i_mj_m}^{r_1,\dots,\ts r_m}
(F_{i_1j_1}x^{r_1-1})\cdots (F_{i_mj_m}x^{r_m-1}),
\eeq
summed over the indices $i_a,j_a,r_a$ such that $r_1+\cdots+r_m=d+m$.

Consider the element $C^{\tss\prime}\in \Y(\agot)$ given by the formula
\ben
C^{\tss\prime}=\sum {C\ts}_{i_1j_1,\dots,\ts i_mj_m}^{r_1,\dots,\ts r_m}
\tau_{i_1j_1}^{(r_1)}\cdots \tau_{i_mj_m}^{(r_m)},
\een
where the summation is taken over the same set of indices as in \eqref{C}
with the same coefficients. Then the image of $C^{\tss\prime}$ in the $d$-th
component of the graded algebra $\gr\Y(\agot)$ coincides with $\psi(C)$.
So, it suffices to show that $\deg C^{\tss\prime}=d$.

Applying the standard transposition \eqref{sttr}
to the third copy of $\End\CC^N$ in the
Yang--Baxter equation \eqref{ybe} and using \eqref{rinv} we come to
the relation
\beql{ybet}
R_{12}(u-v)\ts R_{13}^{\ts\prime}(-u)\ts R_{23}^{\ts\prime}(-v)=
R_{23}^{\ts\prime}(-v)\ts R_{13}^{\ts\prime}(-u)\ts R_{12}(u-v),
\eeq
where
\ben
R^{\ts\prime}(u)=1-\frac{P^{\ts\prime}}{u}+\frac{Q^{\ts\prime}}{u-\kappa}
\een
with the transposition applied to the first (or second) copy of $\End\CC^N$.
Hence, by the defining relations \eqref{RTT} of the algebra $\X(\agot)$
we conclude that the mapping $T(u)\mapsto R^{\ts\prime}(-u)$ defines a representation
of $\X(\agot)$ in the space $\CC^N$. Taking its composition with
the automorphism \eqref{shift} we obtain for any $c\in\CC$
the representation $\si_c:T(u)\mapsto R^{\ts\prime}(-u+c)$. Equivalently, in terms
of the generating series \eqref{tiju} we have
\beql{vectrep}
\si_c:t_{ij}(u)\mapsto
\de_{ij}+e_{ij}\ts (u-c)^{-1}-\theta_{ij}\ts e_{-j,-i}\ts (u+\kappa-c)^{-1}.
\eeq
Since the transpositions \eqref{tr} and
\eqref{sttr} commute,
using \eqref{cu} and the relations
\ben
(Q^{\ts\prime})^2=1,\qquad
P^{\ts\prime}Q^{\ts\prime}=Q^{\ts\prime}P^{\ts\prime}=\pm P^{\ts\prime},
\qquad (P^{\ts\prime})^2=N\ts P^{\ts\prime},
\een
we derive that the image of $z(u)$ under $\si_c$ is given by
\ben
\si_c:z(u)\mapsto 1-\frac{1}{(u-c+\kappa)^2}.
\een
There exists a unique series $f_c(u)\in 1+u^{-1}\CC[[u^{-1}]]$ such that
\ben
f_c(u)\ts f_c(u+\kappa)=\frac{(u-c+\kappa)^2}{(u-c+\kappa)^2-1}.
\een
Then
$
\si_c:y(u)^{-1}\mapsto f_c(u)
$
so that due to
\eqref{vectrep}, for the image of the series $\tau_{ij}(u)$ under
$\si_c$ we have
\beql{vectreptau}
\si_c:\tau_{ij}(u)\mapsto f_c(u)\Big(
\de_{ij}+e_{ij}\ts (u-c)^{-1}-\theta_{ij}\ts e_{-j,-i}\ts (u+\kappa-c)^{-1}
\Big).
\eeq
Observe that the coefficient of the series $f_c(u)$ at $u^{-k}$ is
a polynomial in $c$ of degree $\leqslant k-1$. Therefore,
taking the coefficient at $u^{-r}$ in \eqref{vectreptau} we find that
the image of $\tau_{ij}^{(r)}$ under $\si_c$ is a polynomial
in $c$ of degree $\leqslant r-1$ with coefficients in $\End\CC^N$. Moreover,
the coefficient of this polynomial at $c^{r-1}$ coincides
with $\rho(F_{ij})$.

Using Proposition~\ref{prop:yhopf}, we can construct a representation
of $\Y(\agot)$ in the space $(\CC^N)^{\ot\ts l}$ by
\ben
\si_{c_1,\dots,\ts c_l}=\si_{c_1}\ot\cdots\ot\si_{c_l},\qquad c_i\in\CC.
\een
The image of the element $C^{\tss\prime}$ under $\si_{c_1,\dots,\ts c_l}$
is a polynomial in $c_1,\dots,c_l$ of degree $\leqslant d$.
Moreover, the homogeneous component of degree $d$ of this
polynomial coincides with
$D=\rho_{c_1,\dots,\ts c_l}(C)$.
By Lemma~\ref{lem:kernels}, there exist values of the
parameters $c_1,\dots,c_l$
such that $D\ne 0$.
This implies that the element $C^{\tss\prime}$ has degree $d$ and so
$\psi(C)\ne 0$.
\epf

The following is an analog of the Poincar\'e--Birkhoff--Witt
theorem for the algebra $\Y(\agot)$.
It is immediate from Theorem~\ref{thm:isom}.

\bco\label{cor:pbwy}
Given any total ordering on the set of generators
$
\tau_{ij}^{(r)}
$
with
\ben
i+j>0,\qquad r\geqslant 1,\qquad\text{in the orthogonal case},
\een
and
\ben
i+j\geqslant 0,\qquad r\geqslant 1,\qquad\text{in the symplectic case},
\een
the ordered monomials in the generators form a basis of $\Y(\agot)$.
\qed
\eco

\bre\label{rem:polyn}
The algebra $\Y(\agot)$ admits another filtration defined by
setting the degree of the generator $\tau_{ij}^{(r)}$
to be equal to $r$. It follows from
Corollary~\ref{cor:quoti} that the associated graded
algebra $\wt{\text{gr}}\ts \Y(\agot)$ is commutative.
Let $\wt{\tau}_{ij}^{\ts(r)}$ denote the image of $\tau_{ij}^{(r)}$
in the $r$-th component of $\wt{\text{gr}}\ts \Y(\agot)$.
By Corollary~\ref{cor:pbwy},
the graded algebra
$\wt{\rm{gr}}\ts \Y(\agot)$
is isomorphic to the algebra of polynomials in the variables
$\wt{\tau}_{ij}^{\ts(r)}$, where the indices $i,j,r$ are subject to
the same conditions as in Corollary~\ref{cor:pbwy}.
\qed
\ere

Recall that $\ZX(\agot)$ is the subalgebra of $\X(\agot)$
generated by the coefficients $z_i$ of the series $z(u)$.

\bco\label{cor:centerx}
\begin{itemize}
\item[{\rm(i)}]
The center of the algebra $\Y(\agot)$ is trivial.
\item[{\rm(ii)}]
The center
of the algebra $\X(\agot)$ coincides with $\ZX(\agot)$.
\item[{\rm(iii)}]
The coefficients $z_1,z_2,\dots$ of the series $z(u)$ are
algebraically independent over $\CC$, so that the subalgebra
$\ZX(\agot)$ of $\X(\agot)$ is isomorphic to the algebra
of polynomials in countably many variables.
\end{itemize}
\eco

\bpf It is well known that the center of
the universal enveloping algebra $\U\big(\agot[x]\big)$
is trivial; see e.g. \cite[Proposition~2.12]{mno:yc}.
So {\rm(i)} and {\rm(ii)} follow from Theorem~\ref{thm:isom}.
It is implied by the proof of Theorem~\ref{thm:isom} that
the elements $y_1,y_2,\dots$ of the series $y(u)$
are algebraically independent over $\CC\subset\Y(\agot)$.
Hence so are the elements $z_i$, $i\geqslant 1$.
\epf

We shall also need the following version of the Poincar\'e--Birkhoff--Witt
theorem for the algebra $\X(\agot)$.

\bco\label{cor:pbwx}
Given any total ordering on the set of elements
$
t_{ij}^{(r)}
$
and $z_r$ with
\ben
i+j>0,\qquad r\geqslant 1,\qquad\text{in the orthogonal case},
\een
and
\ben
i+j\geqslant 0,\qquad r\geqslant 1,\qquad\text{in the symplectic case},
\een
the ordered monomials in these elements form a basis of $\X(\agot)$.
\eco

\bpf
By Theorems~\ref{thm:tenpr}, \ref{thm:isom}
and Corollary~\ref{cor:centerx}(iii),
the graded algebra
$\gr\X(\agot)$
is isomorphic to the tensor product of
the universal enveloping algebra $\U(\agot[x])$ and the algebra
of polynomials $\CC[\ze_1,\ze_2,\dots]$ in indeterminates
$\ze_r$. An isomorphism is given by
\ben
\bar{t}_{ij}^{\ts(r)}\mapsto F_{ij}\ts x^{r-1}+\frac12\ts\de_{ij}\ts\ze_r,
\een
so that $\ze_r$ is the image of $\bar z_r$; see \eqref{degtrel}.
This implies the statement.
\epf

\bpr\label{prop:emb}
The assignment
\beql{emb}
F_{ij}\mapsto \tau_{ij}^{(1)}
\eeq
defines an embedding $\U(\agot)\hookrightarrow \Y(\agot)$,
while the assignment
\beql{embx}
F_{ij}\mapsto \frac12\ts \Big(t_{ij}^{(1)}-\theta_{ij}\ts t_{-j,-i}^{(1)}\Big)
\eeq
defines an embedding $\U(\agot)\hookrightarrow \X(\agot)$.
\epr

\bpf
The defining relations \eqref{defreltau}
and \eqref{ztau} of $\Y(\agot)$ imply that the map \eqref{emb}
is a homomorphism. Its injectivity follows from Corollary~\ref{cor:pbwy}.
Furthermore,
by \eqref{tauu} we have $\tau_{ij}^{(1)}=t_{ij}^{(1)}-\de_{ij}\ts y_1$.
It remains to observe that $2y_1=z_1=t_{i\tss i}^{(1)}+t_{-i,-i}^{(1)}$
for any $i$ and $t_{ij}^{(1)}=-\theta_{ij}\ts t_{-j,-i}^{(1)}$ for $i\ne j$ by
\eqref{zt}.
\epf

\section{Isomorphisms for low rank Yangians}\label{sec:isom}
\setcounter{equation}{0}

Recall that the Yangian $\Y(\gl_N)$ for the general linear Lie algebra
$\gl_N$ is defined as a unital associative algebra
with countably many generators $T_{ij}^{(1)},\ T_{ij}^{(2)},\dots$ where
$1\leqslant i,j\leqslant N$,
and the defining relations
\beql{defyang}
[T^{(r+1)}_{ij}, T^{(s)}_{kl}]-[T^{(r)}_{ij}, T^{(s+1)}_{kl}]=
T^{(r)}_{kj} T^{(s)}_{il}-T^{(s)}_{kj} T^{(r)}_{il},
\eeq
where $r,s\geqslant 0$ and  $T^{(0)}_{ij}=\delta_{ij}$.
Equivalently, these relations can be written as
\beql{defequiv}
[T^{(r)}_{ij}, T^{(s)}_{kl}] =\sum_{a=1}^{\min(r,s)}
\Big(T^{(a-1)}_{kj} T^{(r+s-a)}_{il}-T^{(r+s-a)}_{kj} T^{(a-1)}_{il}\Big).
\eeq
Introducing the generating
series,
\ben
T_{ij}(u) = \delta_{ij} + T^{(1)}_{ij} u^{-1} + T^{(2)}_{ij}u^{-2} +
\cdots\in\Y(\gl_N)[[u^{-1}]],
\een
we can also write \eqref{defyang} in the form
\beql{defreltij}
(u-v)\ts [T_{ij}(u),T_{kl}(v)]=T_{kj}(u)T_{il}(v)-T_{kj}(v)T_{il}(u).
\eeq
Equivalently, using the notation of Section~\ref{sec:def} and introducing the
matrices
\beql{Rugl}
R^{\tss\circ}(u)=1-{P}\ts{u^{-1}}
\eeq
and
\beql{Tugl}
T^{\tss\circ}(u)=\sum_{i,j=1}^N e_{ij}\ot T_{ij}(u)\in
\End\CC^N\ot \Y(\gl_N)[[u^{-1}]],
\eeq
we can present the defining relations in the form
of an $RTT$-relation
\beql{RTTgl}
R^{\tss\circ}(u-v)\ts T^{\tss\circ}_1(u)\ts T^{\tss\circ}_2(v)=
T^{\tss\circ}_2(v)\ts T^{\tss\circ}_1(u)\ts R^{\tss\circ}(u-v);
\eeq
cf. \eqref{RTT}. We use the superscript "$\scriptstyle{\circ}$" here to distinguish
the objects related to $\Y(\gl_N)$ from those related
to the algebra $\X(\agot)$.

The Yangian $\Y(\gl_N)$ is a Hopf algebra with the coproduct
\beql{Deltagl}
\Delta: T_{ij}(u)\mapsto
\sum_{k=1}^N
\,
T_{ik}(u)\ot T_{kj}(u).
\eeq

An ascending filtration on $\Y(\gl_N)$ can be defined by
setting
\beql{deggl} \deg T_{ij}^{(r)}=r-1.
\eeq
Let $\overline T_{ij}^{\ts(r)}$ denote the image of the generator
$T_{ij}^{(r)}$ in the $(r-1)$-th component of the
associated graded algebra $\gr\Y(\gl_N)$. We have an algebra
isomorphism
\beql{isomgrgl}
\U(\gl_N[x])\to \gr\Y(\gl_N),\qquad E_{ij}\ts x^{r-1}\mapsto
\overline T_{ij}^{\ts(r)}.
\eeq

The assignment
\beql{eval}
\ev:
T_{ij}(u)\mapsto \delta_{ij}+E_{ij}\ts u^{-1}
\eeq
defines a surjective homomorphism $\Y(\gl_N)\to\U(\gl_N)$. Moreover,
the assignment
$
E_{ij}\mapsto T_{ij}^{(1)}
$
defines an embedding $\U(\gl_N)\hookrightarrow \Y(\gl_N)$.

For any series
$
g(u)\in 1+u^{-1}\CC[[u^{-1}]]
$
consider the automorphism
of $\Y(\gl_N)$ defined by
\beql{multyang}
T_{ij}(u)\mapsto g(u)\ts T_{ij}(u).
\eeq
The Yangian for $\sll_N$ is the subalgebra
$\Y(\sll_N)$ of $\ts\Y(\gl_N)$
which consists of the elements stable under all automorphisms
\eqref{multyang}.

The algebra $\Y(\gl_N)$ is isomorphic to the tensor product
of its subalgebras
\beql{decompysln}
\Y(\gl_N)=\ZY(\gl_N)\ot  \Y(\sll_N),
\eeq
where $\ZY(\gl_N)$ denotes the center of the algebra $\Y(\gl_N)$.
The subalgebra $\ZY(\gl_N)$ is generated by the coefficients
of the series $D(u)$ called the quantum determinant.
In the case $N=2$ it takes the form
\ben
D(u)=T_{11}(u)\ts T_{22}(u-1)-T_{21}(u)\ts T_{12}(u-1).
\een
Define the series $d(u)$ with coefficients in $\ZY(\gl_2)$
by the relation $d(u)\ts d(u-1)=D(u)$. Then all the coefficients of the
series $\mathcal T_{ij}(u)=d(u)^{-1}\ts T_{ij}(u)$
belong to the subalgebra $\Y(\sll_2)$. The series
$\mathcal T_{ij}(u)$ satisfy the relations
\beql{dreltau}
(u-v)\ts [\mathcal T_{ij}(u),\mathcal T_{kl}(v)]=
\mathcal T_{kj}(u)\ts\mathcal T_{il}(v)-\mathcal T_{kj}(v)\ts\mathcal T_{il}(u)
\eeq
and
\beql{qdetone}
\mathcal T_{11}(u)\ts \mathcal T_{22}(u-1)-
\mathcal T_{21}(u)\ts \mathcal T_{12}(u-1)=1
\eeq
which are defining relations for the algebra $\Y(\sll_2)$.
In other words, the Yangian $\Y(\sll_2)$ is isomorphic to the quotient of
$\Y(\gl_2)$ by the ideal generated by all the coefficients of  $D(u)$.

For more details on the algebraic
structure of the Yangians $\Y(\gl_N)$ and $\Y(\sll_N)$
see e.g. \cite{mno:yc}, \cite{bk:pp}.

\subsection{Extended Yangian $\X(\spa_2)$}\label{subsec:c1}

Observe that if $N=2$ then in the symplectic case
the operators $P$ and $Q$ satisfy $P+Q=1$; see \eqref{P} and \eqref{Q}.
Therefore, for the corresponding $R$-matrix \eqref{Ru} we have
\ben
R(u)=\frac{u-1}{u-2}\ts \Big(1-\frac{2\ts P}{u}\Big)
=\frac{u-1}{u-2}\cdot R^{\tss\circ}(u/2).
\een
This implies the following isomorphism where we adopt the convention
of Section~\ref{sec:def} for numbering the rows and columns
of $2\times2$ matrices by the indices $\{-1,1\}$.

\bpr\label{prop:c1}
The mapping
\beql{isomc1}
t_{ij}(u)\mapsto T_{ij}(u/2),
\qquad i,j\in\{-1,1\}
\eeq
defines an isomorphism
$\phi:\X(\spa_2)\to\Y(\gl_2)$.
\epr

\bpf
This is immediate from the defining relations \eqref{RTT} and \eqref{RTTgl}.
\epf

\bco\label{cor:restc1}
The restriction of the isomorphism \eqref{isomc1} to the subalgebra
$\Y(\spa_2)$ of $\X(\spa_2)$ induces an isomorphism
$\Y(\spa_2)\to\Y(\sll_2)$.
\eco

\bpf
Recall that the subalgebra $\Y(\spa_2)$ consists of the elements
stable under all automorphisms of $\X(\spa_2)$ of the form
\eqref{mult}. However, given a series $f(u)$ in $u^{-1}$ with complex
coefficients, the mapping \eqref{isomc1} takes
$f(u)\ts t_{ij}(u)$ to $f(u)\ts T_{ij}(u/2)$.
So, we have the relation
$
\phi\circ\mu_f=\mu_g\circ\phi,
$
and hence $\mu_f\circ\phi^{-1}=\phi^{-1}\circ\mu_g$,
where $g(u)$ is the series in $u^{-1}$ defined by $g(u)=f(2u)$.
Thus, the image
of $\Y(\spa_2)$ under the isomorphism $\phi$
coincides with the subalgebra
$\Y(\sll_2)$ of $\Y(\gl_2)$, yielding
the desired isomorphism.
\epf

\bco\label{cor:evc1}
The mapping
\beql{evalc1}
\ev:
T(u)\mapsto 1+F\ts u^{-1}
\eeq
defines a surjective homomorphism $\X(\spa_2)\to\U(\spa_2)$.
\eco

\bpf
The composition of the evaluation homomorphism \eqref{eval} and
the natural projection $\gl_N\to\sll_N$ yields a homomorphism
$\Y(\gl_N)\to\U(\sll_N)$. For $N=2$ it takes the form
\ben
\bal
T_{-1,-1}(u)&\mapsto 1+\big(E_{-1,-1}-E_{1,1}\big)\ts(2u)^{-1},\qquad
T_{-1,1}(u)\mapsto E_{-1,1}u^{-1},\\
T_{1,1}(u)&\mapsto 1+\big(E_{1,1}-E_{-1,-1}\big)\ts(2u)^{-1},
\qquad
T_{1,-1}(u)\mapsto E_{1,-1}u^{-1}.
\eal
\een
Applying the isomorphism of Proposition~\ref{prop:c1} and using the
generators $F_{ij}$ of $\spa_2\cong\sll_2$ we get a homomorphism
$\X(\spa_2)\to\U(\spa_2)$ given by
\ben
\ev:t_{ij}(u)\mapsto \de_{ij}+F_{ij}\ts u^{-1},\qquad i,j\in\{-1,1\}.
\een
Obviously, it is surjective.
\epf

\subsection{Extended Yangian $\X(\oa_3)$}\label{subsec:b1}

We shall now use a more standard notation for the generators
of the Yangian $\Y(\gl_2)$, where the indices $i,j$ in the defining
relations \eqref{defyang} and
\eqref{defreltij} run over the set $\{1,2\}$.
Consider the
vector space $\CC^2$ with its canonical basis $e_1,e_2$ and denote by $V$
the three-dimensional subspace of $\CC^2\ot\CC^2$ spanned by the vectors
\ben
v_{-1}=e_1\ot e_1,
\qquad v_0=\frac{1}{\sqrt 2}\ts\big(e_1\ot e_2+e_2\ot e_1\big),
\qquad v_{1}=-e_2\ot e_2.
\een
We identify $V$ with $\CC^3$ regarding $v_{-1},v_{0},v_{1}$ as its
canonical basis. In particular, the operators $P_{\tss V}$ and $Q_{\tss V}$
in $V\ot V$
will be given by the respective formulas \eqref{P} and \eqref{Q}
so that, for instance, $P_{\tss V}(v_0\ot v_1)=v_1\ot v_0$.
Similarly, we regard the generator matrix $T(u)=(t_{ij}(u))$
as an element of $\End V\ot\X(\oa_3)[[u^{-1}]]$.
Note that the operator
$(1+P)/2$ is a projection of $(\CC^2)^{\ot 2}$ to the subspace
$V$. Due to \eqref{RTTgl}, we have
\ben
\frac{1+P}{2}\cdot T^{\tss\circ}_1(2u)\ts T^{\tss\circ}_2(2u+1)
=T^{\tss\circ}_2(2u+1)\ts T^{\tss\circ}_1(2u)\cdot \frac{1+P}{2},
\een
because $R^{\tss\circ}(-1)=1+P$. Therefore, we may regard
each side of this relation as an element of $\End V\ot\Y(\gl_2)[[u^{-1}]]$.

\bpr\label{prop:b1}
The mapping
\beql{hommatr}
T(u)\mapsto
\frac{1+P}{2}\cdot T^{\tss\circ}_1(2u)\ts T^{\tss\circ}_2(2u+1)
\eeq
defines an isomorphism
$\phi:\X(\oa_3)\to\Y(\gl_2)$. More explicitly, the images of the generators
under the isomorphism are given by the formulas
\ben
\bal
t_{-1,-1}(u)&\mapsto T_{11}(2u)\ts T_{11}(2u+1)\\
t_{-1,0}(u)&\mapsto \frac{1}{\sqrt 2}\ts
\Big(T_{11}(2u)\ts T_{12}(2u+1)+T_{12}(2u)\ts T_{11}(2u+1)\Big)\\
t_{-1,1}(u)&\mapsto -T_{12}(2u)\ts T_{12}(2u+1)\\
t_{0,-1}(u)&\mapsto \frac{1}{\sqrt 2}\ts
\Big(T_{11}(2u)\ts T_{21}(2u+1)+T_{21}(2u)\ts T_{11}(2u+1)\Big)\\
t_{0,0}(u)&\mapsto
T_{11}(2u)\ts T_{22}(2u+1)
+T_{21}(2u)\ts T_{12}(2u+1)\\
t_{0,1}(u)&\mapsto -\frac{1}{\sqrt 2}\ts\Big(T_{12}(2u)\ts T_{22}(2u+1)
+T_{22}(2u)\ts T_{12}(2u+1)\Big)\\
t_{1,-1}(u)&\mapsto -T_{21}(2u)\ts T_{21}(2u+1)\\
t_{1,0}(u)&\mapsto -\frac{1}{\sqrt 2}\ts
\Big(T_{21}(2u)\ts T_{22}(2u+1)+T_{22}(2u)\ts T_{21}(2u+1)\Big)\\
t_{1,1}(u)&\mapsto T_{22}(2u)\ts T_{22}(2u+1).
\eal
\een
\epr

\bpf
We start by showing that the mapping defines an algebra homomorphism.
We use a version of the well known fusion procedure for $R$-matrices;
see e.g. \cite{aacdfr:gb} and references therein.

Consider the tensor product space $(\CC^2)^{\ot 4}$.
As in \eqref{ybe}, we use subscripts of the $R$-matrix \eqref{Rugl}
or the permutation operator $P\in\End(\CC^2)^{\ot 2}$
to indicate the copies of $\CC^2$ where the operator acts.
In the following we consider $V\ot V$ as a natural subspace
of $(\CC^2)^{\ot 2}\ot(\CC^2)^{\ot 2}$.
Obviously, the operator
\ben
\frac{1+P_{12}}{2}\cdot \frac{1+P_{34}}{2}=\frac{1}{4}\cdot
R^{\tss\circ}_{12}(-1)\ts R^{\tss\circ}_{34}(-1)
\een
is a projection of $(\CC^2)^{\ot 2}\ot(\CC^2)^{\ot 2}$ to the subspace
$V\ot V$.
Let us set
\beql{rv}
R_{\tss V}(u)=\frac{1+P_{12}}{2}\cdot\frac{1+P_{34}}{2}\cdot
R^{\tss\circ}_{14}(2u-1)\ts
R^{\tss\circ}_{13}(2u)\ts R^{\tss\circ}_{24}(2u)\ts R^{\tss\circ}_{23}(2u+1).
\eeq
Since the $R$-matrix
$R^{\tss\circ}(u)$ satisfies the
Yang--Baxter equation \eqref{ybe}, we have the following
equivalent expression for $R_{\tss V}(u)$,
\beql{rvequiv}
R_{\tss V}(u)=R^{\tss\circ}_{23}(2u+1)\ts R^{\tss\circ}_{13}(2u)\ts
R^{\tss\circ}_{24}(2u)\ts R^{\tss\circ}_{14}(2u-1)\cdot
\frac{1+P_{12}}{2}\cdot \frac{1+P_{34}}{2}.
\eeq
Clearly, the subspace $V\ot V$ is stable under the operator $R_{\tss V}(u)$.

\ble\label{lem:rmatrices}
We have the equality of operators in $V\ot V$,
\beql{equarr}
R_{\tss V}(u)=
\frac{2u-1}{2u+1}\cdot
\Big(1-\frac{P_{\tss V}}{u}+\frac{Q_{\tss V}}{u-1/2}\Big).
\eeq
\ele

\bpf
Using the formulas of the kind
\ben
(1+P_{12})\ts P_{14}\ts P_{24}=(1+P_{12})\ts P_{14}
\een
and
\ben
(1+P_{12})\ts (1+P_{34})\ts P_{14}\ts P_{23}
=(1+P_{12})\ts (1+P_{34})\ts  P_{13}\ts P_{24},
\een
it is easy to get a simplified expression for the operator $R_{\tss V}(u)$,
\ben
R_{\tss V}(u)=\frac{1+P_{12}}{2}\cdot \frac{1+P_{34}}{2}\cdot
\Big(1-\frac{P_{14}+P_{24}+P_{13}+P_{23}}{2u+1}+
\frac{P_{13}\ts P_{24}}{u\ts(2u+1)}\Big).
\een
The restriction of $R_{\tss V}(u)$ to the subspace $V\ot V$
is given by
\beql{simpop}
1-\frac{P_{14}+P_{24}+P_{13}+P_{23}}{2u+1}+
\frac{P_{13}\ts P_{24}}{u\ts(2u+1)}
\eeq
so that the proof of the lemma is completed by the application
of \eqref{simpop} to all basis vectors $v_i\ot v_j$ of $V\ot V$.
For instance, we have
\ben
\bal
R_{\tss V}(u)(v_{-1}\ot v_{-1})&=
\Big(1-\frac{P_{14}+P_{24}+P_{13}+P_{23}}{2u+1}+
\frac{P_{13}\ts P_{24}}{u\ts(2u+1)}\Big)(e_1\ot e_1\ot e_1\ot e_1)\\
&=\frac {(u-1)(2u-1)}{u\ts(2u+1)}\cdot e_1\ot e_1\ot e_1\ot e_1.
\eal
\een
Clearly, the application of the operator on the right hand side of
\eqref{equarr} to the vector $v_{-1}\ot v_{-1}$ gives the same result.
The remaining cases are verified by the same calculation.
\epf

By the lemma, the element $R_V(u)$
coincides with the $R$-matrix \eqref{Ru} for $\agot=\oa_3$,
up to a scalar factor. So, in order to verify that
the mapping \eqref{hommatr} defines a homomorphism
$\X(\oa_3)\to\Y(\gl_2)$ we need to show that
the relation
\beql{rttcheck}
R_V(u-v)\ts T_{1'}(u)\ts T_{2'}(v)=T_{2'}(v)\ts T_{1'}(u)\ts R_V(u-v)
\eeq
remains valid when $T(u)$ is replaced by its image.
Here we use primed indices to indicate the copies of the space $V$
in the tensor product $V\ot V$. We reserve unprimed indices
for the copies of $\CC^2$ in the tensor product $(\CC^2)^{\ot 4}$.
The left hand side of \eqref{rttcheck} reads
\ben
\bal
\frac{1+P_{12}}{2}\cdot\frac{1+P_{34}}{2}&\cdot
R^{\tss\circ}_{14}(2u-2v-1)\ts
R^{\tss\circ}_{13}(2u-2v)\ts R^{\tss\circ}_{24}(2u-2v)\ts
R^{\tss\circ}_{23}(2u-2v+1)\\
{}\times\frac{1+P_{12}}{2}&\cdot T^{\tss\circ}_1(2u)\ts T^{\tss\circ}_2(2u+1)
\cdot
\frac{1+P_{34}}{2}\cdot T^{\tss\circ}_3(2v)\ts T^{\tss\circ}_4(2v+1).
\eal
\een
Writing the product of $R$-matrices in the equivalent form
\eqref{rvequiv}, we simplify this expression to
\ben
\bal
\frac{1+P_{12}}{2}\cdot\frac{1+P_{34}}{2}&\cdot
R^{\tss\circ}_{14}(2u-2v-1)\ts
R^{\tss\circ}_{13}(2u-2v)\ts R^{\tss\circ}_{24}(2u-2v)\ts
R^{\tss\circ}_{23}(2u-2v+1)\\
{}\times{}& T^{\tss\circ}_1(2u)\ts T^{\tss\circ}_2(2u+1)
\ts T^{\tss\circ}_3(2v)\ts T^{\tss\circ}_4(2v+1).
\eal
\een
Now apply the $RTT$-relation \eqref{RTTgl} repeatedly to bring this
expression to the form
\ben
\bal
&\frac{1+P_{12}}{2}\cdot\frac{1+P_{34}}{2}\cdot
T^{\tss\circ}_3(2v)\ts T^{\tss\circ}_4(2v+1)\ts
T^{\tss\circ}_1(2u)\ts T^{\tss\circ}_2(2u+1)\\
{}\times{}&R^{\tss\circ}_{14}(2u-2v-1)\ts
R^{\tss\circ}_{13}(2u-2v)\ts R^{\tss\circ}_{24}(2u-2v)\ts
R^{\tss\circ}_{23}(2u-2v+1).
\eal
\een
Finally, since $R^{\tss\circ}_{12}(-1)/2$ is a projection,
we derive from \eqref{RTTgl} that
\ben
\frac{1+P_{12}}{2}\cdot T^{\tss\circ}_1(2u)\ts T^{\tss\circ}_2(2u+1)
=\frac{1+P_{12}}{2}\cdot T^{\tss\circ}_1(2u)\ts T^{\tss\circ}_2(2u+1)
\cdot\frac{1+P_{12}}{2}.
\een
Using the same property of $R^{\tss\circ}_{34}(-1)/2$ we obtain that
the resulting expression coincides with the right hand side
of \eqref{rttcheck}, where $T(u)$ is replaced with its image
in accordance with \eqref{hommatr}.

The explicit images of the generators of $\X(\oa_3)$ are found by
taking the matrix elements in \eqref{hommatr}.
Indeed, the application of $T(u)$ to the basis vector $v_{-1}$ of $V$
gives
\ben
T(u)(v_{-1})=t_{-1,-1}(u)\ts v_{-1}+t_{0,-1}(u)\ts v_{0}+
t_{1,-1}(u)\ts v_{1},
\een
while
\ben
\frac{1+P_{12}}{2}\cdot T^{\tss\circ}_1(2u)\ts T^{\tss\circ}_2(2u+1)(v_{-1})=
\frac{1+P_{12}}{2}\cdot T^{\tss\circ}_1(2u)\ts T^{\tss\circ}_2(2u+1)(e_1\ot e_1)
\een
\ben
\bal
{}&=\frac12\ts\sum_{a,b=1}^2 T_{a1}(2u)\ts T_{b1}(2u+1)\ts(e_a\ot e_b+e_b\ot e_a)
=T_{11}(2u)\ts T_{11}(2u+1)\ts v_{-1}\\
{}&+\frac{1}{\sqrt 2}\ts
\Big(T_{11}(2u)\ts T_{21}(2u+1)+T_{21}(2u)\ts T_{11}(2u+1)\Big)\ts v_0
-T_{21}(2u)\ts T_{21}(2u+1)\ts v_1.
\eal
\een
This agrees with the formulas for the images of the series $t_{a,-1}(u)$
for $a=-1,0,1$ given in the statement.
The remaining formulas are verified in the same way.
Note also that the image of the series $t_{0,0}(u)$ can be equivalently
written as
\ben
T_{12}(2u)\ts T_{21}(2u+1)
+T_{22}(2u)\ts T_{11}(2u+1)
\een
due to the defining relation in the Yangian $\Y(\gl_2)$.

In order to complete the proof of the proposition we now verify
that the homomorphism
$\X(\oa_3)\to\Y(\gl_2)$ given by \eqref{hommatr}
is bijective. Taking the coefficient at $u^{-r}$ in $t_{-1,-1}(u)$
we find that for any $r\geqslant 1$
\ben
t_{-1,-1}^{(r)}\mapsto 2^{-r+1}\ts T_{11}^{(r)}+
A_{r-1}(T_{11}^{(1)},\dots,T_{11}^{(r-1)}),
\een
where $A_{r-1}$ stands for a quadratic polynomial in the generators
$T_{11}^{(1)},\dots,T_{11}^{(r-1)}$. The obvious induction on $r$ shows that
each generator $T_{11}^{(r)}$ with $r\geqslant 1$ belongs to the image
of the homomorphism. Similarly, taking the image of $t_{1,1}^{(r)}$
we find that each generator $T_{22}^{(r)}$
with $r\geqslant 1$ also belongs to the image.
Then taking the images of $t_{-1,0}^{(r)}$ and $t_{0,-1}^{(r)}$
we derive the same property of
the generators $T_{12}^{(r)}$ and $T_{21}^{(r)}$ with $r\geqslant 1$.
This proves that the homomorphism is surjective.

Finally, observe that the homomorphism preserves the respective
filtrations on $\X(\oa_3)$ and $\Y(\gl_2)$.
Hence, we have a homomorphism of the associated graded algebras
$\gr\X(\oa_3)\to \gr\Y(\gl_2)$. It suffices to show that this homomorphism
is injective. Identifying $\gr\Y(\gl_2)$ with the universal enveloping algebra
$\U(\gl_2[x])$ via the isomorphism \eqref{isomgrgl}, we get
\ben
\bar t_{0,1}^{\ts(r)}\mapsto -\frac{1}{\sqrt 2}\ts E_{12}\ts (x/2)^{r-1},
\qquad
\bar t_{1,0}^{\ts(r)}\mapsto -\frac{1}{\sqrt 2}\ts E_{21}\ts (x/2)^{r-1},
\qquad
\bar t_{1,1}^{\ts(r)}\mapsto E_{22}\ts (x/2)^{r-1}
\een
and
\ben
\bar z_{r}\mapsto (E_{11}+E_{22})\ts (x/2)^{r-1}.
\een
Therefore, the injectivity
of the homomorphism follows from Corollary~\ref{cor:pbwx}.
\epf

\bco\label{cor:restb1}
The restriction of the isomorphism $\phi:\X(\oa_3)\to\Y(\gl_2)$
to the subalgebra $\Y(\oa_3)$ induces an isomorphism
$\Y(\oa_3)\to\Y(\sll_2)$.
\eco

\bpf
Recall that the subalgebra $\Y(\oa_3)$ consists of the elements
stable under all automorphisms of $\X(\oa_3)$ of the form
\eqref{mult}. For any series $f(u)$ of the form \eqref{fu}
there exists a unique series
\ben
g(u)=1+g_1u^{-1}+g_2u^{-2}+\cdots\in\CC[[u^{-1}]]
\een
such that $f(u)=g(2u)\ts g(2u+1)$. By Proposition~\ref{prop:b1},
we have the relation
$\phi\circ\mu_f=\mu_g\circ\phi$, and hence
$\mu_f\circ\phi^{-1}=\phi^{-1}\circ\mu_g$.
This implies that the image
of $\Y(\oa_3)$ under the isomorphism $\phi$ coincides with the subalgebra
$\Y(\sll_2)$ of $\Y(\gl_2)$ thus yielding
the desired isomorphism.
\epf

Let us denote by $c$ the Casimir element for the Lie algebra
$\oa_3$,
\ben
c=\frac12\big(F_{11}^2-F_{11}\big)+F_{10}F_{01}.
\een
In the following we use notation \eqref{F}.

\bco\label{cor:evb1}
The mapping
\beql{evalb1}
\ev:
T(u)\mapsto 1+\frac{F}{u}+\frac{F^2-c\ts 1}{u\ts(2u-1)}
\eeq
defines a surjective homomorphism $\X(\oa_3)\to\U(\oa_3)$.
\eco

\bpf
Writing the homomorphism $\Y(\gl_2)\to\U(\sll_2)$ used
in the proof of Corollary~\ref{cor:evc1} in the current notation
we get
\ben
\bal
T_{11}(u)&\mapsto 1+\big(E_{11}-E_{22}\big)\ts(2u)^{-1},\qquad
T_{12}(u)\mapsto E_{12}u^{-1},\\
T_{22}(u)&\mapsto 1+\big(E_{22}-E_{11}\big)\ts(2u)^{-1},
\qquad
T_{21}(u)\mapsto E_{21}u^{-1}.
\eal
\een
Composing this with the isomorphism $\sll_2\to\oa_3$ given by
\ben
E_{11}-E_{22}\mapsto 2F_{-1,-1},\qquad
E_{12}\mapsto\sqrt2\ts F_{-1,0},\qquad
E_{21}\mapsto\sqrt2\ts F_{0,-1},
\een
we get another homomorphism $\Y(\gl_2)\to\U(\oa_3)$ such that
\begin{alignat}{2}
T_{11}(u)&\mapsto 1+F_{-1,-1}\ts u^{-1},\qquad
&&T_{12}(u)\mapsto \sqrt2\ts  F_{-1,0}\ts u^{-1},
\non\\
T_{22}(u)&\mapsto 1+F_{1,1}\ts u^{-1},
\qquad
&&T_{21}(u)\mapsto \sqrt2\ts  F_{0,-1}\ts u^{-1}.
\non
\end{alignat}
Finally, compose the isomorphism of Proposition~\ref{prop:b1}
with the shift automorphism $t_{ij}(u)\mapsto t_{ij}(u-1/2)$
of $\X(\oa_3)$ and use the above formulas to get
a homomorphism $\X(\oa_3)\to\U(\oa_3)$. It remains to verify that
the resulting formulas for the images of $t_{ij}(u)$ agree with
\eqref{evalb1}. This can be done by an easy straightforward calculation.
For instance, for the image of $t_{0,0}(u)$ we calculate
\beql{calcoo}
\bal
t_{0,0}(u)&\mapsto
T_{11}(2u-1)\ts T_{22}(2u)
+T_{21}(2u-1)\ts T_{12}(2u)\\
{}&{}\mapsto
\Big(1+\frac{F_{-1,-1}}{2u-1}\Big)\Big(1+\frac{F_{1,1}}{2u}\Big)
+2\cdot\frac{F_{0,-1}}{2u-1}\cdot\frac{F_{-1,0}}{2u}.
\eal
\eeq
On the other hand, formula \eqref{evalb1} gives
\ben
\bal
t_{0,0}(u)&\mapsto 1+\frac{2F_{0,-1}F_{-1,0}+2F_{0,1}F_{1,0}-
F_{11}^2+F_{11}-2\ts F_{10}F_{01}}{2u\ts (2u-1)}\\
{}&{}=1+\frac{2F_{0,-1}F_{-1,0}-
F_{11}^2-F_{11}}{2u\ts (2u-1)}.
\eal
\een
Clearly, this agrees with \eqref{calcoo}. All the remaining cases
are verified by a similar and even shorter calculation.
Obviously, the homomorphism \eqref{evalb1} is surjective.
\epf

\subsection{Extended Yangian $\X(\oa_4)$}\label{subsec:d2}

We shall need the tensor product algebra $\Y(\gl_2)\ot\Y(\gl_2)$.
In order to distinguish the two copies of the algebra $\Y(\gl_2)$,
we denote the corresponding generator series respectively
by $T_{ij}(u)$ and $T^{\ts\prime}_{ij}(u)$
for the first and second copies, where $i,j\in\{1,2\}$.
We also identify $T_{ij}(u)\ot 1$ with $T_{ij}(u)$ and
$1\ot T^{\ts\prime}_{ij}(u)$ with $T^{\ts\prime}_{ij}(u)$.
As before, we combine the series $T_{ij}(u)$ and $T^{\ts\prime}_{ij}(u)$
into the matrices  $T^{\tss\circ}(u)$ and $T^{\tss\circ\ts\prime}(u)$,
respectively.

The algebra $\Y(\gl_2)\ot\Y(\gl_2)$ is naturally equipped
with an ascending filtration, where the degrees of the elements on
each copy of $\Y(\gl_2)$ are defined by \eqref{deggl}.

Consider the
vector space $\CC^2$ with its canonical basis $e_1,e_2$ and
set $V=\CC^2\ot\CC^2$.
We identify $V$ with $\CC^4$ regarding the vectors
\ben
v_{-2}=e_1\ot e_1,
\qquad v_{-1}=e_1\ot e_2,
\qquad v_{1}=e_2\ot e_1,
\qquad v_{2}=-e_2\ot e_2
\een
as the canonical basis of $V$. Then
$T^{\tss\circ}_1(u)\ts T^{\tss\circ\ts\prime}_2(u)$
may be regarded as an element of
$\End V\ot\big(\Y(\gl_2)\ot\Y(\gl_2)\big)[[u^{-1}]]$.
The operators $P_{\tss V}$ and $Q_{\tss V}$
in $V\ot V$ are given by the respective formulas \eqref{P} and \eqref{Q}.
We shall regard the matrix $T(u)=(t_{ij}(u))$
as an element of $\End V\ot\X(\oa_4)[[u^{-1}]]$.

\bpr\label{prop:d2}
The mapping
\beql{hommatrd}
T(u)\mapsto
T^{\tss\circ}_1(u)\ts T^{\tss\circ\ts\prime}_2(u),
\eeq
defines an embedding
$\psi:\X(\oa_4)\hookrightarrow\Y(\gl_2)\ot\Y(\gl_2)$.
More explicitly, the images of the generators
under the embedding are given by the formulas
\begin{align}
t_{-2,-2}(u)&\mapsto \phantom{-}T_{11}(u)\ts T^{\ts\prime}_{11}(u),
\quad&
t_{-2,-1}(u)&\mapsto \phantom{-}T_{11}(u)\ts T^{\ts\prime}_{12}(u),\non\\
t_{-2,1}(u)&\mapsto \phantom{-}T_{12}(u)\ts T^{\ts\prime}_{11}(u),
\quad&
t_{-2,2}(u)&\mapsto -T_{12}(u)\ts T^{\ts\prime}_{12}(u),\non\\
t_{-1,-2}(u)&\mapsto \phantom{-}T_{11}(u)\ts T^{\ts\prime}_{21}(u),
\quad&
t_{-1,-1}(u)&\mapsto \phantom{-}T_{11}(u)\ts T^{\ts\prime}_{22}(u),\non\\
t_{-1,1}(u)&\mapsto \phantom{-}T_{12}(u)\ts T^{\ts\prime}_{21}(u),
\quad&
t_{-1,2}(u)&\mapsto -T_{12}(u)\ts T^{\ts\prime}_{22}(u),\non\\
t_{1,-2}(u)&\mapsto \phantom{-}T_{21}(u)\ts T^{\ts\prime}_{11}(u),
\quad&
t_{1,-1}(u)&\mapsto \phantom{-}T_{21}(u)\ts T^{\ts\prime}_{12}(u),\non\\
t_{1,1}(u)&\mapsto \phantom{-}T_{22}(u)\ts T^{\ts\prime}_{11}(u),
\quad&
t_{1,2}(u)&\mapsto -T_{22}(u)\ts T^{\ts\prime}_{12}(u),\non\\
t_{2,-2}(u)&\mapsto -T_{21}(u)\ts T^{\ts\prime}_{21}(u),
\quad&
t_{2,-1}(u)&\mapsto -T_{21}(u)\ts T^{\ts\prime}_{22}(u),\non\\
t_{2,1}(u)&\mapsto -T_{22}(u)\ts T^{\ts\prime}_{21}(u),
\quad&
t_{2,2}(u)&\mapsto \phantom{-}T_{22}(u)\ts T^{\ts\prime}_{22}(u).\non
\end{align}
\epr

\bpf
We start by showing that the mapping defines an algebra homomorphism.
Identifying $V\ot V$ with
the tensor product space $(\CC^2)^{\ot4}$, we set
\beql{rvd}
R_{\tss V}(u)=R^{\tss\circ}_{13}(u)\ts
R^{\tss\circ}_{24}(u).
\eeq

\ble\label{lem:rmatricesd}
We have the equality of operators in $V\ot V$,
\beql{equarrd}
R_{\tss V}(u)=
\frac{u-1}{u}\cdot
\Big(1-\frac{P_{\tss V}}{u}+\frac{Q_{\tss V}}{u-1}\Big).
\eeq
\ele

\bpf
We have
\ben
R_{\tss V}(u)=
\Big(1-\frac{P_{13}}{u}\Big)\Big(1-\frac{P_{24}}{u}\Big)
=\frac{u-1}{u}\Big(1-\frac{P_{13}P_{24}}{u}
+\frac{(1-P_{13})(1-P_{24})}{u-1}\Big).
\een
It remains to note that $P_{\tss V}=P_{13}P_{24}$
and $Q_{\tss V}=(1-P_{13})(1-P_{24})$.
This is verified by the application
of the operators to all basis vectors $v_i\ot v_j$ of $V\ot V$.
For instance, by the definition of $Q_{\tss V}$,
\ben
Q_{\tss V}(v_{-2}\ot v_2)=v_{-2}\ot v_2+v_{-1}\ot v_1+
v_{1}\ot v_{-1}+v_{2}\ot v_{-2},
\een
while
\ben
\bal
&{}(1-P_{13})(1-P_{24})(v_{-2}\ot v_2)=
(1-P_{13})(1-P_{24})(-e_1\ot e_1\ot e_2\ot e_2)\\
{}&{}
=-e_1\ot e_1\ot e_2\ot e_2+e_1\ot e_2\ot e_2\ot e_1+e_2\ot e_1\ot e_1\ot e_2
-e_2\ot e_2\ot e_1\ot e_1,
\eal
\een
which clearly coincides with $Q_{\tss V}(v_{-2}\ot v_2)$. The remaining
relations are verified in the same way.
\epf

By the lemma, the element $R_V(u)$
coincides with the $R$-matrix \eqref{Ru} for $\agot=\oa_4$,
up to a scalar factor. So, in order to verify that
the mapping \eqref{hommatrd} defines a homomorphism
$\X(\oa_4)\to\Y(\gl_2)\ot \Y(\gl_2)$ we need to show that
the relation
\beql{rttcheckd}
R_V(u-v)\ts T_{1'}(u)\ts T_{2'}(v)=T_{2'}(v)\ts T_{1'}(u)\ts R_V(u-v)
\eeq
remains valid when $T(u)$ is replaced by its image.
The primed indices are used here to indicate the copies of the space $V$
in the tensor product $V\ot V$.
The left hand side of \eqref{rttcheckd} reads
\ben
R^{\tss\circ}_{13}(u-v)\ts R^{\tss\circ}_{24}(u-v)\ts
T^{\tss\circ}_1(u)\ts T^{\tss\circ\ts\prime}_2(u)\ts
T^{\tss\circ}_3(v)\ts T^{\tss\circ\ts\prime}_4(v).
\een
Applying the $RTT$-relation \eqref{RTTgl} twice, we bring this
expression to the form
\ben
T^{\tss\circ}_3(v)\ts T^{\tss\circ\ts\prime}_4(v)\ts
T^{\tss\circ}_1(u)\ts T^{\tss\circ\ts\prime}_2(u)\ts
R^{\tss\circ}_{13}(u-v)\ts R^{\tss\circ}_{24}(u-v)
\een
which coincides with the right hand side
of \eqref{rttcheckd}, where $T(u)$ is replaced with its image
in accordance with \eqref{hommatrd}.

The explicit images of the generators of $\X(\oa_4)$ are found by
taking the matrix elements in \eqref{hommatrd}.
Indeed, the application of $T(u)$ to the basis vector $v_{-2}$ of $V$
gives
\ben
T(u)(v_{-2})=t_{-2,-2}(u)\ts v_{-2}+t_{-1,-2}(u)\ts v_{-1}+
t_{1,-2}(u)\ts v_{1}+t_{2,-2}(u)\ts v_{2},
\een
while
\ben
T^{\tss\circ}_1(u)\ts T^{\tss\circ\ts\prime}_2(u)(v_{-2})=
T^{\tss\circ}_1(u)\ts T^{\tss\circ\ts\prime}_2(u)(e_1\ot e_1)
=\sum_{a,b=1}^2 T_{a1}(u)\ts T^{\ts\prime}_{b1}(u)\ts(e_a\ot e_b)
\een
\ben
=T_{11}(u)\ts T^{\ts\prime}_{11}(u)\ts v_{-2}
+T_{11}(u)\ts T^{\ts\prime}_{21}(u)\ts v_{-1}+
T_{21}(u)\ts T^{\ts\prime}_{11}(u)\ts v_{1}-
T_{21}(u)\ts T^{\ts\prime}_{21}(u)\ts v_{2}.
\een
This agrees with the formulas for the images of the series $t_{a,-2}(u)$
for $a=-2,-1,1,2$ given in the statement.
The remaining formulas are verified in the same way.

In order to demonstrate that the homomorphism $\psi$ is injective,
observe that it preserves the respective
filtrations on $\X(\oa_4)$ and $\Y(\gl_2)\ot \Y(\gl_2)$.
Hence, we have a homomorphism of the associated graded algebras
\ben
\gr\X(\oa_4)\to \gr\big(\Y(\gl_2)\ot\Y(\gl_2)\big).
\een
Identifying the graded algebra $\gr\Y(\gl_2)$ with $\U(\gl_2[x])$
via the isomorphism \eqref{isomgrgl}, we get a homomorphism
\ben
\gr\X(\oa_4)\to \U(\gl_2[x])\ot\U(\gl_2[y])
\een
so that
\ben
\bal
\bar t_{-1,2}^{\ts(r)}&\mapsto -E_{12}\ts x^{r-1},
\qquad
\bar t_{1,2}^{\ts(r)}\mapsto -E_{12}\ts y^{r-1},
\qquad
\bar t_{1,1}^{\ts(r)}\mapsto E_{22}\ts x^{r-1}+E_{11}\ts y^{r-1},\\
\qquad
\bar t_{2,-1}^{\ts(r)}&\mapsto -E_{21}\ts x^{r-1},
\qquad
\bar t_{2,1}^{\ts(r)}\mapsto -E_{21}\ts y^{r-1},
\qquad
\bar t_{2,2}^{\ts(r)}\mapsto E_{22}\ts x^{r-1}+E_{22}\ts y^{r-1},
\eal
\een
and
\ben
\bar z_{r}\mapsto (E_{11}+E_{22})\ts x^{r-1}+(E_{11}+E_{22})\ts y^{r-1}.
\een
Therefore, the injectivity
of $\psi$ follows from Corollary~\ref{cor:pbwx}.
\epf

Due to the presentation of the Yangian $\Y(\sll_2)$ provided by
\eqref{dreltau} and \eqref{qdetone} we have a natural projection
$\Y(\gl_2)\to\Y(\sll_2)$ defined by the mapping
$T_{ij}(u)\mapsto\mathcal T_{ij}(u)$. Applying this projection
to the first or second copy of $\Y(\gl_2)$ in the tensor product algebra
$\Y(\gl_2)\ot\Y(\gl_2)$ and taking its composition with the
embedding $\psi$ we get homomorphisms
\ben
\chi^{(1)}:\X(\oa_4)\to\Y(\sll_2)\ot\Y(\gl_2),\qquad
\chi^{(2)}:\X(\oa_4)\to\Y(\gl_2)\ot\Y(\sll_2).
\een

\bco\label{cor:restdisom}
The homomorphisms $\chi^{(1)}$ and $\chi^{(2)}$ are bijective.
\eco

\bpf
We only consider $\chi^{(1)}$, the proof for $\chi^{(2)}$ is similar.
By the formulas of Proposition~\ref{prop:d2} we have
\ben
\bal
\chi^{(1)}:{}&{}t_{-2,-2}(u)\ts t_{1,1}(u-1)-t_{1,-2}(u)\ts t_{-2,1}(u-1)
\mapsto\\
&\Big(\mathcal T_{11}(u)\ts \mathcal T_{22}(u-1)-
\mathcal T_{21}(u)\ts \mathcal T_{12}(u-1)\Big)\ts
T^{\ts\prime}_{11}(u)\ts T^{\ts\prime}_{11}(u-1)
=T^{\ts\prime}_{11}(u)\ts T^{\ts\prime}_{11}(u-1).
\eal
\een
Therefore, all the coefficients of the series $T^{\ts\prime}_{11}(u)$
belong to the image of $\chi^{(1)}$. Hence, so do the coefficients of
$\mathcal T_{ij}(u)$ with $i,j\in\{1,2\}$. This implies that
$\chi^{(1)}$ is surjective.
To verify the injectivity of $\chi^{(1)}$ we use the same
argument as in the proof of Proposition~\ref{prop:d2}.
Namely, $\chi^{(1)}$ induces a homomorphism of the
associated graded algebras
\ben
\gr\X(\oa_4)\to \U(\sll_2[x])\ot\U(\gl_2[y])
\een
and the argument is completed by the application of
Corollary~\ref{cor:pbwx}.
\epf

\bco\label{cor:restd2}
The restriction of each isomorphism $\chi^{(1)}$ and $\chi^{(2)}$
to the subalgebra $\Y(\oa_4)$ induces an isomorphism
$\Y(\oa_4)\to\Y(\sll_2)\ot \Y(\sll_2)$.
\eco

\bpf
Again, we only consider the isomorphism $\chi^{(1)}$.
The subalgebra $\Y(\oa_4)$ consists of the elements
stable under all automorphisms of $\X(\oa_4)$ of the form
\eqref{mult}. For any formal series $f(u)$ of the form \eqref{fu}
consider the automorphism $\wt{\mu}_f$
of the algebra $\Y(\sll_2)\ot \Y(\gl_2)$ defined by
\ben
\wt{\mu}_f: \mathcal T_{ij}(u)\mapsto \mathcal T_{ij}(u),\qquad
T^{\ts\prime}_{ij}(u)\mapsto f(u)\ts T^{\ts\prime}_{ij}(u).
\een
By the definition of $\chi^{(1)}$, we have the relation
$\chi^{(1)}\circ\mu_f=\wt{\mu}_f\circ\chi^{(1)}$.
This implies that if $y\in \Y(\oa_4)$ then $\chi^{(1)}(y)$ is stable under
the automorphisms $\wt{\mu}_f$ for all series $f(u)$.
Hence, the image of the subalgebra $\Y(\oa_4)$
of $\X(\oa_4)$ under the isomorphism $\chi^{(1)}$ coincides with the subalgebra
$\Y(\sll_2)\ot\Y(\sll_2)$ of $\Y(\sll_2)\ot\Y(\gl_2)$ thus providing
the desired isomorphism.
\epf

Let us denote by $c$ the following Casimir element for the Lie algebra
$\oa_4$,
\ben
c=\frac12\big(F_{11}^2+F_{22}^2\big)-F_{22}+F_{21}F_{12}+F_{2,-1}F_{-1,2}.
\een
In the following we use notation \eqref{F}.

\bco\label{cor:evd2}
The mapping
\beql{evald2}
\ev:
T(u)\mapsto 1+\frac{F}{u}+\frac{F^2-F-c\ts 1}{2u^2}
\eeq
defines a surjective homomorphism $\X(\oa_4)\to\U(\oa_4)$.
\eco

\bpf
Consider the isomorphism $\sll_2\oplus\sll_2\to\oa_4$ given by
\ben
E_{11}-E_{22}\mapsto -F_{11}-F_{22},\qquad
E_{12}\mapsto F_{-2,1},\qquad
E_{21}\mapsto F_{1,-2},
\een
and
\ben
E^{\ts\prime}_{11}-E^{\ts\prime}_{22}\mapsto F_{11}-F_{22},\qquad
E^{\ts\prime}_{12}\mapsto F_{-2,-1},\qquad
E^{\ts\prime}_{21}\mapsto F_{-1,-2},
\een
where the primes indicate the basis elements of the second copy
of $\sll_2$. Applying the homomorphism $\Y(\gl_2)\to\U(\sll_2)$ used
in the proof of Corollary~\ref{cor:evb1},
we get a homomorphism $\Y(\gl_2)\ot\Y(\gl_2)\to\U(\oa_4)$ such that
\begin{alignat}{2}
T_{11}(u)&\mapsto 1-\frac{F_{11}+F_{22}}{2}\ts u^{-1},\qquad
&&T_{12}(u)\mapsto F_{-2,1}\ts u^{-1},
\non\\
T_{22}(u)&\mapsto 1+\frac{F_{11}+F_{22}}{2}\ts u^{-1},
\qquad
&&T_{21}(u)\mapsto F_{1,-2}\ts u^{-1}
\non
\end{alignat}
and
\begin{alignat}{2}
T^{\ts\prime}_{11}(u)&\mapsto 1+\frac{F_{11}-F_{22}}{2}\ts u^{-1},\qquad
&&T^{\ts\prime}_{12}(u)\mapsto F_{-2,-1}\ts u^{-1},
\non\\
T^{\ts\prime}_{22}(u)&\mapsto 1-\frac{F_{11}-F_{22}}{2}\ts u^{-1},
\qquad
&&T^{\ts\prime}_{21}(u)\mapsto F_{-1,-2}\ts u^{-1}.
\non
\end{alignat}
Using the isomorphism of Proposition~\ref{prop:d2}
we get a homomorphism $\X(\oa_4)\to\U(\oa_4)$. It remains to verify that
the resulting formulas for the images of $t_{ij}(u)$ agree with
\eqref{evald2}. This can be done by an easy straightforward calculation.
For instance, for the image of $t_{-2,-2}(u)$ we calculate
\beql{calcood2}
\bal
t_{-2,-2}(u)\mapsto
T_{11}(u)\ts T^{\ts\prime}_{11}(u)
{}&{}\mapsto
\Big(1-\frac{F_{11}+F_{22}}{2}\ts u^{-1}\Big)
\Big(1+\frac{F_{11}-F_{22}}{2}\ts u^{-1}\Big)\\
{}&{}=
1+F_{-2,-2}\ts u^{-1}+\frac{F_{-2,-2}^2-F_{-1,-1}^2}{4}\ts u^{-2}.
\eal
\eeq
On the other hand, formula \eqref{evald2} gives
\ben
t_{-2,-2}(u)\mapsto 1+F_{-2,-2}\ts u^{-1}+
\frac{F_{-2,-2}^2+F_{-2,-1}F_{-1,-2}+F_{-2,1}F_{1,-2}
-F_{-2,-2}-c}{2u^2}
\een
which agrees with \eqref{calcood2}. All the remaining cases
are verified by a similar calculation.
Obviously, the homomorphism \eqref{evald2} is surjective.
\epf

\bre\label{rem:evala}
The respective compositions of the evaluation homomorphisms
provided by Corollaries~\ref{cor:evc1}, \ref{cor:evb1} and \ref{cor:evd2}
with the shift automorphism $\tau_a$ given by \eqref{shift} yields
the homomorphisms $\ev_a=\ev\circ\tau_a$ with the evaluation parameter $a$.
\qed
\ere

\section{Representations of the extended Yangians}\label{sec:rep}
\setcounter{equation}{0}

Here we introduce the highest weight representations for the extended
Yangians $\X(\agot)$, where as before, $\agot=\oa_{2n+1}$, $\spa_{2n}$
or $\oa_{2n}$.
We show by a standard argument that
finite-dimensional irreducible
representations of $\X(\agot)$ are highest weight representations.
Then we give necessary and sufficient conditions for the irreducible
highest weight representations to be finite-dimensional.
In particular, we obtain an alternative
proof of Drinfeld's classification theorem
for the finite-dimensional irreducible
representations of the Yangians $\Y(\agot)$.

\subsection{Highest weight representations}\label{subsec:hw}

A representation $V$ of the algebra $\X(\agot)$
is called a {\it highest weight representation\/}
if there exists a nonzero vector
$\xi\in V$ such that $V$ is generated by $\xi$,
\beql{trianb}
\begin{alignedat}{2}
t_{ij}(u)\ts\xi&=0 \qquad &&\text{for}
\quad -n\leqslant i<j\leqslant n, \qquad \text{and}\\
t_{ii}(u)\ts\xi&=\la_i(u)\ts\xi \qquad &&\text{for}
\quad -n\leqslant i\leqslant n,
\end{alignedat}
\eeq
for some formal series
\beql{laiu}
\la_i(u)=1+\la_i^{(1)}u^{-1}+\la_i^{(2)}u^{-2}+
\cdots,\qquad
\la_i^{(r)}\in\CC,
\eeq
where the value $i=0$ only occurs in the case $\agot=\oa_{2n+1}$.
The vector $\xi$ is called the {\it highest vector\/}
of $V$ and the tuple $\la(u)=(\la_{-n}(u),\dots,\la_n(u))$
of the formal series is the {\it highest weight\/} of $V$.

Let us identify the elements
$F_{ij}\in\agot$ with their
images in $\X(\agot)$ under the embedding \eqref{embx}.
The defining relations
\eqref{defrel} imply
\ben
[t_{ij}^{(1)},t_{kl}(u)]=
\de_{kj}\ts t_{il}(u)-\de_{il}\ts t_{kj}(u)
-\de_{k,-i}\ts\theta_{ij}\ts t_{-j,l}(u)
+\de_{l,-j}\ts\theta_{ij}\ts t_{k,-i}(u).
\een
Also, due to \eqref{zt} we have
\ben
t_{ij}^{(1)}+\theta_{ij}\ts t_{-j,-i}^{(1)}=\de_{ij}\ts z_1.
\een
Therefore, $F_{ij}=t_{ij}^{(1)}-\de_{ij}\ts z_1/2$.
Since $z_1$ is central in $\X(\agot)$, this gives
\beql{tijF}
[F_{ij},t_{kl}(u)]=
\de_{kj}\ts t_{il}(u)-\de_{il}\ts t_{kj}(u)
-\de_{k,-i}\ts\theta_{ij}\ts t_{-j,l}(u)
+\de_{l,-j}\ts\theta_{ij}\ts t_{k,-i}(u).
\eeq

Take the linear span of the elements
$F_{11},\dots,F_{nn}$ as the Cartan subalgebra $\mathfrak{h}$ of $\agot$
and consider the standard triangular decomposition of $\agot$.
Then the nonzero elements $F_{ij}$ with $i<j$ are the positive root vectors.
The corresponding positive roots are
\ben
-\ve_i-\ve_j,\quad\ve_i-\ve_j\quad\text{with}\quad 1\leqslant i<j\leqslant n
\een
for $\agot=\oa_{2n}$,
\ben
-2\ts\ve_i\quad\text{with}\quad 1\leqslant i\leqslant n\fand
-\ve_i-\ve_j,\quad\ve_i-\ve_j\quad\text{with}\quad 1\leqslant i<j\leqslant n
\een
for $\agot=\spa_{2n}$, and
\ben
-\ve_i\quad\text{with}\quad 1\leqslant i\leqslant n\fand
-\ve_i-\ve_j,\quad\ve_i-\ve_j\quad\text{with}\quad 1\leqslant i<j\leqslant n
\een
for $\agot=\oa_{2n+1}$, where $\ve_i$ denotes the element of
$\mathfrak{h}^*$ defined by $\ve_i(F_{jj})=\de_{ij}$.
The standard partial ordering on the set of weights
of any $\agot$-module
is now defined as follows.
If $\al$ and $\be$ are two weights, then
$\al$ precedes $\be$ if $\be-\al$
is a $\ZZ_+$-linear combination of the positive roots.

\bth\label{thm:fdhw}
Every finite-dimensional irreducible representation $V$ of the algebra $\X(\agot)$
is a highest weight representation. Moreover, $V$ contains a unique,
up to a constant factor, highest vector.
\eth

\bpf
Introduce the subspace $V^{\tss0}$ of $V$ by
\beql{vo}
V^{\tss0}=\{\eta\in V\ |\ t_{ij}(u)\ts\eta=0,\qquad -n\leqslant i<j\leqslant n\}.
\eeq
We show first that $V^{\tss0}$ is nonzero.
Consider the set of weights
of $V$, where $V$ is regarded as the $\agot$-module defined via the embedding
\eqref{embx}. This set is finite and hence contains a maximal weight $\nu$
with respect to the partial ordering on the set of weights of $V$.
The corresponding weight vector $\eta$ belongs to $V^{\tss0}$.
Indeed, if $i<j$ then by \eqref{tijF} the weight of $t_{ij}(u)\ts\eta$
has the form $\nu+\al$ for a positive root $\al$.
By the maximality of $\nu$,
we have $t_{ij}(u)\ts\eta=0$.

Next, we show that all the operators $t_{kk}(u)$
preserve the subspace $V^{\tss0}$. Consider first the case
$\agot=\oa_{2n+1}$.
In the following argument we write $\equiv$ for an equality
of operators in $V^{\tss0}$.
Due to \eqref{tijF}, it suffices to show that
for any $i$ and $k$ we have
\beql{tiipotkk}
t_{i,i+1}(u)\ts t_{kk}(v)\equiv 0.
\eeq
Suppose first that $i<k$. Then \eqref{tiipotkk}
is immediate from \eqref{defrel} except for the cases
$i=-k$ and $i=-k-1$. In the former case, we have $k>0$ and
so \eqref{defrel} gives
\beql{tkkkv}
t_{-k,-k+1}(u)\ts t_{kk}(v)\equiv
-\frac{1}{u-v-\kappa}
\sum_{p=k}^n t_{-p,-k+1}(u)\ts t_{pk}(v),
\eeq
while for each $p\geqslant k$,
\ben
t_{-p,-k+1}(u)\ts t_{pk}(v)\equiv
-\frac{1}{u-v-\kappa}
\sum_{q=k}^n t_{-q,-k+1}(u)\ts t_{qk}(v).
\een
Hence, $t_{-p,-k+1}(u)\ts t_{pk}(v)\equiv t_{-k,-k+1}(u)\ts t_{kk}(v)$.
So, \eqref{tkkkv} implies
\ben
\left(1+\frac{n-k+1}{u-v-\kappa}\right)\ts t_{-k,-k+1}(u)\ts t_{kk}(v)\equiv 0
\een
and thus,
$
t_{-k,-k+1}(u)\ts t_{kk}(v)\equiv 0
$
verifying \eqref{tiipotkk}.

Similarly, in the case $i=-k-1$ we have $k\geqslant 0$ and so
\beql{tkkkvm}
t_{-k-1,-k}(u)\ts t_{kk}(v)\equiv
\frac{1}{u-v-\kappa}
\sum_{p=k+1}^n t_{kp}(v)\ts t_{-k-1,-p}(u).
\eeq
For each $p\geqslant k+1$ we have
\ben
t_{kp}(v)\ts t_{-k-1,-p}(u)\equiv
-[t_{-k-1,-p}(u),t_{kp}(v)]\equiv
-\frac{1}{u-v-\kappa}
\sum_{q=k+1}^n t_{kq}(v)\ts t_{-k-1,-q}(u).
\een
Therefore, $t_{kp}(v)\ts t_{-k-1,-p}(u)\equiv -t_{-k-1,-k}(u)\ts t_{kk}(v)$.
So, \eqref{tkkkvm} gives
\ben
\left(1+\frac{n-k}{u-v-\kappa}\right)\ts t_{-k-1,-k}(u)\ts t_{kk}(v)\equiv 0
\een
implying \eqref{tiipotkk} in the case under consideration.

Suppose now that $i\geqslant k$. We can write
\ben
t_{i,i+1}(u)\ts t_{kk}(v)\equiv -[t_{kk}(v),\ts t_{i,i+1}(u)].
\een
Now \eqref{tiipotkk}
is immediate from \eqref{defrel} except for the cases
$i=-k$ and $i=-k-1$. In the former case, we have $k\leqslant 0$ and
so \eqref{defrel} gives
\beql{tkkkvmo}
t_{-k,-k+1}(u)\ts t_{kk}(v)\equiv
\frac{1}{v-u-\kappa}
\sum_{p=-k+1}^n t_{-p,k}(v)\ts t_{p,-k+1}(u),
\eeq
while for each $p\geqslant -k+1$,
\ben
t_{-p,k}(v)\ts t_{p,-k+1}(u)\equiv
-\frac{1}{v-u-\kappa}
\sum_{q=-k+1}^n t_{-q,k}(v)\ts t_{q,-k+1}(u).
\een
Hence, $t_{-p,k}(v)\ts t_{p,-k+1}(u)\equiv -t_{-k,-k+1}(u)\ts t_{kk}(v)$.
So, \eqref{tkkkvmo} implies
\ben
\left(1+\frac{n+k}{v-u-\kappa}\right)\ts t_{-k,-k+1}(u)\ts t_{kk}(v)\equiv 0
\een
verifying \eqref{tiipotkk}.

Finally, let $i=-k-1$. Then $k<0$ and
\beql{tkkkvmom}
t_{-k-1,-k}(u)\ts t_{kk}(v)\equiv -[t_{kk}(v),\ts t_{-k-1,-k}(u)]
\equiv
-\frac{1}{v-u-\kappa}
\sum_{p=-k}^n t_{-k-1,p}(u)\ts  t_{k,-p}(v).
\eeq
For each $p\geqslant -k$ we have
\ben
t_{-k-1,p}(u)\ts  t_{k,-p}(v)\equiv -[t_{k,-p}(v),\ts t_{-k-1,p}(u)]
\equiv -\frac{1}{v-u-\kappa}
\sum_{q=-k}^n t_{-k-1,q}(u)\ts  t_{k,-q}(v).
\een
Therefore, $t_{-k-1,p}(u)\ts  t_{k,-p}(v)\equiv t_{-k-1,-k}(u)\ts t_{kk}(v)$.
So, \eqref{tkkkvmom} gives
\ben
\left(1+\frac{n+k+1}{v-u-\kappa}\right)\ts t_{-k-1,-k}(u)\ts t_{kk}(v)\equiv 0
\een
completing the proof of \eqref{tiipotkk}.

For the Lie algebras $\agot=\spa_{2n}$ and $\oa_{2n}$ the argument is
essentially the same as in
the previous case. If $\agot=\spa_{2n}$, then
due to \eqref{tijF}, it suffices to show that
\eqref{tiipotkk} holds for $i\in\{-n,\dots,-2,1,\dots,n-1\}$ and all $k$,
together with the relation
\beql{tmoo}
t_{-1,1}(u)\ts t_{kk}(v)\equiv 0.
\eeq
This relation is immediate from \eqref{defrel} for $k>1$ and $k<-1$;
for the latter we apply \eqref{defrel}
to the commutator $[t_{kk}(v),t_{-1,1}(u)]$. If $k=1$ or $k=-1$
then the claim is verified by a calculation similar
to the cases \eqref{tkkkv} and \eqref{tkkkvmo}, respectively.

If $\agot=\oa_{2n}$, then it is sufficient to verify
\eqref{tiipotkk} for $i\in\{-n,\dots,-2,1,\dots,n-1\}$ and all $k$,
together with the relations \eqref{tmoo} and
$
t_{-1,2}(u)\ts t_{kk}(v)\equiv 0.
$
The calculation is again a repetition
of the one for $\agot=\oa_{2n+1}$.

Now we verify that all the operators $t_{i\tss i}^{(r)}$ on the space $V^0$
with $i\in\{-n,\dots,n\}$ and $r\geqslant 1$ comprise a commutative family.
First of all, by \eqref{defrel} we have
$
[t_{ii}(u),t_{ii}(v)]=0
$
for any $i\ne 0$. Furthermore,
for any $i<j$ such that $i+j\ne 0$ we have
\ben
(u-v)\ts [t_{ii}(u),t_{jj}(v)]=t_{ji}(u)\ts t_{ij}(v)-t_{ji}(v)\ts t_{ij}(u)
\een
and so, $[t_{ii}(u),t_{jj}(v)]\equiv 0$ as operators on $V^0$.
Next, for any $0\leqslant i\leqslant j$ set
\ben
A_{ij}=t_{-j,-i}(u)\ts t_{ji}(v)-t_{ij}(v)\ts t_{-i,-j}(u),
\een
where the value $i=0$ only occurs in the case $\agot=\oa_{2n+1}$.
By \eqref{defrel}, we get
\beql{aooequiv}
A_{00}\equiv \frac{1}{u-v}\ts A_{00}
-\frac{1}{u-v-\kappa}\ts\sum_{j=0}^n \ts A_{0j},
\eeq
and for any $i>0$
\beql{aiiequiv}
A_{ii}\equiv -\frac{1}{u-v-\kappa}\ts\sum_{j=i}^n \ts A_{ij},
\eeq
as operators
on $V^0$, while for $0\leqslant i<j$ we have
\ben
A_{ij}\equiv -\frac{1}{u-v-\kappa}\ts\sum_{k=i}^n \ts A_{ik}
-\frac{1}{u-v-\kappa}\ts\sum_{l=j}^n \ts A_{jl}.
\een
This implies
\ben
A_{ij}\equiv A_{ii}-A_{jj}
\een
for $0<i<j$, and
\ben
A_{0j}=\frac{u-v-1}{u-v}\ts A_{00}-A_{jj}
\een
for $j>0$.
Hence, \eqref{aiiequiv} gives
\ben
\Big(\ts 1+\frac{n-i+1}{u-v-\kappa}\ts\Big) A_{ii}-\frac{1}{u-v-\kappa}\ts
\sum_{j=i+1}^n \ts A_{jj}\equiv 0,
\een
thus proving that $A_{ii}=[t_{-i,-i}(u),t_{ii}(v)]\equiv 0$ for all $i>0$
by an obvious induction.
Moreover, in the case $\agot=\oa_{2n+1}$, we derive
from \eqref{aooequiv} that $A_{00}=[t_{00}(u),t_{00}(v)]\equiv 0$.

Since the operators $t_{i\tss i}^{(r)}$ on $V^{0}$ are pairwise
commuting, they have a simultaneous eigenvector
$\xi\in V^{0}$.
Then $\xi$ satisfies the conditions \eqref{trianb}. Moreover,
since $V$ is irreducible, the submodule $\X(\agot)\ts\xi$ must coincide with $V$
so that $V$ is a highest weight module over $\X(\agot)$.
In particular,
$\xi$ is an $\agot$-weight vector with a certain weight $\nu$.

Finally, since the central elements $z_r$ act on $V$ as scalar
operators, Corollary~\ref{cor:pbwx} implies that
the vector space $V$ is spanned by the elements
\ben
t_{j_1i_1}^{(r_1)}\dots t_{j_mi_m}^{(r_m)}\ts \xi,\qquad m\geqslant 0,
\een
with $j_a>i_a$ and $r_a\geqslant 1$.
Hence, by \eqref{tijF}
the $\agot$-weight space $V_{\nu}$
is one-dimensional and spanned by the vector $\xi$.
Moreover, if $\rho$ is a weight of $V$ and $\rho\ne\nu$ then
$\rho$ strictly precedes $\nu$.
This proves that the highest vector $\xi$ of $V$ is determined
uniquely, up to a constant factor.
\epf

Given any tuple $\la(u)=(\la_{-n}(u),\dots,\la_n(u))$
of formal series of the form \eqref{laiu}, we define
the {\it Verma module\/} $M(\la(u))$ as the quotient of $\X(\agot)$ by
the left ideal generated by all the coefficients of the series $t_{ij}(u)$
with $-n\leqslant i<j\leqslant n$, and $t_{ii}(u)-\la_i(u)$ for
$i=-n,\dots,n$. As we shall see below,
the Verma module
$M(\la(u))$ can be trivial for some $\la(u)$.
In the non-trivial case,
the Verma module $M(\la(u))$ is a highest weight representation
of $\X(\agot)$ with the highest weight $\la(u)$ and
the highest vector $1_{\la}$ which is the canonical
image of the element $1\in \X(\agot)$.
Moreover, any highest weight representation
of $\X(\agot)$ with the highest weight $\la(u)$
is isomorphic to a quotient
of $M(\la(u))$.
Regarding
$M(\la(u))$ as an $\agot$-module,
we obtain the weight space decomposition
\ben
M(\la(u))=\underset{\nu}{\bigoplus}\ts M(\la(u))_{\nu},
\een
summed over all $\agot$-weights $\nu=(\nu_1,\dots,\nu_n)$
of $M(\la(u))$, where
\ben
M(\la(u))_{\nu}=
\{\eta\in M(\la(u))\ |\ F_{ii}\ts\eta=\nu_i\ts\eta,\quad i=1,\dots,n\}.
\een
By \eqref{tijF}, the set of weights of $M(\la(u))$ coincides with that of
the $\agot$-Verma module
with the highest weight $\la^{(1)}=(\la^{(1)}_1,\dots,\la^{(1)}_n)$.
This set consists of all weights of the form
$\la^{(1)}-\om$, where $\om$ is a $\ZZ_+$-linear combination
of the positive roots.

One easily shows that any submodule $K$ of
a non-trivial Verma module $M(\la(u))$
admits the weight space decomposition
\ben
K=\underset{\nu}{\bigoplus}\ts K_{\nu},\qquad K_{\nu}=K\cap M(\la(u))_{\nu}.
\een
This implies that
the sum of all proper submodules
is the unique maximal proper submodule
of $M(\la(u))$. The {\it irreducible highest weight representation\/}
$L(\la(u))$ of $\X(\agot)$ with the highest weight $\la(u)$ is defined as the
quotient of the Verma module $M(\la(u))$
by the unique maximal proper submodule.

\bpr\label{prop:actzu}
Let $V$ be a highest weight representation of $\X(\agot)$ with
the highest weight $\la(u)=(\la_{-n}(u),\dots,\la_n(u))$
with some series \eqref{laiu}. Then the coefficients
of the series $z(u)$ act on $V$ as scalar operators determined by
$z(u)|_{V}=\la_{-n}(u+\kappa)\ts\la_n(u)$.
\epr

\bpf
Let $\xi$ be the highest vector of $V$. Then $V=\X(\agot)\ts\xi$ so that
$z(u)$ acts on $V$ as a scalar function determined by its action on $\xi$.
However, taking $k=l=n$ in \eqref{zt} we get
\beql{zun}
z(u)=\sum_{i=-n}^n\theta_{ni}\ts t_{-i,-n}(u+\kappa)\ts t_{in}(u).
\eeq
Therefore, $z(u)\ts\xi=\la_{-n}(u+\kappa)\ts\la_n(u)\ts\xi$.
\epf

\subsection{Representations of low rank Yangians}\label{subsec:rlr}

Using the results on representations of the Yangian $\Y(\gl_2)$
(Tarasov~\cite{t:sq, t:im}; see also \cite[Chapter~12]{cp:gq}, \cite{m:fd}),
and the isomorphisms constructed in Section~\ref{sec:isom}, we describe here
the finite-dimensional irreducible representations
of the extended Yangians $\X(\oa_3)$, $\X(\spa_2)$ and $\X(\oa_4)$.
For the sake of completeness, we also reproduce a simplified version of
Tarasov's classification theorem for the representations of $\Y(\gl_2)$;
cf. \cite{m:fd}.

We shall use the notation for the generators of $\Y(\gl_2)$ introduced
in Section~\ref{sec:isom}.
A representation $L$ of the Yangian $\Y(\gl_2)$
is called a {\it highest weight representation\/} if there exists a nonzero vector
$\ze\in L$ such that $L$ is generated by $\ze$
and the following relations hold
\begin{alignat}{2}
\label{rtthwr}
T_{12}(u)\ts\ze&=0 \qquad &&\text{and}\\
\label{rttiihwr}
T_{ii}(u)\ts\ze&=\mu_i(u)\ts\ze \qquad &&\text{for} \quad i=1,2.
\end{alignat}
for some formal series
\beql{muaiu}
\mu_i(u)=1+\mu_i^{(1)}u^{-1}+\mu_i^{(2)}u^{-2}+\dots,\qquad\mu_i^{(r)}\in\CC.
\eeq
The vector $\ze$ is called the {\it highest vector\/}
of $L$, and the pair
$\mu(u)=\big(\mu_1(u),\mu_2(u)\big)$
is the {\it highest weight\/} of $L$. A standard argument,
similar to the one used in Section~\ref{subsec:hw} (see e.g. \cite{m:fd}),
shows that every finite-dimensional irreducible representation
of $\Y(\gl_2)$ is a highest weight representation.
Given any pair of series $\mu(u)=\big(\mu_1(u),\mu_2(u)\big)$,
the corresponding Verma module $M(\mu(u))$ for $\Y(\gl_2)$ is the quotient of
$\Y(\gl_2)$ by the left ideal generated by all the coefficients
of the series
$T_{12}(u)$ and $T_{ii}(u)-\mu_i(u)$ for
$i=1,2$. When the components of $\mu(u)$ satisfy
the condition $\mu_1(u)\ts\mu_2(u-1)=1$ then
$M(\mu(u))$ may also be regarded as a module over the Yangian $\Y(\sll_2)$.

The $\Y(\gl_2)$-module $M(\mu(u))$ has a unique irreducible quotient
$L(\mu(u))$. Thus, any finite-dimensional irreducible representation
of $\Y(\gl_2)$ is isomorphic to $L(\mu(u))$ for a pair
$\mu(u)=\big(\mu_1(u),\mu_2(u)\big)$. It remains to describe the highest
weights $\mu(u)$ which correspond to finite-dimensional modules $L(\mu(u))$.
This is given by the following theorem due to Tarasov~\cite{t:sq, t:im}
in Drinfeld's version \cite{d:nr}.

\bth\label{thm:classgltwo}
The irreducible highest weight representation
$L(\mu(u))$
of $\Y(\gl_2)$ is
finite-dimensional if and only if there exists a monic
polynomial $P(u)$ in $u$ such that
\beql{dpone}
\frac{\mu_1(u)}{\mu_2(u)}=\frac{P(u+1)}{P(u)}.
\eeq
In this case, $P(u)$ is unique.
\eth

\bpf
We shall need the following lemma.

\ble\label{lem:polyn}
If $\dim L(\mu(u))<\infty$ then
there exists a formal series
\ben
f(u)=1+f_1u^{-1}+f_2u^{-2}+\dots,\qquad f_r\in\CC,
\een
such that
$f(u)\mu_1(u)$ and
$f(u)\mu_2(u)$ are polynomials in $u^{-1}$.
\ele

\bpf
By twisting the action of $\Y(\gl_2)$ on $L(\mu(u))$
by the automorphism \eqref{multyang}
with $g(u)=\mu_2(u)^{-1}$,
we obtain a module over $\Y(\gl_2)$ which is isomorphic to
the irreducible highest weight representation
$L(\nu(u),1)$ with $\nu(u)=\mu_1(u)/\mu_2(u)$.
So, we may assume without loss of generality that the highest
weight of $L(\mu(u))$ has the form $\mu(u)=(\nu(u),1)$.
Let $\ze$ denote the highest vector of
the Verma module $M(\nu(u),1)$.
Since $\dim L(\nu(u),1)<\infty$,
the vectors $T_{21}^{(i)}\ze\in M(\nu(u),1)$ with $i\geqslant 1$ are
linearly dependent modulo the maximal proper
submodule $K$ of $M(\nu(u),1)$.
Hence, $M(\nu(u),1)$
contains a nonzero vector $\xi\in K$ of the form
\ben
\xi=\sum_{i=1}^m c_i\ts T_{21}^{(i)}\ze,\qquad c_i\in\CC.
\een
Here $m$ is a positive integer and we may assume that $c_m\ne 0$.
Then we have $T_{12}^{(r)}\xi=0$ for all $r\geqslant 1$ because
otherwise the highest vector $\ze$ would belong to $K$. Write
\ben
\nu(u)=1+\nu^{(1)}u^{-1}+\nu^{(2)}u^{-2}+\dots,\qquad \nu^{(i)}\in\CC.
\een
By the defining relations \eqref{defequiv}, in $M(\nu(u),1)$ we have
\ben
T_{12}^{(r)}\ts T_{21}^{(i)}\ze=
\sum_{a=1}^{\min(r,i)}
\Big(T^{(a-1)}_{22} T^{(r+i-a)}_{11}-T^{(r+i-a)}_{22} T^{(a-1)}_{11}\Big)\ze
=\nu^{(r+i-1)}\tss \ze.
\een
Hence, for all $r\geqslant 1$ we have the relations
\ben
\sum_{i=1}^m c_i\tss \nu^{(r+i-1)}=0.
\een
They imply
\ben
\nu(u)\tss(c_1+c_2\tss u+\dots+c_m\tss u^{m-1})
=(b_1+b_2\tss u+\dots+b_m\tss u^{m-1})
\een
for some coefficients $b_i\in\CC$ with $b_m=c_m$. Thus,
taking now
\ben
f(u)=c_m^{-1}\sum_{i=1}^m c_i\ts u^{-m+i}
\een
we conclude that both $f(u)\tss\nu(u)$ and $f(u)\tss 1$ are polynomials
in $u^{-1}$.
\epf

Thus, taking the composition of the
representation of $\Y(\gl_2)$ on $L(\mu(u))$ with an appropriate automorphism
of the form \eqref{multyang}, we can get another highest weight representation
of $\Y(\gl_2)$ where both components of the highest weight
are polynomials in $u^{-1}$.

For any $\al,\be\in\CC$ consider the irreducible highest weight
representation $L(\al,\be)$ of the Lie algebra $\gl_2$
and equip it with a $\Y(\gl_2)$-module structure
via the evaluation homomorphism \eqref{eval}.
Let $\ze$ denote the highest vector of
$L(\al,\be)$. Then
\ben
E_{11}\ts\ze=\al\ts\ze,\qquad
E_{22}\ts\ze=\be\ts\ze,\qquad
E_{12}\ts\ze=0.
\een
Moreover,
if $\al-\be\in\ZZ_+$ then
the vectors $(E_{21})^r\ze$ with
$r=0,1,\dots,\al-\be$ form a basis of $L(\al,\be)$
so that $\dim L(\al,\be)=\al-\be+1$. If
$\al-\be\notin\ZZ_+$ then a basis of $L(\al,\be)$
is formed by the vectors $(E_{21})^r\ze$, where
$r$ runs over all nonnegative integers.

Now let $\mu_1(u)$ and $\mu_2(u)$ be polynomials in $u^{-1}$
of degree not more than $k$. Write the decompositions
\beql{laalbedec}
\bal
\mu_1(u)&=(1+\al_1u^{-1})\dots (1+\al_ku^{-1}),
\\
\mu_2(u)&=(1+\be_1u^{-1})\dots (1+\be_ku^{-1}),
\eal
\eeq
where the constants $\al_i$ and $\be_i$ are complex numbers
(some of them are zero if the degree of the corresponding polynomial
is strictly less than $k$).

For any $\Y(\gl_2)$-modules $L_1$ and $L_2$, their tensor product
$L_1\ot L_2$ is equipped with a $\Y(\gl_2)$-module structure
defined by the coproduct \eqref{Deltagl}.

\ble\label{lem:decomptp}
Re-number the parameters $\al_i$ and $\be_i$ if necessary, so that
for every index $i=1,\dots, k-1$ the following condition holds:
if the multiset $\{\al_p-\be_q\ |\ i\leqslant p,q \leqslant k\}$
contains nonnegative integers, then $\al_i-\be_i$
is minimal amongst them. Then the representation $L(\mu_1(u),\mu_2(u))$
of $\ts\Y(\gl_2)$ is isomorphic to the tensor product module
\beql{tpll}
L(\al_1,\be_1)\ot L(\al_2,\be_2)\ot\dots\ot L(\al_k,\be_k).
\eeq
\ele

\bpf
Let us denote the module \eqref{tpll} by $L$ and
let $\ze_i$ be the highest vector of
$L(\al_i,\be_i)$ for $i=1,\dots,k$.
Using the definition of the coproduct on $\Y(\gl_2)$ we derive
that the cyclic span $\Y(\gl_2)\tss\ze$ of
the vector
$
\ze=\ze_1\ot\dots\ot\ze_k
$
is a highest weight module
with the highest weight $(\mu_1(u),\mu_2(u))$. Therefore,
the proposition will follow if we prove that
the module $L$ is irreducible.

We claim that any vector $\xi\in L$
satisfying $T_{12}(u)\tss\xi=0$ is proportional to $\ze$.
We shall prove this claim
by induction on $k$.
This is obvious for $k=1$ so suppose that $k\geqslant 2$.
Write any such vector $\xi$, which is assumed
to be nonzero, in the form
\ben
\xi=\sum_{r=0}^p(E_{21})^r \ze_1\ot \xi_r
\quad
\text{where}
\quad
\xi_r\in L(\al_2,\be_2)\ot\dots\ot L(\al_k,\be_k)
\een
and $p$ is some non-negative integer. Moreover,
if $\al_1-\be_1\in\ZZ_+$ then we may and will assume that
$p\leqslant\al_1-\be_1$.
We also assume that $\xi_p\ne 0$.
Applying $T_{12}(u)$ to $\xi$, with the use of \eqref{Deltagl} we get
\beql{srpte}
\sum_{r=0}^p\Big(T_{11}(u)(E_{21})^r \ze_1\ot T_{12}(u)\tss \xi_r
+T_{12}(u)(E_{21})^r \ze_1\ot T_{22}(u)\tss \xi_r\Big)=0.
\eeq
Using the definition
of the Yangian action on $L(\al_1,\be_1)$
and commutation relations in $\gl_2$, we obtain
\ben
T_{11}(u)(E_{21})^r \ze_1=(1+E_{11}u^{-1})(E_{21})^r \ze_1
=(1+(\al_1-r)\tss u^{-1})(E_{21})^r \ze_1,
\een
and
\ben
T_{12}(u)(E_{21})^r \ze_1=u^{-1}E_{12}\tss(E_{21})^r \ze_1=
u^{-1}\ts r(\al_1-\be_1-r+1)
(E_{21})^{r-1} \ze_1.
\een
Hence, taking the coefficient at $(E_{21})^p \ze_1$ in \eqref{srpte} gives
\ben
(1+(\al_1-p)\tss u^{-1})\ts T_{12}(u)\tss\xi_p=0,
\een
implying the relation $T_{12}(u)\tss\xi_p=0$. By the induction hypothesis,
applied to the $\Y(\gl_2)$-module
$
L(\al_2,\be_2)\ot\dots\ot L(\al_k,\be_k),
$
the vector
$\xi_p$ must be proportional to $\ze_2\ot\dots\ot\ze_k\tss$.
Therefore, using \eqref{Deltagl} we get
\beql{tuetap}
T_{22}(u)\tss\xi_p=(1+\be_2u^{-1})\dots (1+\be_ku^{-1}) \tss\xi_p.
\eeq
In order to complete the proof of the claim it now suffices to show
that $p$ must be equal to zero.
Suppose by way of contradiction that $p\geqslant 1$. Then
taking the coefficient at
$(E_{21})^{p-1} \ze_1$ in \eqref{srpte} we derive
\ben
(1+(\al_1-p+1)u^{-1})\ts T_{12}(u)\tss\xi_{p-1}+
u^{-1}\ts p\tss(\al_1-\be_1-p+1)\ts T_{22}(u)\tss\xi_p=0.
\een
Hence, multiplying by $u^k$ and taking into account \eqref{tuetap}
we get
\ben
(u+\al_1-p+1)u^{k-1}T_{12}(u)\tss\xi_{p-1}+
\ts p\tss(\al_1-\be_1-p+1)(u+\be_2)\dots (u+\be_k)\tss\xi_p=0.
\een
Now observe that the vector
$u^{k-1}T_{12}(u)\tss\xi_{p-1}$
depends on $u$ polynomially. This follows by an easy induction
with the use of \eqref{Deltagl}.
So,
taking the value $u=-\al_1+p-1$ we obtain the relation
\ben
p\tss(\al_1-\be_1-p+1)(\al_1-\be_2-p+1)
\dots (\al_1-\be_k-p+1)=0.
\een
But this is impossible due to the conditions on the parameters
$\al_i$ and $\be_i$. Thus, $p$ must be zero and the claim follows.

Suppose now that $M$ is a nonzero submodule of $L$.
Then $M$ must contain a nonzero vector $\xi$ such that
$T_{12}(u)\tss\xi=0$. Indeed, this follows from
the fact that
the set of $\gl_2$-weights of $L$ has
an upper boundary. The above argument thus shows that $M$
contains the vector $\ze$. It remains to prove that
the cyclic span $K=\Y(\gl_2)\tss\ze$ coincides with $L$.

Denote by $\varkappa$ the anti-automorphism of the algebra $\Y(\gl_2)$,
defined by
\beql{varsigma}
\varkappa:\ t_{ij}(u)\mapsto t_{3-i,3-j}(-u).
\eeq
Consider the vector space $L^*$ dual to $L$. That is,
$L^*$ is spanned by all linear maps $\sigma:L\to \CC$
satisfying the condition that the linear span of the vectors $\eta\in L$
such that $\sigma(\eta)\ne 0$, is finite-dimensional.
Equip $L^*$ with a $\Y(\gl_2)$-module structure
by setting
\ben
(y\ts\sigma)(\eta)=\sigma(\varkappa(y)\ts\eta)
\quad
\text{for}
\quad
y\in\Y(\gl_2)
\fand
\sigma\in L^*,\ \eta\in L.
\een
It is easy to see that
the dual module $L(\al,\be)^*$ to the evaluation
module $L(\al,\be)$ is isomorphic to
$L(-\be,-\al)$. Moreover,
the $\Y(\gl_2)$-module $L^*$ is isomorphic to the tensor
product module
\ben
L(-\be_1,-\al_1)\ot\dots\ot L(-\be_k,-\al_k).
\een
This is deduced from the fact that
the anti-automorphism $\varkappa$ commutes with the
coproduct $\Delta$, where $\varkappa$ is extended to
$\Y(\gl_2)\ot \Y(\gl_2)$ by
$\varkappa(x\ot y)=\varkappa(x)\ot\varkappa(y)$
for $x,y\in\Y(\gl_2)$.
Furthermore, the highest vector $\ze_i^*$ of the module
$L(-\be_i,-\al_i)\cong L(\al_i,\be_i)^*$ can be identified with
the element of $L(\al_i,\be_i)^*$ such that
$\ze_i^*(\ze_i)=1$
and $\ze_i^*(\eta_i)=0$ for all weight vectors $\eta_i\in L(\al_i,\be_i)$
whose weights are different from $(\al_i,\be_i)$.

Suppose now that the submodule $K$ of $L$ is proper and consider
its annihilator
\ben
\text{Ann\ts} K=\{\xi^*\in L^*\ |\ \xi^*(\eta)=0
\quad\text{for all}\quad\eta\in K\}.
\een
Then $\text{Ann\ts} K$
is a nonzero submodule of $L^*$, which does not
contain the vector $\ze_1^*\ot\dots\ot\ze_k^*$.
However, this contradicts the claim verified in the first part of the proof,
because the condition on the parameters $\al_i$ and $\be_i$
remain satisfied after we replace each $\al_i$ by $-\be_i$ and each
$\be_i$ by $-\al_i$.
\epf

By this lemma, all differences $\al_i-\be_i$ must be
nonnegative integers because the representation $L(\la_1(u),\la_2(u))$
is finite-dimensional. Then the polynomial
\beql{dptpepr}
P(u)=\prod_{i=1}^k (u+\be_i)(u+\be_i+1)\dots (u+\al_i-1)
\eeq
obviously satisfies \eqref{dpone}.

Conversely, suppose \eqref{dpone} holds for
a polynomial $P(u)=(u+\ga_1)\dots(u+\ga_p)$.
Set
\ben
\bal
\nu_1(u) &=(1+(\ga_1+1)u^{-1})\dots (1+(\ga_p+1)u^{-1}),\\
\nu_2(u) &=(1+\ga_1 u^{-1})\dots (1+\ga_p u^{-1}),
\eal
\een
and consider the tensor product module
\ben
L=L(\ga_1+1,\ga_1)\ot L(\ga_2+1,\ga_2)\ot\dots\ot L(\ga_p+1,\ga_p)
\een
of $\Y(\gl_2)$. Obviously, this module is finite-dimensional.
The cyclic $\Y(\gl_2)$-span
of the tensor product of the highest vectors of $L(\ga_i+1,\ga_i)$
is a highest weight module with the highest weight $(\nu_1(u),\nu_2(u))$.
Since this submodule is finite-dimensional, then so is its
irreducible quotient $L(\nu_1(u),\nu_2(u))$. Since
\ben
\frac{\nu_1(u)}{\nu_2(u)}=\frac{\mu_1(u)}{\mu_2(u)},
\een
there exists an automorphism of $\Y(\gl_2)$ of the form
\eqref{multyang}
such that its composition with the representation
$L(\nu_1(u),\nu_2(u))$ is isomorphic to
$L(\mu_1(u),\mu_2(u))$. Thus, the latter
is also finite-dimensional.

Finally, suppose that $Q(u)$ is another
monic polynomial in $u$ and
\ben
\frac{P(u+1)}{P(u)}=\frac{Q(u+1)}{Q(u)}.
\een
This means that the ratio $P(u)/Q(u)$ is periodic in $u$
which is only possible for $P(u)=Q(u)$.
\epf

The polynomial $P(u)$ is called the {\it Drinfeld polynomial\/}
of the representation $L(\mu(u))$.

We now apply Theorem~\ref{thm:classgltwo} to the low rank
extended Yangians.

\bco\label{cor:fdc1}
The Verma module $M(\la(u))$ over $\X(\spa_2)$ is non-trivial for
any highest weight $\la(u)=(\la_{-1}(u),\la_1(u))$.
Moreover, the $\X(\spa_2)$-module $L(\la(u))$ is finite-dimensional
if and only if there exists a monic
polynomial $P(u)$ in $u$ such that
\beql{dponec1}
\frac{\la_{-1}(u)}{\la_1(u)}=\frac{P(u+2)}{P(u)}.
\eeq
In this case, $P(u)$ is unique.
\eco

\bpf
This is immediate from Proposition~\ref{prop:c1}
and Theorem~\ref{thm:classgltwo}.
\epf

The evaluation homomorphism provided by Corollary~\ref{cor:evc1}
allows one to regard any irreducible $\spa_2$-module $V(\mu)$
as an $\X(\spa_2)$-module. The corresponding evaluation module
is immediately identified with an irreducible highest weight module.

\bpr\label{prop:evalc1}
The evaluation module $V(\mu)$ over $\X(\spa_2)$ is isomorphic
to $L(\la(u))$ with
\ben
\la_{-1}(u)=1-\mu_1\ts u^{-1}\Fand \la_{1}(u)=1+\mu_1\ts u^{-1}.
\vspace{-1.3\baselineskip}
\een
\qed
\epr

\bco\label{cor:fdb1}
The Verma module $M(\la(u))$ over $\X(\oa_3)$ is non-trivial
if and only if
the highest weight $\la(u)=\big(\la_{-1}(u),\la_0(u),\la_1(u)\big)$
satisfies the condition
\beql{condexi}
\la_{-1}(u-1/2)\ts\la_1(u)=\la_{0}(u-1/2)\ts\la_0(u).
\eeq
Moreover, if this condition holds then
the $\X(\oa_3)$-module $L(\la(u))$ is finite-dimensional
if and only if there exists a monic
polynomial $P(u)$ in $u$ such that
\beql{dponeb1}
\frac{\la_{0}(u)}{\la_1(u)}=
\frac{P(u+1/2)}{P(u)}.
\eeq
In this case, $P(u)$ is unique.
\eco

\bpf
Let the Verma module $M(\la(u))$ be non-trivial.
By Proposition~\ref{prop:b1},
we may regard $M(\la(u))$ as a $\Y(\gl_2)$-module.
In particular,
we have
\ben
T_{11}(2u)\ts T_{11}(2u+1)\ts 1_{\la}=\la_{-1}(u)\ts 1_{\la},
\een
where $1_{\la}$ is the highest vector of
$M(\la(u))$. This implies that $1_{\la}$ is an eigenvector
for $T_{11}(u)$, that is, $T_{11}(u)\ts 1_{\la}=\mu_{1}(u)\ts 1_{\la}$
for a certain series $\mu_1(u)$. Moreover, this series satisfies
\beql{muone}
\mu_1(2u)\ts \mu_1(2u+1)=\la_{-1}(u).
\eeq
Similarly, $T_{22}(u)\ts 1_{\la}=\mu_{2}(u)\ts 1_{\la}$
for a series $\mu_2(u)$ satisfying
\beql{mutwo}
\mu_2(2u)\ts \mu_2(2u+1)=\la_{1}(u).
\eeq
Furthermore, by the defining relations \eqref{defreltij}
we have
\ben
T_{12}(2u)\ts T_{22}(2u+1)+T_{22}(2u)\ts T_{12}(2u+1)
=2\ts T_{12}(2u+1)\ts T_{22}(2u).
\een
Since $t_{0,1}(u)\ts 1_{\la}=0$ we derive that $T_{12}(u)\ts 1_{\la}=0$.
Hence, using the action of $t_{0,0}(u)$ on $1_{\la}$ we also get
\beql{muonetwo}
\mu_1(2u)\ts \mu_2(2u+1)=\la_0(u).
\eeq
This gives the condition \eqref{condexi}.

Conversely, if the condition \eqref{condexi} holds for
a highest weight $\la(u)$ then there exist series
$\mu_1(u)$ and $\mu_2(u)$ satisfying \eqref{muone}, \eqref{mutwo}
and \eqref{muonetwo}. Consider the Verma module $M(\mu_1(u),\mu_2(u))$
over $\Y(\gl_2)$. Using the formulas of Proposition~\ref{prop:b1},
we find that the highest vector $1_{\mu}\in M(\mu_1(u),\mu_2(u))$
satisfies the conditions \eqref{trianb} for the action
of $\X(\oa_3)$.

The argument of the first part of the proof shows that,
regarded as a $\Y(\gl_2)$-module,
the module $L(\la(u))$ is isomorphic to $L(\mu_1(u),\mu_2(u))$
with $\mu_1(u)$ and $\mu_2(u)$ satisfying \eqref{muone}, \eqref{mutwo}
and \eqref{muonetwo}. So writing the relation of
Theorem~\ref{thm:classgltwo} in terms of the series $\la_i(u)$,
we get the desired condition.
\epf

The evaluation homomorphism provided by Corollary~\ref{cor:evb1}
allows one to regard any irreducible $\oa_3$-module $V(\mu)$
as an $\X(\oa_3)$-module.

\bpr\label{prop:evalb1}
The evaluation module $V(\mu)$ over $\X(\oa_3)$ is isomorphic
to $L(\la(u))$ with
\ben
\bal
\la_{-1}(u)&=\frac{(2u-\mu_1)(2u-\mu_1-1)}{2u\ts(2u-1)},\\
\la_0(u)&=\frac{(2u+\mu_1)(2u-\mu_1-1)}{2u\ts(2u-1)},\\
\la_{1}(u)&=\frac{(2u+\mu_1)(2u+\mu_1-1)}{2u\ts(2u-1)}.
\eal
\een
\epr

\bpf
This is immediate from Corollary~\ref{cor:evb1}, as the Casimir element
$c$ acts on $V(\mu)$ as multiplication by the scalar $(\mu_1^2-\mu_1)/2$.
\epf

\bco\label{cor:fdd2}
The Verma module $M(\la(u))$ over $\X(\oa_4)$ is non-trivial
if and only if
the highest weight $\la(u)=\big(\la_{-2}(u),\la_{-1}(u),\la_1(u),\la_2(u)\big)$
satisfies the condition
\beql{condexid2}
\la_{-2}(u)\ts\la_2(u)=\la_{-1}(u)\ts\la_1(u).
\eeq
Moreover, if this condition holds then
the $\X(\oa_4)$-module $L(\la(u))$ is finite-dimensional
if and only if there exist monic
polynomials $P(u)$ and $Q(u)$ in $u$ such that
\beql{dponed2}
\frac{\la_{-1}(u)}{\la_2(u)}=
\frac{P(u+1)}{P(u)}
\Fand
\frac{\la_{1}(u)}{\la_2(u)}=
\frac{Q(u+1)}{Q(u)}.
\eeq
In this case, $P(u)$ and $Q(u)$ are determined uniquely.
\eco

\bpf
Suppose that the Verma module $M(\la(u))$ over $\X(\oa_4)$
is non-trivial. Using the isomorphism
$\chi^{(1)}$
provided by Corollary~\ref{cor:restdisom},
we shall regard
$M(\la(u))$ as a module over the algebra $\Y(\sll_2)\ot\Y(\gl_2)$.
As was seen in the proof of Corollary~\ref{cor:restdisom},
\ben
\chi^{(1)}:t_{-2,-2}(u)\ts t_{1,1}(u-1)-t_{1,-2}(u)\ts t_{-2,1}(u-1)
\mapsto T^{\ts\prime}_{11}(u)\ts T^{\ts\prime}_{11}(u-1).
\een
This implies that $1_{\la}$ is an eigenvector
for $T^{\ts\prime}_{11}(u)$, that is,
$T^{\ts\prime}_{11}(u)\ts 1_{\la}=\mu'_{1}(u)\ts 1_{\la}$
for a certain series $\mu'_1(u)$. Similarly,
$T^{\ts\prime}_{22}(u)\ts 1_{\la}=\mu'_{2}(u)\ts 1_{\la}$
for a series $\mu'_2(u)$. Then, by the formulas of
Proposition~\ref{prop:d2}, we also have
\ben
\mathcal T_{11}(u)\ts 1_{\la}=\mu_{1}(u)\ts 1_{\la}\Fand
\mathcal T_{22}(u)\ts 1_{\la}=\mu_{2}(u)\ts 1_{\la}
\een
for some series $\mu_1(u)$ and $\mu_2(u)$. Moreover, we have the relations
\beql{mucondla}
\bal
\la_{-2}(u)&=\mu_{1}(u)\ts\mu'_{1}(u),\qquad&
\la_{-1}(u)&=\mu_{1}(u)\ts\mu'_{2}(u),\\
\la_{1}(u)&=\mu_{2}(u)\ts\mu'_{1}(u),\qquad&
\la_{2}(u)&=\mu_{2}(u)\ts\mu'_{2}(u),
\eal
\eeq
which imply \eqref{condexid2}. Conversely, if \eqref{condexid2} holds
for some series $\la_i(u)$,
then there exist series $\mu_i(u)$ and $\mu'_i(u)$ satisfying
\eqref{mucondla} together with the condition $\mu_1(u)\ts\mu_2(u-1)=1$.
Consider the $\Y(\sll_2)\ot\Y(\gl_2)$-module
$M(\mu_1(u),\mu_2(u))\ot M(\mu'_1(u),\mu'_2(u))$. The vector
$1_{\mu}\ot 1_{\mu'}$ satisfies the conditions \eqref{trianb}
for the action of the series $t_{ij}(u)$ thus proving that
the $X(\oa_4)$-module $M(\la(u))$ is non-trivial.

Finally, the argument of the first part of the proof shows that,
regarded as a $\Y(\sll_2)\ot\Y(\gl_2)$-module,
the module $L(\la(u))$ is isomorphic to
$L(\mu_1(u),\mu_2(u))\ot L(\mu'_1(u),\mu'_2(u))$
with the $\mu_i(u)$ and $\mu'_i(u)$ satisfying \eqref{mucondla}.
By Theorem~\ref{thm:classgltwo}, the module
$L(\mu_1(u),\mu_2(u))\ot L(\mu'_1(u),\mu'_2(u))$ is finite-dimensional
if and only if there exist monic polynomials $P(u)$ and $Q(u)$ in $u$
such that
\ben
\frac{\mu_{1}(u)}{\mu_2(u)}=
\frac{P(u+1)}{P(u)}\Fand
\frac{\mu'_{1}(u)}{\mu'_2(u)}=
\frac{Q(u+1)}{Q(u)}.
\een
Writing these formulas in terms of the $\la_i(u)$ we get the desired
conditions.
\epf

The evaluation homomorphism provided by Corollary~\ref{cor:evd2}
allows one to regard any irreducible $\oa_4$-module $V(\mu)$
as an $\X(\oa_4)$-module.

\bpr\label{prop:evald2}
The evaluation module $V(\mu)$ over $\X(\oa_4)$ is isomorphic
to $L(\la(u))$ with
\ben
\bal
\la_{-2}(u)&=\frac{(2u-\mu_1-\mu_2)(2u+\mu_1-\mu_2)}{4u^2},\\
\la_{-1}(u)&=\frac{(2u-\mu_1-\mu_2)(2u-\mu_1+\mu_2)}{4u^2},\\
\la_{1}(u)&=\frac{(2u+\mu_1-\mu_2)(2u+\mu_1+\mu_2)}{4u^2},\\
\la_{2}(u)&=\frac{(2u-\mu_1+\mu_2)(2u+\mu_1+\mu_2)}{4u^2}.
\eal
\een
\epr

\bpf
This follows from Corollary~\ref{cor:evd2}, as the Casimir element
$c$ acts on $V(\mu)$ as multiplication by the scalar $(\mu_1^2+\mu_2^2)/2-\mu_2$.
\epf

\bre\label{rem:evalmod}
More general evaluation modules $V(\mu)_{a}$ with $a\in\CC$ over $\X(\agot)$
for $\agot=\spa_2$, $\oa_3$ and $\oa_4$ can be obtained by using
the respective evaluation homomorphisms $\ev_a:\X(\agot)\to\U(\agot)$
instead of $\ev$; see Remark~\ref{rem:evala}. Then
$V(\mu)_{a}$ will be isomorphic to the irreducible highest weight module
$L(\la(u))$, where the components $\la_i(u)$ are found from
the formulas of Propositions~\ref{prop:evalc1}, \ref{prop:evalb1}
or \ref{prop:evald2} by replacing $u$ with $u-a$.
\ere

\subsection{Classification theorems}\label{subsec:clth}

Our goal here is to prove classification theorems for the
finite-dimensional irreducible representations of the extended
Yangians $\X(\agot)$ for $\agot=\oa_{2n+1}$, $\spa_{2n}$,
and $\oa_{2n}$. The corresponding theorem
for the Yangian $\Y(\sll_N)$
implies that every finite-dimensional irreducible representation
of $\Y(\sll_N)$ is isomorphic to a subquotient of the tensor
product of the fundamental representations \cite{d:nr},
\cite[Chapter~12]{cp:gq}. We shall use the following version of the
well-known construction of the fundamental representations
of $\Y(\sll_N)$. They are obtained by restriction from the corresponding
representation of $\Y(\gl_N)$ which is obtained by a simple particular
case of the fusion procedure; see e.g. \cite{c:ni}, \cite{n:yc}.
The vector space $\CC^N$ carries an irreducible representation
of $\Y(\gl_N)$ with the action of the generators given by
\ben
T_{ij}(u)\mapsto \de_{ij}+e_{ij}\ts u^{-1},\qquad i,j\in\{1,\dots, N\},
\een
where the $e_{ij}$ denote the standard matrix units. So,
\ben
T_{ij}(u)\ts e_k=\de_{ij}\ts e_k+\de_{jk}\ts e_i\ts u^{-1},
\een
where $e_1,\dots,e_N$ denote the canonical basis of $\CC^N$.
Since
for any $b\in\CC$ the mapping
$T_{ij}(u)\mapsto T_{ij}(u-b)$ defines an automorphism of $\Y(\gl_N)$,
using the coproduct \eqref{Deltagl}, we can equip the tensor product
$(\CC^N)^{\ot m}$ with the action of $\Y(\gl_N)$ by the rule
\begin{multline}\label{tijtpreigln}
T_{ij}(u)\ts (e_{i_1}\ot\cdots\ot e_{i_m})=\\
\sum_{a_1,\dots,a_{m-1}=1}^N T_{ia_1}(u)\ts e_{i_1}\ot T_{a_1a_2}(u+1)\ts e_{i_2}
\ot\cdots\ot T_{a_{m-1}j}(u+m-1)\ts e_{i_m}.
\end{multline}
For any $1\leqslant m<N$ set
\ben
\xi_m=\sum_{\si\in\Sym_m} \sgn\si\cdot e_{\si(1)}\ot\cdots\ot e_{\si(m)} \in
(\CC^N)^{\ot m}.
\een
Then $\xi_m$ has the properties
\beql{tijxigln}
T_{ij}(u)\ts \xi_m=0\qquad\text{for all}\quad 1\leqslant i<j\leqslant N
\eeq
and
\ben
T_{i\tss i}(u)\ts \xi_m=\begin{cases}
\dfrac{u+m}{u+m-1}\ts\xi_m \qquad&\text{if}\quad 1\leqslant i\leqslant m,\\[1em]
\xi_m\qquad&\text{if}\quad m+1\leqslant i\leqslant N.
\end{cases}
\een
Thus, the vector $\xi_m$ generates a highest weight module
over $\Y(\gl_N)$ whose irreducible quotient is isomorphic
to a fundamental module; see \cite[Chapter~12]{cp:gq}, \cite{m:fd}.

Consider the extended Yangian $\X(\agot')$ for the subalgebra
$\agot'$ of $\agot$ of rank $n-1$. That is,
\ben
\agot'=\oa_{2n-1},\  \spa_{2n-2},\
\oa_{2n-2} \quad\text{respectively for}\quad
\agot=\oa_{2n+1},\  \spa_{2n},\
\oa_{2n}.
\een
Note that $\X(\agot')$ is not
a natural subalgebra of $\X(\agot)$.
Let $V$ be an $\X(\agot)$-module.
Set
\ben
\bal
V^+=\{\eta\in V\ |\ {}&{}t_{k,n}(u)\ts\eta=0&\qquad&\text{for}\quad k<n
\Fand\\
{}&{}t_{-n,k}(u)\ts\eta=0&\qquad&\text{for}\quad k>-n\}.
\eal
\een

\ble\label{lem:vplus}
The subspace $V^+$ is stable under all operators $t_{ij}(u)$
with the condition $-n+1\leqslant i,j\leqslant n-1$. Moreover, these operators
form a representation of the algebra $\X(\agot')$ on $V^+$, where
each operator $t_{ij}(u)$ is the image of the generator series
of $\X(\agot')$ with the same name.
\ele

\bpf
For any $\eta\in V^+$ we have the following relations modulo
elements of $V^+$ which are implied by
\eqref{defrel}: if $k<n$ and $-n+1\leqslant i,j\leqslant n-1$ then
\ben
t_{kn}(v)\ts t_{ij}(u)\ts\eta\equiv -[t_{ij}(u),t_{kn}(v)]\ts\eta
\equiv \frac{\de_{k,-i}}{u-v-\kappa}\ts
\theta_{i,-n}t_{-n,j}(u)\ts t_{nn}(v)\ts\eta.
\een
However, applying again \eqref{defrel}, we find that
\ben
t_{-n,j}(u)\ts t_{nn}(v)\ts\eta\equiv
-\frac{1}{u-v-\kappa}\ts t_{-n,j}(u)\ts t_{nn}(v)\ts\eta.
\een
Therefore, $t_{-n,j}(u)\ts t_{nn}(v)\ts\eta\equiv 0$ implying
$t_{kn}(v)\ts t_{ij}(u)\ts\eta\equiv 0$. A similar calculation
shows that for any $k>-n$ and $-n+1\leqslant i,j\leqslant n-1$
we also have
$t_{-n,k}(v)\ts t_{ij}(u)\ts\eta\equiv 0$ proving the first part of
the lemma.

In order to prove the second part, suppose
that the indices $i,j,k,l$ satisfy the condition
$-n+1\leqslant i,j,k,l\leqslant n-1$. Then by \eqref{defrel} for any
$\eta\in V^+$ we have
\ben
\bal[]
[\tss t_{ij}(u),t_{kl}(v)]\ts\eta=\frac{1}{u-v}&
\Big(t_{kj}(u)\ts t_{il}(v)-t_{kj}(v)\ts t_{il}(u)\Big)\ts\eta\\
{}-\frac{1}{u-v-\kappa}
&\Big(\de_{k,-i}\sum_{p=-n}^n\theta_{ip}\ts t_{pj}(u)\ts t_{-p,l}(v)-
\de_{l,-j}\sum_{p=-n}^n\theta_{jp}\ts t_{k,-p}(v)\ts t_{ip}(u)\Big)\ts\eta.
\eal
\een
Writing the right hand side modulo $V^+$, we get
\ben
\bal[]
\frac{1}{u-v}&
\Big(t_{kj}(u)\ts t_{il}(v)-t_{kj}(v)\ts t_{il}(u)\Big)\ts\eta\\
{}-\frac{1}{u-v-\kappa}
&\Big(\de_{k,-i}\sum_{p=-n+1}^{n-1}\theta_{ip}\ts t_{pj}(u)\ts t_{-p,l}(v)-
\de_{l,-j}\sum_{p=-n+1}^{n-1}\theta_{jp}\ts t_{k,-p}(v)\ts t_{ip}(u)\Big)\ts\eta\\
-\frac{1}{u-v-\kappa}
&\Big(\de_{k,-i}\ts\theta_{i,-n}\ts t_{-n,j}(u)\ts t_{n,l}(v)-
\de_{l,-j}\ts\theta_{j,-n}\ts t_{k,n}(v)\ts t_{i,-n}(u)\Big)\ts\eta.
\eal
\een
Applying again \eqref{defrel}, we obtain
\ben
\bal
t_{-n,j}(u)\ts t_{n,l}(v)\ts\eta\equiv
-\frac{1}{u-v-\kappa}&
\sum_{p=-n+1}^{n-1}\theta_{-n,p}\ts t_{pj}(u)\ts t_{-p,l}(v)\ts\eta\\
-\frac{1}{u-v-\kappa}&
\Big(t_{-n,j}(u)\ts t_{n,l}(v)-
\de_{l,-j}\ts\theta_{j,-n}\ts t_{n,n}(v)\ts t_{-n,-n}(u)\Big)\ts\eta.
\eal
\een
Hence,
\ben
\bal
t_{-n,j}(u)\ts t_{n,l}(v)\ts\eta\equiv
-\frac{1}{u-v-\kappa+1}&
\sum_{p=-n+1}^{n-1}\theta_{-n,p}\ts t_{pj}(u)\ts t_{-p,l}(v)\ts\eta\\
+\frac{1}{u-v-\kappa+1}\ts{}&{}
\de_{l,-j}\ts\theta_{j,-n}\ts t_{n,n}(v)\ts t_{-n,-n}(u)\ts\eta.
\eal
\een
Similarly,
$t_{k,n}(v)\ts t_{i,-n}(u)\ts\eta\equiv -[t_{i,-n}(u),t_{k,n}(v)]\ts\eta$
and
\ben
\bal[]
[t_{i,-n}(u),t_{k,n}(v)]\ts\eta\equiv
-\frac{1}{u-v-\kappa}\ts{}&{}
\de_{k,-i}\ts\theta_{i,-n}\ts t_{-n,-n}(u)\ts t_{nn}(v)\ts\eta\\
{}+\frac{1}{u-v-\kappa}\ts{}
&{}\Big(\sum_{p=-n+1}^{n-1}\theta_{-n,p}\ts t_{k,-p}(v)\ts t_{ip}(u)
+t_{k,n}(v)\ts t_{i,-n}(u)\Big)\ts\eta
\eal
\een
which gives
\ben
\bal[]
t_{k,n}(v)\ts t_{i,-n}(u)\ts\eta\equiv
\frac{1}{u-v-\kappa+1}\ts{}&{}
\de_{k,-i}\ts\theta_{i,-n}\ts t_{-n,-n}(u)\ts t_{nn}(v)\ts\eta\\
{}-\frac{1}{u-v-\kappa+1}\ts{}
&{}\sum_{p=-n+1}^{n-1}\theta_{-n,p}\ts t_{k,-p}(v)\ts t_{ip}(u)\ts\eta.
\eal
\een
Combining these expressions, we come to the following relation
\ben
\bal[]
&[\tss t_{ij}(u),t_{kl}(v)]\ts\eta\equiv\frac{1}{u-v}
\Big(t_{kj}(u)\ts t_{il}(v)-t_{kj}(v)\ts t_{il}(u)\Big)\ts\eta\\
{}&{}-\frac{1}{u-v-\kappa+1}
\Big(\de_{k,-i}\sum_{p=-n+1}^{n-1}\theta_{ip}\ts t_{pj}(u)\ts t_{-p,l}(v)-
\de_{l,-j}\sum_{p=-n+1}^{n-1}\theta_{jp}\ts t_{k,-p}(v)\ts t_{ip}(u)\Big)\ts\eta\\
{}&{}+\frac{1}{(u-v-\kappa)(u-v-\kappa+1)}\ts
\de_{k,-i}\ts\de_{l,-j}\ts\theta_{ij}
\ts [t_{-n,-n}(u),t_{nn}(v)]\ts\eta.
\eal
\een
Finally, by \eqref{defrel},
\ben
[t_{-n,-n}(u),t_{nn}(v)]\ts\eta\equiv
-\frac{1}{u-v-\kappa}\ts [t_{-n,-n}(u),t_{nn}(v)]\ts\eta,
\een
so that $[t_{-n,-n}(u),t_{nn}(v)]\ts\eta\equiv 0$.
This yields the desired relations between the operators $t_{ij}(u)$
on $V^+$ since $\kappa-1=\kappa'$ coincides with the value
of the parameter $\kappa$ for the Lie algebra $\agot'$.
\epf

\bpr\label{prop:nontrvm}
The Verma module $M(\la(u))$ over $\X(\agot)$ is non-trivial
if and only if the components of the highest weight $\la(u)$ satisfy
the conditions
\beql{nontrvm}
\frac{\la_{-n+i-1}(u+\kappa-i)}{\la_{-n+i}(u+\kappa-i)}
=\frac{\la_{n-i}(u)}{\la_{n-i+1}(u)}
\eeq
for $i=1,\dots,n-1$ if $\agot=\oa_{2n}$ or $\spa_{2n}$,
and for $i=1,\dots,n$ if $\agot=\oa_{2n+1}$.
\epr

\bpf
Suppose first that $M(\la(u))$ is non-trivial.
We use induction on $n$ taking Corollaries~\ref{cor:fdc1},
\ref{cor:fdb1} and \ref{cor:fdd2} as the induction base.
Let us apply $t_{-n,-n+1}(u)\ts t_{n,n-1}(v)$ to the highest
vector $1_{\la}$ of $M(\la(u))$. By \eqref{defrel} we have
\ben
\bal
t_{-n,-n+1}(u)\ts t_{n,n-1}(v)&\ts 1_{\la}=\\
-\frac{1}{u-v-\kappa}\ts&\Big(t_{-n,-n+1}(u)\ts t_{n,n-1}(v)
+\la_{-n+1}(u)\ts \la_{n-1}(v)-\la_{-n}(u)\ts \la_{n}(v)\Big)\ts 1_{\la},
\eal
\een
which implies
\ben
(u-v-\kappa+1)\ts t_{-n,-n+1}(u)\ts t_{n,n-1}(v)\ts 1_{\la}
=\la_{-n}(u)\ts \la_{n}(v)\ts 1_{\la}
-\la_{-n+1}(u)\ts \la_{n-1}(v)\ts 1_{\la}.
\een
Putting $u=v+\kappa-1$ and replacing $v$ by $u$ we obtain
\eqref{nontrvm} for $i=1$. Furthermore, by Lemma~\ref{lem:vplus},
the subspace $M(\la(u))^+$ of $M(\la(u))$ is a module
over $\X(\agot')$. The highest vector $1_{\la}$ belongs
to $M(\la(u))^+$ and generates a highest weight $\X(\agot')$-module
with the highest weight $(\la_{-n+1}(u),\dots,\la_{n-1}(u))$.
So, the remaining conditions hold by the induction hypothesis.

Conversely, suppose that $\la(u)$ satisfies the conditions.
Consider the left ideal $I$ of the algebra
$\X(\agot)$ generated by
the coefficients of the series $t_{ij}(u)$ with $i<j$ where
$i+j>0$ or $i+j\geqslant 0$ for the orthogonal or symplectic case,
respectively; and by the coefficients of the series
$t_{ii}(u)-\la_i(u)$
for $i=1,\dots,n$ and $z(u)-\la_{-n}(u+\kappa)\ts\la_n(u)$.
By Corollary~\ref{cor:pbwx}, the quotient $\wt M(\la(u))=\X(\agot)/I$
is non-trivial. Let $1_{\la}$ be the image of $1\in\X(\agot)$ in the quotient.
It suffices to verify that the vector $1_{\la}$ satisfies all
the conditions \eqref{trianb}.
Now we use Corollary~\ref{cor:pbwx}
again. Let us choose the total
ordering on the elements $t_{ij}^{(r)}$ and $z_r$ with
the conditions on the indices as in the statement
of the corollary, in such a way that any element
$t_{ij}^{(r)}$ with $i>j$ precedes any element $t_{kk}^{(s)}$ while the latter
precedes any element of the form $t_{ij}^{(r)}$ with $i<j$.
We shall regard $\X(\agot)$ as the adjoint $\agot$-module with the action
defined on the generators by \eqref{tijF}.
For any pair $k<l$ and any $r\geqslant 1$
write the element $t_{kl}^{(r)}$
as a linear combination of the ordered monomials.
The $\agot$-weight of each of the monomials coincides
with the $\agot$-weight of $t_{kl}^{(r)}$.
Then the relation
$t_{kl}^{(r)}\ts 1_{\la}=0$ follows because the vector
$1_{\la}$ is annihilated by any monomial occurring in the combination.
The same argument shows that $1_{\la}$ is an eigenvector for
the action of any element $t_{kk}^{(s)}$. Thus, the
$\X(\agot)$-module $\wt M(\la(u))$ is a Verma module $M(\wt\la(u))$.
It remains to verify that its highest weight $\wt\la(u)$
coincides with $\la(u)$.
This holds for the components of $\wt\la(u)$ with positive subscripts
by the definition of $\wt M(\la(u))$. Furthermore, since
$z(u)\ts 1_{\la}=\la_{-n}(u+\kappa)\ts\la_n(u)\ts 1_{\la}$,
\eqref{zun} implies that
$t_{-n,-n}(u)\ts 1_{\la}=\la_{-n}(u)\ts 1_{\la}$. So, $\wt\la_{-n}(u)=\la_{-n}(u)$.
By the first part of the proof, since the Verma module $M(\wt\la(u))$
is non-trivial, the conditions \eqref{nontrvm} must hold
for the components of $\wt\la(u)$. This shows that $\wt\la(u)=\la(u)$,
and thus $M(\la(u))$ is non-trivial.
\epf

\bco\label{cor:existlu}
The irreducible highest weight module $L(\la(u))$
over $\X(\agot)$ exists if and only if the conditions
\eqref{nontrvm} hold.
\eco

\bpf
If $L(\la(u))$ exists then the conditions \eqref{nontrvm} are derived
by repeating the argument of the first part of the proof of
Proposition~\ref{prop:nontrvm}. Conversely, if the conditions hold
then the Verma module $M(\la(u))$ is non-trivial
by Proposition~\ref{prop:nontrvm}. Therefore, the irreducible
quotient $L(\la(u))$ of $M(\la(u))$ exists.
\epf

We are now in a position to prove the classification theorem for
finite-dimensional irreducible representations
of the extended Yangian $\X(\agot)$.

\bth\label{thm:fdim}
Every finite-dimensional irreducible $\X(\agot)$-module
is isomorphic to
$L(\la(u))$ where $\la(u)$ satisfies the conditions
\eqref{nontrvm} and there exist monic polynomials
$P_1(u),\dots,P_n(u)$ in $u$ such that
\beql{lainoli}
\frac{\la_{i-1}(u)}{\la_{i}(u)}=\frac{P_i(u+1)}{P_i(u)},\qquad
\text{for}\quad i=2,\dots,n
\eeq
and also
\ben
\frac{\la_{0}(u)}{\la_{1}(u)}=\frac{P_1(u+1/2)}{P_1(u)},
\qquad
\text{if}\quad \agot=\oa_{2n+1},
\een
\ben
\frac{\la_{-1}(u)}{\la_{1}(u)}=\frac{P_1(u+2)}{P_1(u)},
\qquad
\text{if}\quad \agot=\spa_{2n},
\een
\ben
\frac{\la_{-1}(u)}{\la_{2}(u)}=\frac{P_1(u+1)}{P_1(u)},
\qquad
\text{if}\quad \agot=\oa_{2n}.
\een
Conversely, if \eqref{nontrvm} and
the above conditions on the highest weight $\la(u)$
are satisfied then $L(\la(u))$ exists and has finite dimension.
\eth

The polynomials $P_1(u),\dots,P_n(u)$ are called the {\it Drinfeld
polynomials\/} corresponding to the finite-dimensional representation
$L(\la(u))$.

\bpf
Due to Theorem~\ref{thm:fdhw}, every finite-dimensional
irreducible $\X(\agot)$-module is isomorphic to
$L(\la(u))$ for some highest weight $\la(u)$.
Then $\la(u)$ must satisfy \eqref{nontrvm}
by Corollary~\ref{cor:existlu}
since $L(\la(u))$ exists.
Now we argue by induction on $n$ taking
Corollaries~\ref{cor:fdc1},
\ref{cor:fdb1} and \ref{cor:fdd2} as the induction base.
Observe that if $n\geqslant 2$ then by \eqref{defrel}, the mapping
\ben
T_{ij}(u)\mapsto t_{i+n-2,j+n-2}(u),\qquad i,j\in\{1,2\}
\een
defines a homomorphism $\Y(\gl_2)\to \X(\agot)$.
So, $L(\la(u))$ can be regarded as a $\Y(\gl_2)$-module.
The highest vector $1_{\la}\in L(\la(u))$ then satisfies
\ben
T_{11}(u)\ts 1_{\la}=\la_{n-1}(u)\ts 1_{\la},\qquad
T_{22}(u)\ts 1_{\la}=\la_{n}(u)\ts 1_{\la},\qquad
T_{12}(u)\ts 1_{\la}=0.
\een
Since the cyclic span $\Y(\gl_2)\ts 1_{\la}$ is finite-dimensional,
we derive from Theorem~\ref{thm:classgltwo} that there exists
a monic polynomial $P_n(u)$ such that \eqref{lainoli} holds
for $i=n$. Furthermore, by Lemma~\ref{lem:vplus}, the subspace
$L(\la(u))^+$ is a module
over $\X(\agot')$. The highest vector $1_{\la}$ belongs
to $L(\la(u))^+$ and generates a highest weight $\X(\agot')$-module
with the highest weight $(\la_{-n+1}(u),\dots,\la_{n-1}(u))$.
Since the cyclic span $\X(\agot')\ts 1_{\la}$ is finite-dimensional,
the remaining conditions on the $\la_i(u)$ hold by the induction hypothesis.

Suppose now that the highest weight $\la(u)$ satisfies the
given conditions. Then $L(\la(u))$ exists by Corollary~\ref{cor:existlu}.
We need to show that $\dim L(\la(u))<\infty$.
Observe that the $n$-tuple of Drinfeld polynomials corresponding to
an $\X(\agot)$-module $L(\la(u))$ determines the highest weight
$\la(u)$ up to a simultaneous multiplication of all components
$\la_i(u)$ by a series $f(u)\in 1+u^{-1}\CC[[u^{-1}]]$.
On the other hand, the composition of the action of
$\X(\agot)$ on $L(\la(u))$ with the automorphism \eqref{mult}
yields a representation of $\X(\agot)$ isomorphic to
$L(\la'(u))$, where the components of $\la'(u)$ are given by
$\la'_i(u)=f(u)\ts\la_i(u)$. Therefore, it suffices to prove that
a particular module $L(\la(u))$ corresponding to an arbitrary $n$-tuple
of Drinfeld polynomials is finite-dimensional.

We shall use the coproduct \eqref{Delta} to equip the tensor
product of two $\X(\agot)$-modules with an $\X(\agot)$-module structure.

\ble\label{lem:prodten}
Let $L(\la(u))$ and $L(\mu(u))$ be two irreducible
highest weight modules over $\X(\agot)$ with
\ben
\la(u)=(\la_{-n}(u),\dots,\la_n(u))\Fand
\mu(u)=(\mu_{-n}(u),\dots,\mu_n(u)).
\een
Then the tensor product $1_{\la}\ot\tss 1_{\mu}$ of the highest vectors
of $L(\la(u))$ and $L(\mu(u))$ generates a highest weight
submodule $V$ over $\X(\agot)$ in $L(\la(u))\ot L(\mu(u))$
with the highest weight
\beql{hwlamu}
(\la_{-n}(u)\ts\mu_{-n}(u),\dots,\la_n(u)\ts\mu_n(u)).
\eeq
Moreover, if the modules $L(\la(u))$ and $L(\mu(u))$ are
finite-dimensional with the corresponding $n$-tuples of Drinfeld
polynomials $(P_1(u),\dots,P_n(u))$ and $(Q_1(u),\dots,Q_n(u))$,
respectively,
then the $n$-tuple of Drinfeld polynomials corresponding to
the irreducible quotient of $V$ is
$(P_1(u)\ts Q_1(u),\dots,P_n(u)\ts Q_n(u))$.
\ele

\bpf
It follows easily from \eqref{Delta} that the vector $\xi=1_{\la}\ot\tss 1_{\mu}$
satisfies \eqref{trianb} with the highest weight \eqref{hwlamu}.
The second statement now follows from the relations defining the Drinfeld
polynomials.
\epf

By the lemma, we only need to show that if an irreducible highest weight
module $L(\la(u))$ corresponds to an $n$-tuple of Drinfeld polynomials
of the form
$P_j(u)=1$ for all $j\ne i$ and $P_i(u)=u-b$ for certain $i\in\{1,\dots,n\}$
and $b\in\CC$, then $\dim L(\la(u))<\infty$.
Furthermore, the composition of the action of
$\X(\agot)$ on $L(\la(u))$ with an automorphism of the form \eqref{shift}
yields a representation of $\X(\agot)$ whose $n$-tuple of Drinfeld
polynomials is $P_j(u)=1$ for all $j\ne i$ and $P_i(u)=u-a-b$.
Thus, it suffices to prove the claim for all values of the index
$i$ and a certain particular value of $b\in\CC$.

Consider the representation of $\X(\agot)$ on $\CC^N$ defined in
\eqref{vectrep} with $c=0$ so that
\ben
t_{ij}(u)\mapsto
\de_{ij}+e_{ij}\ts u^{-1}-\theta_{ij}\ts e_{-j,-i}\ts (u+\kappa)^{-1}.
\een
Equip the tensor product $(\CC^N)^{\ot m}$
with an $\X(\agot)$-action by
\begin{multline}\label{tijtprei}
t_{ij}(u)\ts (e_{i_1}\ot\cdots\ot e_{i_m})=\\
\sum_{a_1,\dots,a_{m-1}=-n}^n t_{ia_1}(u)\ts e_{i_1}\ot t_{a_1a_2}(u+1)\ts e_{i_2}
\ot\cdots\ot t_{a_{m-1}j}(u+m-1)\ts e_{i_m},
\end{multline}
where we use the coproduct \eqref{Delta} on $\X(\agot)$ and
the automorphism \eqref{shift}.
For any $1\leqslant m\leqslant n$ set
\ben
\xi_m=\sum_{\si\in\Sym_m} \sgn\si\cdot
e_{-n-1+\si(1)}\ot\cdots\ot e_{-n-1+\si(m)} \in
(\CC^N)^{\ot m}.
\een
We claim that $\xi_m$ satisfies
\beql{tijxi}
t_{ij}(u)\ts \xi_m=0\qquad\text{for all}\quad -n\leqslant i<j\leqslant n
\eeq
and
\beql{tiixi}
t_{ii}(u)\ts \xi_m=\begin{cases}
\dfrac{u+m}{u+m-1}\ts\xi_m \qquad&\text{if}
\quad -n\leqslant i\leqslant -n+m-1,\\[1em]
\ \xi_m\qquad&\text{if}\quad -n+m\leqslant i\leqslant n-m,\\[1em]
\dfrac{u+\kappa-1}{u+\kappa}\ts\xi_m \qquad&\text{if}
\quad n-m+1\leqslant i\leqslant n.
\end{cases}
\eeq
Denote by $P^{(m)}$ the operator in $(\CC^N)^{\ot m}$ which acts
on the basis vectors by
\ben
P^{(m)}\ts(e_{i_1}\ot\cdots\ot e_{i_m})=e_{i_m}\ot\cdots\ot e_{i_1}.
\een
We have $P^{(m)}(\xi_m)=\al\ts\xi_m$, where $\al=1$ or $-1$.
The definition \eqref{tijtprei} implies the following
relation for the action of $\X(\agot)$ on $(\CC^N)^{\ot m}$,
\beql{ttravp}
\theta_{ij}\ts t_{-j,-i}(u)=P^{(m)}\ts t_{ij}(-u-\kappa-m+1) \ts P^{(m)}.
\eeq
Due to \eqref{tijF}, in order to verify \eqref{tijxi}
in the case $\agot=\oa_{2n+1}$,
it therefore suffices to consider the
values $j=i+1$ with $-n\leqslant i\leqslant -1$.
Since the expression for the vector $\xi_m$ only involves
the tensor products $e_{i_1}\ot\cdots\ot\ts e_{i_m}$
with negative subscripts $i_k$, we may assume that
the summation indices $a_1,\dots,a_{m-1}$
in \eqref{tijtprei} are all negative. Indeed, $t_{ia_1}(u)\ts e_{i_1}=0$
unless $a_1<0$ implying $t_{a_1a_2}(u+1)\ts e_{i_2}=0$ unless $a_2<0$ etc.
However, in this case the formula \eqref{tijtprei} takes the same form
as its $\Y(\gl_N)$-counterpart \eqref{tijtpreigln} if we take into account
the convention on the basis vector indices. Therefore,
the relations $t_{i,i+1}(u)\ts \xi_m=0$ and, hence \eqref{tijxi},
are implied by the corresponding property \eqref{tijxigln}
of the vector $\xi_m$ in the case of $\Y(\gl_N)$.
Moreover, this argument also proves \eqref{tiixi} for the non-positive
values of $i$. The application of \eqref{ttravp} completes
the proof of \eqref{tiixi}.

The same argument applies to the cases $\agot=\spa_{2n}$ and $\agot=\oa_{2n}$
which also shows that
$t_{-1,1}(u)\ts \xi_m=0$ together with $t_{-1,2}(u)\ts \xi_m=0$
for $\agot=\oa_{2n}$.

Thus, in the case $\agot=\spa_{2n}$
for any $m\in\{1,\dots,n-1\}$
the vector $\xi_m$ generates a highest weight submodule
of $(\CC^N)^{\ot m}$ whose $n$-tuple of Drinfeld polynomials is
$P_j(u)=1$ for all $j\ne m$ and $P_m(u)=u+\kappa-1$,
while $\xi_n$ generates a highest weight submodule
of $(\CC^N)^{\ot n}$ whose $n$-tuple of Drinfeld polynomials is
$P_1(u)=u+n-1$ and $P_j(u)=1$ for $j\ne 1$.
This completes the proof
of the theorem in the symplectic case, as
the irreducible highest weight modules over
$\X(\agot)$ with such $n$-tuples of Drinfeld polynomials are finite-dimensional.

Similarly, the proof is also complete in the case $\agot=\oa_{2n}$
and the values $m\in\{1,\dots,n-2\}$, as well as
in the case $\agot=\oa_{2n+1}$ for the values
$m\in\{1,\dots,n-1\}$.
In order to complete the proof in the remaining cases, we
shall use the spinor
representations of the orthogonal Lie algebras.
The spinor representation $V(-1/2,\dots,-1/2)$ of the Lie algebra
$\oa_{2n+1}$ can be realized in the $2^n$-dimensional space $\La_n$ of polynomials
in $n$ anti-commuting variables
$\xi_1,\dots,\xi_n$,
\ben
\La_n=\text{span of}\ \{\xi_{i_1}\dots\xi_{i_k}\ |\
1\leqslant i_1<\cdots<i_k\leqslant n,\quad 0\leqslant k\leqslant n\}.
\een
The generators of $\oa_{2n+1}$ act
on this space as the operators
\beql{fijspin}
\bal
F_{ij}=\xi_i\ts\pa_j-\frac12\ts \de_{ij},&\qquad F_{-j,i}=\pa_i\ts\pa_j,
\qquad F_{j,-i}=\xi_i\ts\xi_j,\\
&F_{0,i}=\frac{1}{\sqrt2}\ts\pa_i,\qquad F_{i,0}=\frac{1}{\sqrt2}\ts\xi_i,
\eal
\eeq
where $i,j\in\{1,\dots,n\}$ and $\pa_i$ is the left derivative over $\xi_i$. The
restriction of $\La_n$ to the subalgebra $\oa_{2n}\subset\oa_{2n+1}$
(spanned by the elements $F_{ij}$ with $i,j\ne 0$)
splits into the direct sum of two irreducible
submodules,
$
\La_n=\La_n^+\oplus\La_n^-,
$
where $\La_n^+$ (respectively, $\La_n^-$) is the subspace of $\La_n$
spanned by the even (respectively, odd) monomials in the generators
$\xi_i$.
We have the isomorphisms
\beql{lapmisom}
\La_n^+\cong V(-1/2,\dots,-1/2)\Fand \La_n^-\cong V(1/2,-1/2,\dots,-1/2).
\eeq
The highest weight vectors of the $\oa_{2n}$-modules
$\La_n^+$ and $\La_n^-$ are, respectively, the vectors $1$ and $\xi_1$.

\ble\label{lem:spinors}
Each spinor representation of $\oa_N$ can be extended to
a representation of the algebra $\X(\oa_N)$ by the rule
\ben
t_{ij}(u)\mapsto \de_{ij}+F_{ij}\ts u^{-1}, \qquad i,j\in\{-n,\dots,n\}.
\een
\ele

\bpf
The claim follows by a direct verification that the images of $t_{ij}(u)$
satisfy the defining relations \eqref{defrel} with the use of the
following identity of operators
in each spinor representation:
\beql{fsqua}
(F^2)_{ij}=\Big(\frac{\kappa}{2}+\frac{1}{4}\Big)
\ts \de_{ij}+\kappa\ts F_{ij},
\eeq
where $F$ is defined in \eqref{F}.
Indeed, in the particular case $i=j=n$,
the identity is verified by a straightforward calculation.
The general case then follows by commuting both sides of this
particular identity with appropriate generators $F_{ij}$.
\epf

The lemma implies that the spinor representation $V(-1/2,\dots,-1/2)$
of $\oa_N$ becomes an irreducible highest weight representation
of $\X(\oa_N)$ with the highest weight $\la(u)$, where
\ben
\la_i(u)=1+\frac12 u^{-1}\quad\text{for}\quad i\leqslant -1,\qquad
\la_i(u)=1-\frac12 u^{-1}\quad\text{for}\quad i\geqslant 1
\een
and $\la_0(u)=1$ (the latter only occurs for $N=2n+1$). The corresponding
$n$-tuple of Drinfeld polynomials is
$(u-1/2,1,\dots,1)$ in both cases $N=2n$ and $N=2n+1$.
Finally, the spinor representation $V(1/2,-1/2,\dots,-1/2)$
of $\oa_{2n}$ becomes an irreducible highest weight representation
of $\X(\oa_{2n})$ with the highest weight $\la(u)$, where
\ben
\bal
\la_i(u)&=1+\frac12 u^{-1}\quad\text{for}\quad i\leqslant -2\fand i=1,\\
\la_i(u)&=1-\frac12 u^{-1}\quad\text{for}\quad i\geqslant 2\fand i=-1.
\eal
\een
The corresponding
$n$-tuple of Drinfeld polynomials is
$(1,u-1/2,1,\dots,1)$.
\epf

Theorem~\ref{thm:fdim} allows us to get another proof of Drinfeld's
classifications theorem for the Yangian modules \cite{d:nr};
cf. \cite[Chapter~12]{cp:gq}.

\bco\label{cor:yangrep}
Any finite-dimensional irreducible representation of the Yangian $\Y(\agot)$
is isomorphic to the restriction of an
$\X(\agot)$-module $L(\la(u))$ to the subalgebra $\Y(\agot)$,
where the components of $\la(u)$
satisfy the conditions of Theorem~\ref{thm:fdim}.
In particular, such representations of $\Y(\agot)$
are parameterized by the tuples $(P_1(u),\dots,P_n(u))$
of monic polynomials in $u$.
\eco

\bpf
By Theorem~\ref{thm:tenpr}, any
finite-dimensional irreducible representation $V$ of $\Y(\agot)$
can be extended to a representation of $\X(\agot)$ where
the elements of the center $\ZX(\agot)$ act as scalar operators.
By Theorem~\ref{thm:fdim}, the $\X(\agot)$-module $V$ is isomorphic
to $L(\la(u))$ for an appropriate highest weight $\la(u)$.
This allows one to attach a tuple of polynomials $(P_1(u),\dots,P_n(u))$
to the $\Y(\agot)$-module $V$.

Conversely, given any $n$-tuple of polynomials $(P_1(u),\dots,P_n(u))$,
there exists a highest weight $\la(u)$ such that the conditions of
Theorem~\ref{thm:tenpr} hold. Moreover, the components of $\la(u)$
are uniquely determined up to simultaneous multiplication by
a formal series in $u^{-1}$. This implies that the corresponding
$\X(\agot)$-module $L(\la(u))$ is determined up to twisting by an
appropriate automorphism \eqref{mult}. However, the subalgebra
$\Y(\agot)$ consists of the elements stable under
all such automorphisms. This yields the desired
parametrization of the representations of $\Y(\agot)$.
\epf

The finite-dimensional irreducible representations $L(\la(u))$
corresponding to the $n$-tuples of Drinfeld polynomials
of the form $(1,\dots,u-a,1,\dots,1)$, where
$a\in\CC$ and $u-a$ is on the $i$-th position,
are called the {\it fundamental representations\/}
of $\X(\agot)$ or $\Y(\agot)$.
The following corollary was established
in the proof of Theorem~\ref{thm:fdim}.

\bco\label{cor:subquo}
Every finite-dimensional irreducible representation of $\Y(\agot)$
is isomorphic to a subquotient of a tensor product
of the fundamental representations.
\qed
\eco

\subsection{Fundamental representations}\label{subsec:fundam}

In this section we give a more explicit description of
the fundamental representations of the algebras $\X(\agot)$ and
$\Y(\agot)$. We shall follow the general approach of the paper
by Chari and Pressley~\cite{cp:fr}. However, contrary to \cite{cp:fr},
we avoid using the theorem describing the singularities of $R$-matrices.

We start with the orthogonal case $\agot=\oa_N$.
The fundamental representations
with the $n$-tuples of Drinfeld polynomials $(u-1/2,1,\dots,1)$
and $(1,u-1/2,1,\dots,1)$ (the latter for $N=2n$ only),
were constructed in the proof of Theorem~\ref{thm:fdim}.

Now let $N=2n+1$. The tensor square of the spinor representation $\La_n$
of $\oa_{2n+1}$ has the following decomposition into irreducibles:
\beql{lalade}
\La_n\ot\La_n\cong\overset{n}{\underset{p=0}\bigoplus}\ts V(\mu^{(p)}),
\eeq
where $\mu^{(p)}=(0,\dots,0,-1,\dots,-1)$ with $p$ zeros. Note that
$V(\mu^{(p)})$ is a fundamental representation of $\oa_{2n+1}$
for any $1\leqslant p\leqslant n-1$. It
corresponds to the fundamental weight $\omega_{n-p}$
in a more standard notation.
The highest weight vector $v_p$ of $V(\mu^{(p)})$ is given in an explicit form
by
\beql{vplan}
v_p=\sum (-1)^{j_1+\cdots+j_l} \ts\xi_{i_1}\cdots\xi_{i_k}
\ot \xi_{j_1}\cdots\xi_{j_l},
\eeq
summed over all partitions of the set $\{1,\dots,p\}$ into the disjoint union of
two subsets $\{i_1,\dots,i_k\}$ and $\{j_1,\dots,j_l\}$ so that $p=k+l$
with $k,l\geqslant 0$ while
$i_1<\dots<i_k$ and $j_1<\dots<j_l$.

By Lemma~\ref{lem:spinors}, we may regard $\La_n$ as an $\X(\oa_{2n+1})$-module.
Furthermore, using the coproduct \eqref{Delta} and the automorphism \eqref{shift},
we can equip $\La_n\ot\La_n$ with an $\X(\oa_{2n+1})$-action by
\beql{tijactlala}
t_{ij}(u)(\eta\ot\zeta)=\sum_{k=-n}^n \Big(\de_{ik}+F_{ik}\ts (u-a)^{-1}\Big)\ts\eta
\ot\Big(\de_{kj}+F_{kj}\ts u^{-1}\Big)\ts\zeta,
\eeq
where $\eta,\zeta\in\La_n$ and $a\in\CC$ is a fixed constant.

\bpr\label{prop:vpsing}
If $a=p-1/2$ then the vector $v_p\in\La_n\ot\La_n$ has the properties
\beql{trianbvp}
t_{ij}(u)\ts v_p=0 \qquad \text{for}
\quad -n\leqslant i<j\leqslant n
\eeq
and
\beql{eigen}
t_{ii}(u)\ts v_p=\begin{cases}\dfrac{(u-p)(u+1/2)}{u\ts(u-p+1/2)}
\ts v_p \qquad &\text{for}
\quad 0\leqslant i\leqslant p,\\[1.5em]
\dfrac{(u-p)(u-1/2)}{u\ts(u-p+1/2)}
\ts v_p \qquad &\text{for}
\quad p+1\leqslant i\leqslant n.
\end{cases}
\eeq
\epr

\bpf
By the definition \eqref{tijactlala}, we have
\ben
t_{ij}^{(1)}\ts(\eta\ot\zeta)=F_{ij}\ts\eta \ot\zeta+\eta\ot F_{ij}\ts\zeta
\een
and
\beql{tijre}
t_{ij}^{(r)}\ts(\eta\ot\zeta)=a^{r-2}\sum_{k=-n}^n
F_{ik}\ts\eta \ot F_{kj}\zeta+a^{r-1}\ts F_{ij}\ts\eta\ot \zeta
\eeq
for $r\geqslant 2$. In particular,
\beql{tijrpone}
t_{ij}^{(r+1)}\ts(\eta\ot\zeta)=a\ts t_{ij}^{(r)}\ts(\eta\ot\zeta)
\eeq
for any $r\geqslant 2$.
Since $v_p$ is the highest weight vector
in the $\oa_{2n+1}$-module $V(\mu^{(p)})$, we have the relations
$t_{ij}^{(1)}\ts v_p=0$ for $-n\leqslant i<j\leqslant n$ and
\ben
t_{ii}^{(1)}\ts v_p=\begin{cases}0 \qquad &\text{for}
\quad 0\leqslant i\leqslant p,\\
-v_p \qquad &\text{for}
\quad p+1\leqslant i\leqslant n.
\end{cases}
\een
Now, \eqref{tijF} implies that
\ben
[F_{i-1,i},t_{i,i}^{(2)}]=t_{i-1,i}^{(2)},\qquad i=1,\dots,n.
\een
Furthermore, taking the $(i-1,i)$ entry in \eqref{cu} and comparing
the coefficients at $u^{-2}$ we get
\ben
t_{i-1,i}^{(2)}-\sum_{k=-n}^n t_{i-1,k}^{(1)}\ts t_{k,i}^{(1)}
+t_{-i,-i+1}^{(2)}-\kappa\ts t_{-i,-i+1}^{(1)}=0.
\een
Hence, \eqref{trianbvp} will follow if we prove that $v_p$ is an eigenvector
for all the operators $t_{ii}^{(2)}$ with $i=1,\dots,n$.
By \eqref{tijre}, we have the following equality of operators
in $\La_n\ot\La_n$,
\ben
t_{ii}^{(2)}=\sum_{k=-n}^n (F_{ik}\ot 1)(F_{ki}\ot 1+1\ot F_{ki})
-(F^2)_{ii}\ot 1+a\ts F_{ii}\ot 1.
\een
Note that each element $F_{ki}\in \oa_{2n+1}$ acts on $\La_n\ot\La_n$ as
the operator
\ben
\Delta(F_{ki})=F_{ki}\ot 1+1\ot F_{ki}.
\een
Due to \eqref{fsqua}, in the spinor representation $\La_n$ we have
$(F^2)_{ii}=n/2+(n-1/2)\ts F_{ii}$. Moreover, we have
$\Delta(F_{ki})\ts v_p=0$ for $k<i$ and for $1\leqslant i<k\leqslant p$.
The latter follows from the fact that each vector
$\Delta(F_{k,k-1})\ts v_p$ with $k\in\{2,\dots,p\}$ is annihilated by
all operators $\Delta(F_{j,j+1})$ and hence must be zero,
as the $\oa_{2n+1}$-module $V(\mu^{(p)})$ is irreducible.
Recalling that $a=p-1/2$ we thus get for any $i\in\{1,\dots,p\}$,
\ben
t_{ii}^{(2)}\ts v_p=\sum_{k=p+1}^n (F_{ik}\ot 1)\ts\Delta(F_{ki})\ts v_p
+(p-n)\ts (F_{ii}\ot 1)\ts v_p-n/2\ts v_p.
\een
Using the expression \eqref{vplan} for $v_p$ and the formulas
\eqref{fijspin} it is now easy to derive
the relation $t_{ii}^{(2)}\ts v_p=-p/2\cdot v_p$.
If $i\in\{p+1,\dots,n\}$ then
\ben
t_{ii}^{(2)}\ts v_p=\sum_{k=i}^n (F_{ik}\ot 1)\ts\Delta(F_{ki})\ts v_p
+(p-n)\ts (F_{ii}\ot 1)\ts v_p-n/2\ts v_p.
\een
Using again \eqref{vplan} and \eqref{fijspin}, we find that
$\Delta(F_{ki})\ts v_p=0$ for $k>i$ which gives
$t_{ii}^{(2)}\ts v_p=(-p/2+1/2)\ts v_p$. Thus, \eqref{trianbvp}
is proved. For any $i>0$ the relation \eqref{eigen}
is now implied by \eqref{tijrpone} with $j=i$. Finally,
we have $t_{00}^{(2)}\ts v_p=-p/2\cdot v_p$ which is verified
by a similar calculation. This implies \eqref{eigen} for $i=0$.
\epf

Due to Proposition~\ref{prop:vpsing},
the cyclic span $W_p=\X(\oa_{2n+1})\ts v_p$ of the highest vector
$v_p\in\La_n\ot\La_n$ is a highest weight module over $\X(\oa_{2n+1})$.
By the following theorem, $W_p$ is irreducible. This module
is finite-dimensional, and if $1\leqslant p\leqslant n-1$
then the corresponding
$n$-tuple of Drinfeld polynomials is $(1,\dots,u-1/2,1,\dots,1)$ with
$u-1/2$ on the $(p+1)$-th position;
see Theorem~\ref{thm:fdim}. So, this yields a construction of the fundamental
representations of $\X(\oa_{2n+1})$ alternative to the one used
in the proof of Theorem~\ref{thm:fdim}. The following is a version
of a result of Chari and Pressley~\cite[Theorem~6.2]{cp:fr}
and earlier results of Ogievetsky, Reshetikhin and Wiegmann~\cite{orw:pc}.
We assume that $1\leqslant p\leqslant n-1$ and $a=p-1/2$.

\bth\label{thm:fundb}
The $\X(\oa_{2n+1})$-module $W_p$ is irreducible. Its restriction
to the universal enveloping algebra $\U(\oa_{2n+1})$ is given by
\ben
W_p|_{\U(\oa_{2n+1})}\cong \overset{[(n-p)/2]}{\underset{i=0}
\bigoplus} \ts V(\mu^{(p+2i)}).
\een
\eth

\bpf
By Corollary~\ref{cor:pbwx} and Proposition~\ref{prop:vpsing}
the vector space $W_p$ is spanned by the elements
\ben
t_{j_1i_1}^{(r_1)}\dots t_{j_mi_m}^{(r_m)}\ts v_p,\qquad m\geqslant 0,
\een
with $j_a>i_a$ and $r_a\geqslant 1$. By \eqref{tijF}, the
$\oa_{2n+1}$-weights of $W_p$ have the form $\mu^{(p)}-\omega$,
where $\omega$ is a $\ZZ_+$-linear combination of the positive
roots; see their description in the beginning of Section~\ref{subsec:hw}.
However, any $\ZZ_+$-linear combination of the positive
roots has the form
$k_1\ts\ve_1+\cdots+k_n\ts\ve_n$, where the $k_i$ are integers
and the sum $k_1+\cdots+k_n$ is a non-positive integer.
Since $\mu^{(p)}-\mu^{(l)}=\ve_{l+1}+\cdots+\ve_{p}$ for $l<p$, we
conclude that, as an $\oa_{2n+1}$-module,
\beql{wpsubset}
W_p\subseteq \overset{n}{\underset{s=p}\bigoplus}\ts V(\mu^{(s)}).
\eeq
We shall now demonstrate that none of the irreducible $\oa_{2n+1}$-modules
of the form $V(\mu^{(s)})$ with $s=p+1,p+3,\dots$ can occur
in the irreducible decomposition of $W_p$. We need the following lemma
which holds for any value of the parameter $a$.

\ble\label{lem:raising}
For any $s\in\{2,\dots,n\}$ in the $\X(\oa_{2n+1})$-module
$\La_n\ot\La_n$ we have
\ben
t_{-s+1,s}^{(2)}\ts v_s=(a-s+1/2)\ts v_{s-2}.
\een
\ele

\bpf
By \eqref{tijre}, we have
\ben
t_{-s+1,s}^{(2)}=\sum_{k=-n}^n (F_{-s+1,k}\ot 1)\ts\Delta(F_{ks})
-(F^2)_{-s+1,s}\ot 1+a\ts F_{-s+1,s}\ot 1.
\een
Furthermore, \eqref{fsqua} implies
$(F^2)_{-s+1,s}=(n-1/2)\ts F_{-s+1,s}$. Moreover,
in the $\oa_{2n+1}$-submodule $V(\mu^{(s)})$ of $\La_n\ot\La_n$ we have
$\Delta(F_{ks})\ts v_s=0$ for $k\leqslant s$. Hence, applying \eqref{fijspin}
we obtain
\ben
t_{-s+1,s}^{(2)}\ts v_s=
\sum_{k=s+1}^n (\pa_k\ts\pa_s\ot 1)(\xi_k\ts\pa_s\ot 1+1\ot \xi_k\ts\pa_s)\ts v_s
+(a-n+1/2)\ts (\pa_s\ts\pa_{s-1}\ot 1)\ts v_{s}.
\een
Finally, using the formula \eqref{vplan} for $v_s$ we come to
\ben
t_{-s+1,s}^{(2)}\ts v_s=(a-s+1/2)\ts (\pa_s\ts\pa_{s-1}\ot 1)\ts v_{s}
=(a-s+1/2)\ts v_{s-2}.
\vspace{-1.5\baselineskip}
\een
\epf

Now, if the irreducible module $V(\mu^{(s)})$ with $s=p+2i-1$ for some
$i\geqslant 1$
occurs in the irreducible decomposition of $W_p$ then
$W_p$ would also contain $V(\mu^{(p-1)})$ by Lemma~\ref{lem:raising}.
But this contradicts \eqref{wpsubset}. Thus,
as an $\oa_{2n+1}$-module,
\beql{wpsubseteven}
W_p\subseteq \overset{[(n-p)/2]}{\underset{i=0}\bigoplus}\ts V(\mu^{(p+2i)}).
\eeq
We now need the following counterpart of Lemma~\ref{lem:raising}.

\ble\label{lem:lower}
Let $s\in\{2,\dots,n\}$.
If $a\ne -s+1/2$ then
the projection of
the vector $t_{s,-s+1}^{(2)}\ts v_{s-2}\in\La_n\ot\La_n$
onto the component $V(\mu^{(s)})$ in the decomposition
\eqref{lalade} is nonzero.
\ele

\bpf
Let us introduce a bilinear form on the vector space $\La_n$ by
\ben
\langle \xi_{i_1}\cdots \xi_{i_k}, \xi_{j_1}\cdots \xi_{j_l}\rangle=\de_{IJ},
\een
where $I=\{i_1,\dots, i_k\}$ and $J=\{j_1,\dots, j_l\}$ are subsets
of $\{1,\dots,n\}$ such that $i_1<\cdots < i_k$ and $j_1<\cdots < j_l$,
with $\de_{IJ}=1$ if $I=J$, and $0$ otherwise.
The form possesses the covariance property with respect to the action
of $\oa_{2n+1}$,
\ben
\langle F_{ij}\ts\eta,\ts\zeta\rangle=\langle \eta,\ts F_{ji}\ts\zeta\rangle,
\qquad \eta,\zeta\in \La_n.
\een
Extend the form $\langle\ ,\ \rangle$ to a bilinear form on
the tensor product space $\La_n\ot\La_n$ by
\ben
\langle \eta_1\ot \eta_2,\zeta_1\ot \zeta_2\rangle=
\langle \eta_1,\zeta_2\rangle
\langle \eta_2,\zeta_1\rangle.
\een
One easily verifies that this form inherits the covariance property.
In particular, the irreducible components $V(\mu^{(s)})$ in the decomposition
\eqref{lalade} are pairwise orthogonal. So the lemma will follow
if we prove that $\langle t_{s,-s+1}^{(2)}\ts v_{s-2},v_s\rangle\ne 0$.
However, a direct calculation with the use of \eqref{tijre}
shows that for any
$\eta,\zeta\in\La_n\ot\La_n$ we have
\ben
\langle\tss t_{ij}^{(2)}\tss\eta,\ts\zeta\rangle=
\langle \eta,\big(t_{ji}^{(2)}+a\ts(1\ot F_{ji}-F_{ji}\ot 1)\big) \tss\zeta\rangle.
\een
Hence, using Lemma~\ref{lem:raising} and the formulas
\eqref{fijspin} we find that
\ben
\bal
\langle t_{s,-s+1}^{(2)}\ts v_{s-2},v_s\rangle
&=\langle v_{s-2},\big(t_{-s+1,s}^{(2)}+
a\ts(1\ot F_{-s+1,s}-F_{-s+1,s}\ot 1)\big)\ts v_s\rangle\\
&=(-a-s+1/2)\ts\langle v_{s-2},\ts v_{s-2}\rangle\ne 0,
\eal
\een
completing the proof of the lemma.
\epf

If $a=p-1/2$ then the condition of Lemma~\ref{lem:lower} is satisfied for
any $s\in\{2,\dots,n\}$. Thus, Lemmas~\ref{lem:raising}
and \ref{lem:lower} imply that the $\X(\oa_{2n+1})$-module $W_p$
is irreducible and its $\oa_{2n+1}$-irreducible decomposition
coincides with the right hand side of \eqref{wpsubseteven}.
\epf

Consider now the case $\agot=\oa_{2n}$. As we mentioned in the previous
section, the restriction of the $\oa_{2n+1}$-module $\La_n$ to the subalgebra
$\oa_{2n}$ splits into the direct sum of two irreducible
submodules,
$
\La_n=\La_n^+\oplus\La_n^-,
$
and we have the isomorphisms \eqref{lapmisom}.
We have the following tensor product decompositions of the
$\oa_{2n}$-modules:
\begin{align}\label{laladed}
\La_n^+\ot\La_n^+&\cong\overset{[n/2]}{\underset{r=0}\bigoplus}\ts V(\mu^{(2r)}),\\
\La_n^+\ot\La_n^-&\cong\overset{[(n-1)/2]}{\underset{r=0}\bigoplus}
\ts V(\mu^{(2r+1)}),
\label{laladedm}
\end{align}
where $\mu^{(p)}=(0,\dots,0,-1,\dots,-1)$ with $p$ zeros. Note that
$V(\mu^{(p)})$ is a fundamental representation of $\oa_{2n}$
for any $2\leqslant p\leqslant n-1$.
The highest weight vector $v_p$ of $V(\mu^{(p)})$
in the decompositions \eqref{laladed} and \eqref{laladedm}
is given by \eqref{vplan} with the following additional restrictions:
both $k$ and $l$ are even for \eqref{laladed} with $p=2r$, while
$k$ is even and $l$ is odd for \eqref{laladedm} with $p=2r+1$.

By Lemma~\ref{lem:spinors}, we may regard $\La_n^+$ and
$\La_n^-$ as $\X(\oa_{2n})$-modules.
As in the previous case, we equip the tensor products
$\La_n^+\ot\La_n^+$ and $\La_n^+\ot\La_n^-$ with
an $\X(\oa_{2n})$-action by
\beql{tijactlalad}
t_{ij}(u)(\eta\ot\zeta)=\sum_{k=-n}^n \Big(\de_{ik}+F_{ik}\ts (u-a)^{-1}\Big)\ts\eta
\ot\Big(\de_{kj}+F_{kj}\ts u^{-1}\Big)\ts\zeta,
\eeq
where $a\in\CC$ is a fixed constant. In the following
proposition we consider the cases of even and odd $p$ simultaneously.
If $p=2r$ then $v_p\in \La_n^+\ot\La_n^+$ and if
$p=2r+1$ then $v_p\in \La_n^+\ot\La_n^-$.

\bpr\label{prop:vpsingd}
If $a=p-1$ then the vector $v_p$ has the properties
\beql{trianbvpd}
t_{ij}(u)\ts v_p=0 \qquad \text{for}
\quad -n\leqslant i<j\leqslant n
\eeq
and
\beql{eigend}
t_{ii}(u)\ts v_p=\begin{cases}\dfrac{(u-p+1/2)(u+1/2)}{u\ts(u-p+1)}
\ts v_p \qquad &\text{for}
\quad -1\leqslant i\leqslant p,\\[1.5em]
\dfrac{(u-p+1/2)(u-1/2)}{u\ts(u-p+1)}
\ts v_p \qquad &\text{for}
\quad p+1\leqslant i\leqslant n.
\end{cases}
\eeq
\epr

\bpf
The proof is essentially the same as for Proposition~\ref{prop:vpsing}
with the use of the relation \eqref{fsqua}. The calculation of
the eigenvalues of the operators $t_{ii}^{(2)}$ on $v_p$ gives
\ben
t_{ii}^{(2)}\ts v_p=\begin{cases}(1/4-p/2)
\ts v_p \qquad &\text{for}
\quad -1\leqslant i\leqslant p,\\
(3/4-p/2)
\ts v_p \qquad &\text{for}
\quad p+1\leqslant i\leqslant n.
\end{cases}
\een
These imply the desired properties.
\epf

The cyclic span $W_p=\X(\oa_{2n})\ts v_p$ of the vector
$v_p$ is a highest weight module over $\X(\oa_{2n})$.
By the following theorem, $W_p$ is irreducible. This module
is finite-dimensional, and if $2\leqslant p\leqslant n-1$
then the corresponding
$n$-tuple of Drinfeld polynomials is $(1,\dots,u-1/2,1,\dots,1)$ with
$u-1/2$ on the $(p+1)$-th position;
see Theorem~\ref{thm:fdim}. So, $W_p$ is a fundamental module over
$\X(\oa_{2n})$. The following is the $\oa_{2n}$-counterpart
of Theorem~\ref{thm:fundb}. We assume that
$2\leqslant p\leqslant n-1$ and $a=p-1$.

\bth\label{thm:fundbd}
The $\X(\oa_{2n})$-module $W_p$ is irreducible. Its restriction
to the universal enveloping algebra $\U(\oa_{2n})$ is given by
\ben
W_p|_{\U(\oa_{2n})}\cong \overset{[(n-p)/2]}{\underset{i=0}
\bigoplus} \ts V(\mu^{(p+2i)}).
\een
\eth

\bpf
Considering the $\oa_{2n}$-weights of $W_p$ and using Corollary~\ref{cor:pbwx},
we conclude that, as an $\oa_{2n}$-module,
\beql{wpsubsetevend}
W_p\subseteq \overset{[(n-p)/2]}{\underset{i=0}\bigoplus}\ts V(\mu^{(p+2i)}).
\eeq
The equality in \eqref{wpsubsetevend}
and irreducibility of the $\X(\oa_{2n})$-module
$W_p$ is implied by
the following two lemmas
which are verified in the same way as their
$\oa_{2n+1}$-counterparts.

\ble\label{lem:raisingd}
For any $s\in\{2,\dots,n\}$ in the $\X(\oa_{2n})$-module
$\La_n^+\ot\La_n^+$ or $\La_n^+\ot\La_n^-$
for even or odd $s$, respectively, we have
\ben
t_{-s+1,s}^{(2)}\ts v_s=(a-s+1)\ts v_{s-2}.
\vspace{-1.5\baselineskip}
\een
\qed
\ele

\ble\label{lem:lowerd}
Let $s\in\{2,\dots,n\}$.
If $a\ne -s+1$ then
the projection of
the vector $t_{s,-s+1}^{(2)}\ts v_{s-2}$
onto the component $V(\mu^{(s)})$ in the decomposition
\eqref{laladed} or \eqref{laladedm}, respectively, is nonzero.
\qed
\ele

In particular, if $a=p-1$ then the condition of Lemma~\ref{lem:lowerd}
is satisfied for any $s\in\{2,\dots,n\}$. This completes the
proof of the theorem.
\epf

We conclude by showing that each fundamental representation
of the Lie algebra $\spa_{2n}$ can be extended to the algebra
$\X(\spa_{2n})$ providing a fundamental representation of the latter.
Due to Theorem~\ref{thm:tenpr}, it suffices
to prove the corresponding statement for the Yangian $\Y(\spa_{2n})$.
We follow the argument of \cite{cp:fr}
adopting it to the presentation
of $\Y(\spa_{2n})$ provided by Corollary~\ref{cor:quoti}.
For any indices $k,l\in\{-n,\dots,n\}$ introduce
the elements $J_{kl}\in\Y(\spa_{2n})$ by
\ben
J_{kl}=\tau_{kl}^{(2)}-\frac12\sum_{i=-n}^n \tau_{ki}^{(1)}\tau_{il}^{(1)}.
\een
We shall identify the universal enveloping algebra $\U(\spa_{2n})$
with a subalgebra of $\Y(\spa_{2n})$ via the embedding
\eqref{emb}. Denote by $\Jc$ the subspace of $\Y(\spa_{2n})$
spanned by all elements $J_{kl}$.

\ble\label{lem:func}
The subspace $\Jc$
is stable under the adjoint action of the Lie algebra $\spa_{2n}$.
Moreover,
the $\spa_{2n}$-module $\Jc$ is isomorphic to the
adjoint representation.
\ele

\bpf
We easily derive from \eqref{tijF} that
\ben
[F_{ij},J_{kl}]=
\de_{kj}\ts J_{il}-\de_{il}\ts J_{kj}
-\de_{k,-i}\ts\theta_{ij}\ts J_{-j,l}
+\de_{l,-j}\ts\theta_{ij}\ts J_{k,-i}.
\een
This proves the first claim.
For the proof of the second,
take the coefficient at $u^{-2}$ in the relation \eqref{ztau}.
This gives
\beql{taukl}
\tau_{kl}^{(2)}+\theta_{kl}\ts\tau_{-l,-k}^{(2)}+\kappa\ts \tau_{kl}^{(1)}
-\sum_{i=-n}^n \tau_{ki}^{(1)}\tau_{il}^{(1)}=0,
\eeq
where we have used the relation $\tau_{kl}^{(1)}+
\theta_{kl}\ts\tau_{-l,-k}^{(1)}=0$. Replacing $k$ and $l$ respectively
by $-l$ and $-k$ in \eqref{taukl}, then multiplying it by $\theta_{kl}$
and adding the result to \eqref{taukl} yields $J_{kl}+\theta_{kl}\ts J_{-l,-k}=0$.
The argument is completed by observing that
$\dim\Jc=\dim\spa_{2n}$ by Corollary~\ref{cor:pbwy}.
\epf

The following lemma is straightforward from the defining relations of
$\Y(\spa_{2n})$ given
in Corollary~\ref{cor:quoti}.

\ble\label{lem:fjgen}
The algebra $\Y(\spa_{2n})$ is generated by the elements
$F_{kl}$ and $J_{kl}$ with $k,l\in\{-n,\dots,n\}$.
\qed
\ele

The fundamental representations of $\spa_{2n}$
are the modules $V(\mu^{(p)})$ where the highest weights
have the form
$\mu^{(p)}=(0,\dots,0,-1,\dots,-1)$ with $p$ zeros, for the values
$p=0,1,\dots,n-1$. In a more common notation, $V(\mu^{(p)})$
corresponds to the fundamental weight $\omega_{n-p}$.
Denote by $W_p(a)$ the fundamental representation
of $\Y(\spa_{2n})$ corresponding to the $n$-tuple of Drinfeld
polynomials $(1,\dots,u-a,1,\dots,1)$ with
$a\in\CC$ and $u-a$ on the $(p+1)$-th position.
By Theorem~\ref{thm:fdim}, the $\Y(\spa_{2n})$-module $W_p(a)$ is isomorphic
to the restriction of the $\X(\spa_{2n})$-module $L(\la(u))$
to the subalgebra $\Y(\spa_{2n})$, where the components of $\la(u)$
are given by
\ben
\la_{i}(u)=\begin{cases}
\dfrac{u-a-p}{u-a-p-1}\qquad&\text{if}
\quad -n\leqslant i\leqslant -p-1,\\[1em]
\ 1\qquad&\text{if}\quad -p\leqslant i\leqslant p,\\[1em]
\dfrac{u-a}{u-a+1} \qquad&\text{if}
\quad p+1\leqslant i\leqslant n
\end{cases}
\een
for $p=1,\dots,n-1$, and
\ben
\la_{i}(u)=\begin{cases}
\dfrac{u-a+1}{u-a}\qquad&\text{if}
\quad -n\leqslant i\leqslant -1,\\[1em]
\dfrac{u-a+1}{u-a+2} \qquad&\text{if}
\quad 1\leqslant i\leqslant n
\end{cases}
\een
for $p=0$.
So, $W_p(a)$ may also be regarded as an $\X(\spa_{2n})$-module.
Recall that the universal enveloping algebra $\U(\spa_{2n})$
is embedded into $\X(\spa_{2n})$ via \eqref{embx}.

The following is essentially a reformulation of a particular case
of \cite[Theorem~6.1]{cp:fr}.

\bth\label{thm:fundsp}
The restriction of $W_p(a)$ to $\U(\spa_{2n})$ is isomorphic
to the fundamental module $V(\mu^{(p)})$. Moreover,
the action of $\Y(\spa_{2n})$ on $V(\mu^{(p)})$ is determined
by the assignment $J_{kl}\mapsto b\ts F_{kl}$ with $b=-(n-p+1)/2+a$.
\eth

\bpf
By Theorem~\ref{thm:fdhw}, the $\X(\spa_{2n})$-module $W_p(a)$
contains a unique, up to a constant factor, highest vector $\xi$.
By the Poincar\'e--Birkhoff--Witt
theorem for $\X(\spa_{2n})$ and the relations \eqref{tijF},
$\xi$ is a unique weight vector of the weight $\mu^{(p)}$
in the $\spa_{2n}$-module $W_p(a)$. Furthermore, the irreducible
decomposition of this module takes the form
\beql{wpdec}
W_p(a)=V(\mu^{(p)})\oplus\underset{\nu}{\bigoplus}\ts c(\nu)\ts V(\nu),
\eeq
summed over the weights $\nu$ strictly preceding $\mu^{(p)}$
with respect to the standard partial ordering on the set
of $\spa_{2n}$-weights, where the $c(\nu)$
are some multiplicities. Consider the $\spa_{2n}$-module homomorphism
\beql{psijvw}
\psi:\Jc\ot V(\mu^{(p)})\to W_p(a)
\eeq
defined by
\ben
\psi:J_{kl}\ot v\mapsto J_{kl}\ts v,
\qquad v\in V(\mu^{(p)}).
\een
By Lemma~\ref{muonetwo}, the $\spa_{2n}$-module $\Jc$
is isomorphic to $V(\rho)$ with $\rho=(0,\dots,0,-2)$.
It is well known that
the irreducible decomposition of $V(\rho)\ot V(\mu^{(p)})$
contains $V(\mu^{(p)})$ with multiplicity one, and
does not contain any modules $V(\nu)$ with $\nu$ strictly
preceding $\mu^{(p)}$; see e.g. \cite{fh:rt}.
Therefore, the homomorphism $\psi$ must be
multiplication by a scalar on the component $V(\mu^{(p)})$
and zero on the other irreducible constituencies
of $V(\rho)\ot V(\mu^{(p)})$. Then by Lemma~\ref{lem:fjgen},
the subspace $V(\mu^{(p)})$ of $W_p(a)$ is stable under the action
of $\Y(\spa_{2n})$ and thus $W_p(a)=V(\mu^{(p)})$ since
$W_p$ is an irreducible $\Y(\spa_{2n})$-module.
This proves the first part of the theorem and shows that
the action of the elements $J_{kl}$ on $V(\mu^{(p)})$
is given by $J_{kl}\mapsto b\ts F_{kl}$ for some $b\in\CC$.
By Lemma~\ref{lem:fjgen}, this determines
the action of $\Y(\spa_{2n})$ on $V(\mu^{(p)})$.
Finally, the exact value of $b$ is found by calculating
the eigenvalue of the operator $J_{nn}$ on the highest vector $\xi$
of $L(\la(u))\cong W_p(a)$. This eigenvalue remains unchanged
if we multiply all components of $\la(u)$ by the formal series
$f(u)\in 1+\CC[[u^{-1}]]\ts u^{-1}$
defined from the relation
\ben
f(u)\ts f(u+\kappa)\ts \la_{-n}(u+\kappa)\ts \la_n(u)=1.
\een
In the case $1\leqslant p\leqslant n-1$ we obtain
\ben
f(u)=1+(n-p)\ts u^{-2}+\cdots.
\een
By Proposition~\ref{prop:actzu}, we have $z(u)=1$ in the
$\X(\spa_{2n})$-module $L\big(f(u)\la(u)\big)$ so that the eigenvalue
of $\tau_{nn}(u)$
on the highest vector of $L\big(f(u)\la(u)\big)$ is $f(u)\ts\la_n(u)$.
This allows one to find the eigenvalue of $\tau_{nn}^{(2)}$
which turns out to be $(n-p)/2-a+1$.
Since the eigenvalue of $\tau_{nn}^{(1)}=F_{nn}$
on the highest vector is $-1$, the eigenvalue of
$J_{nn}$ is $(n-p+1)/2-a$ proving the claim for the case under
consideration. In the case $p=0$ the value of $b$ is found by
the same calculation.
\epf

\end{document}